\documentclass[eqno,12pt]{article}

\usepackage{amssymb}
\usepackage{latexsym}
\usepackage{euscript}
\usepackage{mdwlist}
\usepackage{flafter}
\usepackage{pstricks}
\usepackage{graphicx}
\usepackage{epstopdf} %%Pour inclure des fichiers .eps%%
\usepackage[hidelinks,breaklinks]{hyperref} % Pour les liens hypertextes ; l'option hidelinks permet de ne pas avoir d'encadrements, breaklinks permet de ne pas avoir les reférences qui dépssent sur arxiv. 
\usepackage[latin1]{inputenc}
\usepackage{amsmath}
\usepackage{mathrsfs} %%% Nice caligraphic fonts
\evensidemargin\oddsidemargin
\usepackage{pgf,tikz} % Pour les figures
\usetikzlibrary{shapes} % pour les figures
\usetikzlibrary{arrows} % pour les figures
\usepackage{float} % Pour forcer les figures à se placer là où on leur dit, mettre [H] en argument : \begin{figure}[H]

\begin{document}
\baselineskip=18pt
\setcounter{page}{1}

\renewcommand{\theequation}{\thesection.\arabic{equation}}
\newtheorem{theorem}{Theorem}[section]
\newtheorem{lemma}[theorem]{Lemma}
\newtheorem{definition}[theorem]{Definition}
\newtheorem{proposition}[theorem]{Proposition}
\newtheorem{corollary}[theorem]{Corollary}
\newtheorem{remark}[theorem]{Remark}
\newtheorem{fact}[theorem]{Fact}
\newtheorem{problem}[theorem]{Problem}
\newtheorem{question}[theorem]{Question}
\newtheorem{conjecture}[theorem]{Conjecture}
\newtheorem{claim}[theorem]{Claim}

%%%%%%%%%%%%%%%%%%%%%%%%%%% Equation numberings
\newcommand{\eqnsection}{
\renewcommand{\theequation}{\thesection.\arabic{equation}}
    \makeatletter
    \csname  @addtoreset\endcsname{equation}{section}
    \makeatother}
\eqnsection
%%%%%%%%%%%%%%%%%%%%%%%%%%%

 %%%%%%%%%%%%%% Bbb characters
%%%%%%%%%%%%%% Real numbers
\def\r{{\mathbb R}}
%%%%%%%%%%%%%% Expectation
\def\e{{\mathbb E}}
%%%%%%%%%%%%%% Probability
\def\p{{\mathbb P}}
\def\q{{\mathbb Q}}
%%%%%%%%%%%%%% Law of environment
\def\P{{\bf P}}
\def\E{{\bf E}}
\def\Q{{\bf Q}}
%%%%%%%%%%%%%% Integers
\def\z{{\mathbb Z}}
%%%%%%%%%%%%%% Natural numbers
\def\N{{\mathbb N}}
%%%%%%%%%%%%%% Tree
\def\T{{\mathbb T}}
%%%%%%%%%%%%%% Galton-Watson tree
\def\G{{\mathbb G}}
%%%%%%%%%%%%%% Stopping lines
\def\L{{ \mathscr L}}
%%%%%%%%%%%%%
\def\C{{\mathscr C}}
%%%%%%%%%%%Type 1%%%%
\def\ZZ{{\mathscr Z}}

%%%%%%%%%%%%%%%% Special symbols
%%%%%%%%%%%%%% Exponential
\def\ee{\mathrm{e}}
%%%%%%%%%%%%%% Differentiation
\def\d{\, \mathrm{d}}
%%%%%%%%%%%%%Spine
\def\w{{\tt w}}
%%%%%%%%%%%%%%
\def\law{{\buildrel \mbox{\rm\tiny (law)} \over =}}
\def\wcv{{\buildrel \mbox{\rm\tiny  (law)} \over \longrightarrow}}
\newcommand{\parent}[1]{{\buildrel \leftarrow \over {#1}}} % Flèche retournée au dessus de 

%%%%%%%%%%%%%% Beginning of the text

\vglue30pt

\centerline{\Large\bf  Favorite sites of randomly biased walks } 
 
\centerline{\Large \bf  on a supercritical Galton--Watson tree}

{
\let\thefootnote\relax\footnotetext{\scriptsize Cooperation between  D.C. and Y.H.  was  supported by NSFC 11528101.}
}

\bigskip
\bigskip

 \centerline{by}

\medskip

 \centerline{Dayue Chen\let\thefootnote\relax\footnote{\scriptsize School of Mathematical Sciences, Peking University
Beijing 100871 China, {\tt dayue@math.pku.edu.cn}}, Lo\"{i}c de Raph\'elis\footnote{\scriptsize   UMPA, ENS de Lyon, 46 allée d'Italie,  69364 Lyon Cedex 07, France, {\tt loic.de-raphelis@ens-lyon.fr}},  and Yueyun Hu\let\thefootnote\relax\footnote{\scriptsize LAGA, Universit\'e Paris XIII, 99 avenue J-B Cl\'ement, F-93430 Villetaneuse, France, {\tt yueyun@math.univ-paris13.fr}}}

\medskip

  \centerline{\it Peking University, ENS de Lyon,  and Universit\'e Paris XIII}

\bigskip
 \centerline{November 14, 2016 }
\bigskip
\bigskip

{\leftskip=2truecm \rightskip=2truecm \baselineskip=15pt \small

\noindent{\slshape\bfseries Summary.}  Erd\H os and R\'ev\'esz~\cite{erdos-revesz}  initiated the study of favorite sites by considering the one-dimensional simple random walk.  We  investigate   in this paper the same problem for a class of null-recurrent randomly biased walks on a supercritical Gaton-Watson tree. We prove that there is some parameter  $\kappa \in (1, \infty]$ such that  the set of the favorite sites of the biased walk    is almost surely bounded in the case $\kappa \in (2, \infty]$,  tight in the case $\kappa=2$, and  oscillates between a   neighborhood of the root and the boundary of the range     in the case  $\kappa \in (1, 2)$. Moreover, our results yield a complete answer to the cardinality of the set of favorite sites in the case  $\kappa \in (2, \infty]$.  The proof relies on the exploration of the Markov property of the local times process with respect to the space variable and on a precise tail estimate on the maximum of local times, using a change of measure for  multi-type Galton-Watson trees.

%%\noindent{\slshape\bfseries Summary.} Consider a class of null-recurrent randomly biased walks on a supercritical Gaton-Watson tree. We are interested in the  subdiffusive regime where the asymptotic behaviors of the biased walks are determined by  some parameter $\kappa \in (1, \infty]$. We study the set of the favorite sites of the walks and prove that it is almost surely bounded in the case $\kappa \in (2, \infty]$,  tight in the case $\kappa=2$, and  oscillates between a   neighborhood of the root and the boundary of the range     in the case  $\kappa \in (1, 2)$.   The proof relies on the exploration of the Markov property of the local times process with respect to the space variable and on a precise tail estimate on the maximum of local times. 

\bigskip

\bigskip


} %%%%%% End of narrower
\newpage

\section{Introduction}
  \label{s:intro}

%$\phantom{aob}$

 The study of favorite sites goes back  to Erd\H os and R\'ev\'esz~\cite{erdos-revesz} where they considered the simple    random walk on $\z$,  and    conjectured that 
 
 (a)  the set of favorite sites is tight; 
 
 (b) the cardinality of the set of favorite sites is eventually bounded by $2$.

 We refer to  R\'ev\'esz (\cite{revesz}, Chapter 13) for    a list of ten open problems on the favorite sites.

 Conjecture (b) still remains  open and the best result so far was obtained by T\'oth~\cite{toth}. Conjecture (a) was disproved by Bass and Griffin~\cite{bass-griffin} who showed  the almost sure transience of the favorite sites  for the simple   random walk on $\z$ as well as for the one-dimensional Brownian motion. We note in passing that the exact rate of escape of the favorite sites  is still an open problem.  Later, the transience of the favorite sites was established by Bass,  Eisenbaum and Shi~\cite{BES00}, Marcus~\cite{marcus}, Eisenbaum and Khoshnevisan~\cite{EK02} for  L\'evy processes and even for fairly general Markov processes, and by Hu and Shi~\cite{fav-rwre} for Sinai's one-dimensional random walk in random environment. One may wonder whether the favorite sites are always transient for general ``non-trivial" null-recurrent Markov processes. This was however disproved by Hu and Shi~\cite{yzfavtree}: the set of the favorite sites is tight for a class of randomly biased walks on trees in the slow-movement regime. The present paper is to address the same question in the sub-diffusive regime. The answer is more complicated and is depending on some parameter $\kappa\in (1, \infty]$. For a class of biased walk on tree, conditioned upon the set of non-extinction of the tree, the set of favorite sites will be almost surely bounded if $\kappa > 2$, tight if $\kappa =2$, and may move to infinity almost surely if $1< \kappa < 2$.  As a consequence of our results, we can give a complete answer to the cardinality of the set of favorite sites when $\kappa >2$. 
 
 %We may wonder whether for a general ``non-trivial'' null-recurrent Markov process, the favorite sites are always transient: this was disproved by Hu and Shi~\cite{yzfavtree} where  for a class of  randomly biased walks on trees  in the slow-movement regime, the favorite sites are shown to be tight.  The present paper aims at   the study of the favorite sites of  randomly biased walks on trees in the sub-diffusive regime: for a class of biased walks on trees depending on some parameter $\kappa\in (1, \infty]$,  conditioned upon  the set of non-extinction of the  tree, the favorite sites will be  almost surely bounded if $\kappa>2$, tight if $\kappa=2$, and may escape to infinity almost surely  if $1< \kappa<2$. 
 
Let us define now the model of  the randomly biased walk on trees, a model  introduced by Lyons and Pemantle~\cite{lyons-pemantle}. Let $\T$ be a supercritical Galton--Watson tree, rooted at $\varnothing$. For any vertex $x\in \T \backslash \{ \varnothing\}$, let ${\buildrel \leftarrow \over x}$ be its parent.  Let $\omega := (\omega(x, \cdot), \, x\in \T)$ be a sequence of vectors such that for each vertex $x\in \T$,  $\omega(x, \, y) \ge 0$ for all $y\in \T$ and  $\sum_{y\in \T} \omega(x, \, y) =1$. We assume that $\omega(x, \, y)>0$ if and only if either ${\buildrel \leftarrow \over x}=y$ or ${\buildrel \leftarrow \over y}=x$. Denote by $|x|$ the generation of the vertex $x$ in $\T$. We shall also  use the partial order on the tree: for any $x, y \in \T$, we write $y< x$ if and only if $y$ is an ancestor of $x$ (and $y \le x$ iff $y < x$ or $y=x$).

 For the sake of presentation, we add a specific vertex ${\buildrel \leftarrow \over \varnothing}$, considered as the parent of $\varnothing$. We stress that ${\buildrel \leftarrow \over\varnothing}$ is not a vertex of $\T$, for instance, $\sum_{x\in \T} f(x)$ does not contain the term $f({\buildrel \leftarrow \over\varnothing})$. We define $\omega({\buildrel \leftarrow \over\varnothing}, \varnothing):=1$ and modify the vector $\omega(\varnothing, \cdot)$ such that $\omega(\varnothing, {\buildrel \leftarrow \over \varnothing})>0$ and $\omega(\varnothing, {\buildrel \leftarrow \over \varnothing})+ \sum_{x: {\buildrel \leftarrow \over x}=\varnothing} \omega(\varnothing, x)=1$.

For given $\omega$, the randomly biased walk $(X_n)_{n\ge 0}$ is a Markov chain on $\T\cup\{{\buildrel \leftarrow \over\varnothing}\}$ with transition probabilities $\omega$,  starting from $ \varnothing$; i.e.\ $X_0= \varnothing$ and 
$$
P_\omega \big( X_{n+1} = y \, | \, X_n =x \big) 
= 
\omega(x, \, y).
$$

%%$$A(x):= \frac{\omega({\buildrel\leftarrow \over x}, x)}{\omega({\buildrel\leftarrow \over x}, {\buildrel\Leftarrow \over x})}.$$

For any vertex $x\in \T $, let $(x^{(1)}, \cdots, x^{(\nu_x)})$ be its children, where $\nu_x \ge 0$ is the number of children of $x$. Define ${\bf A}(x) := (A(x^{(i)}), \, 1\le i\le \nu_x)$ by
\begin{equation*}
  A(x^{(i)}) 
  := 
  \frac{\omega(x, \, x^{(i)})}{\omega(x, \, {\buildrel\leftarrow \over x})},
  \qquad 1\le i\le \nu_x \, .
%%  \label{A}
\end{equation*}

%When all $A(x^{(i)})=\lambda$ with some positive constant $\lambda$, the walk is called $\lambda$-biased walk on a Galton-Watson tree whose studies are initiated by  Lyons, Pemantle and Peres~\cite{lyons-pemantle-peres96}. We  refer to Aidekon~\cite{elie-vitesse} for a formula on the speed of the $\lambda$-biased walk in the transient case and for recent references on the $\lambda$-biased walk. 

We denote the vector ${\bf A}(\varnothing)$ by $(A_1, ..., A_{\nu})$. As such, $\nu\equiv \nu_\varnothing$ is the number of children of $\varnothing$. 
When $\nu$ is an integer  (i.e. $\T$ is a regular tree), we  suppose  that $({\bf A}(x))_{ x\in \T}$ are i.i.d.  In general, when $\nu$ is also random, we may construct a marked tree as  in  Neveu~\cite{neveu} such that for any $k\ge 0$, conditionally on $\{{\bf A}(x), |x| \le k\}$, the random variables  $({\bf A}(y))_{|y|=k+1}$ are i.i.d. and distributed as ${\bf A}(\varnothing)$. There is an obvious bijection between $({\bf A}(x))_{ x\in \T}$ and $(\T, \omega)$ and we shall both   notation interchangeably. 

 %The biased walk $(X_n)_{n\ge 0}$ is a  random walk in random environment $(\T, \omega)$. 
 Denote by $\P$ the law of $(\T, \omega)$ and define $\p(\cdot):=\int P_\omega(\cdot) \P(d\omega)$. In the language of random walk in random environment, $P_\omega$ is referred to the quenched probability whereas $\p$ is the annealed probability.

  Assume that $\P(\nu=\infty)=0$, $\E(\nu)\in (1, \infty]$, $\E \big( \sum_{i=1}^\nu \, A_i |\log A_i| \big) <\infty$ and \begin{equation}
  \E \Big( \sum_{i=1}^\nu \, A_i \Big) =1, \qquad \E \Big( \sum_{i=1}^\nu \, A_i \log A_i \Big) <0. \label{hyp1}
  \end{equation}

\noindent  We suppose that     \begin{equation}\label{kappa}
\begin{cases}
\mbox{either there exists a $\kappa\in (1, \infty)$ such that } 
 \E \big( \sum_{i=1}^\nu \, A_i ^\kappa \big) =1 ; 
 \\ 
\mbox{or }   \E \big( \sum_{i=1}^\nu \, A_i ^t \big) < 1 \mbox{ for any $t>1$},  \mbox{ and define   $\kappa:=\infty$ in this case. }
\end{cases}
\end{equation}

\noindent 
For the sake of presentation, we suppose that  \begin{equation}\label{non-lattice}
 \hbox{the support of $\sum_{i=1}^\nu \delta_{\{\log A_i\}} $ is non-lattice when $1< \kappa \le 2$.}
 \end{equation}

\noindent We furthermore assume an integrability condition which is slightly stronger than the usual $X^\kappa \log X$-type condition as in Liu~\cite{liu00}: when $1< \kappa < \infty$, there exists some $\alpha> \kappa $ such that \begin{equation}\label{hyp2} 
 \E \big( \sum_{i=1}^\nu \, A_i\big)^\alpha < \infty, 
\end{equation}

\noindent and when $\kappa=\infty$, we assume that~\eqref{hyp2} holds for some $\alpha>2$. 
% \begin{equation}\label{hyp2} 
% 	\begin{cases}  \E \big( \sum_{i=1}^\nu \, A_i\big)^\kappa + \E \big( \sum_{i=1}^\nu \, A_i ^\kappa \, (\log A_i)^+ \big)< \infty,  \quad & \mbox{ if } 2 < \kappa < \infty, \\
%	\mbox{there exists some $\alpha> \kappa $ such that }  \E \big( \sum_{i=1}^\nu \, A_i\big)^\alpha < \infty,  \quad & \mbox{ if } 1 < \kappa \le 2.
%	\end{cases}
%	\end{equation}

 It is known from Lyons and Pemantle (\cite{lyons-pemantle}), Menshikov and Petritis~\cite{menshikov-petritis} and Faraud~\cite{faraud}  that under~\eqref{hyp1}, $(X_n)_{n\ge 0}$ is null-recurrent. When~\eqref{hyp1} and~\eqref{kappa} are fulfilled,  $(X_n)_{n\ge 0}$ may be   diffusive or subdiffusive. For instance, we have proved in~\cite{yzptrf} that if furthermore $\nu$ equals some integer, then almost surely,  \begin{equation} \label{maxlimsup} \lim_{n\to\infty} \frac1{\log n} \log \max_{0\le i\le n} |X_i| = 1 - \max( {1\over 2}, { 1 \over \kappa}) . \end{equation}

\noindent
 
When $\kappa$ is sufficiently large (say $\kappa \in (5, \infty]$), Faraud~\cite{faraud} proved an invariance principle for  $(|X_n|)_{n\ge 0}$, in line of Peres and Zeitouni~\cite{peres-zeitouni}. Recently,  A\"\i d\'ekon and de Raph\'elis~\cite{AR15} proved that for any $\kappa \in (2, \infty]$, the tree visited by the walk, after renormalization, converges to the Brownian forest. When $1 < \kappa \le 2$, a similar convergence also holds, but towards the stable forest, and the height function of the walk also satisfies a central limit theorem, see~\cite{loic16}.  We refer to Andreoletti and Debs~\cite{andreoletti-debs1, andreoletti-debs2}, Andreoletti and Chen~\cite{andreoletti-chen} for the recent studies of the spread and local times of the biased walk in both subdiffusive and slow-movement regimes. For further detailed references and open problems, see  the survey paper by   Ben Arous and Fribergh~\cite{bf14}. %  on  the biased walks on a random graph.

In this paper, we are interested in the favorite sites of the walk. Let $$L_n(x):= \sum_{i=1}^n 1_{\{ X_i=x\}}, \qquad x \in \T, \, n\ge1, $$ be the local times process of $(X_n)_{n\ge 0}$. The set of the favorite sites is defined as follows: $${\mathbb F}(n) 
:=
 \Big\{ x\in \T: L_n(x)= \max_{y\in \T} L_n(y) \Big\}.$$

Denote by $\P^*$ the probability $\P$ conditioned on the  non-extinction of the Galton-Watson tree $\T$: $$
\P^*(\bullet):=\P\Big(\bullet \, |\,  \mbox{$\T$ is infinite} \Big), 
$$
and denote by $\p^*$  the (annealed) probability conditioned on the set  of non-extinction of  $\T$:  $\p^*(\cdot):=\int P_\omega(\cdot) \P^*(d\omega)$.

 The main  result of this paper is the following description of the favorite sites in the (sub)diffusive regime.

\begin{theorem}\label{t:main} Assume~\eqref{hyp1},~\eqref{kappa},~\eqref{non-lattice} and~\eqref{hyp2}. 

(i) If $\kappa\in (2, \infty]$, then  $\p^*$-almost surely,   \begin{equation}\label{t:kappa>2}  {\mathbb F}(n)  \, \subset\,{\mathscr M}, \qquad \mbox{for all large $n$}, \end{equation} 

\noindent where ${\mathscr M}$ defined in \eqref{MinU}, is almost surely finite. 
% \begin{equation}\label{t:kappa>2} \limsup_{n\to \infty} \sup_{x\in {\mathbb F}(n)} | x| < \infty, \qquad \p^*\mbox{-a.s.} \end{equation} 

(ii) If $\kappa=2$, then $(\sup_{x\in\mathbb{F}(n)} |x|)_{n\ge 1}$ is tight; more precisely, 
\begin{equation}\label{t:kappa=2}
\p\Big(\mathbb{F}(n)\subset \mathscr{M}\Big)\substack{ \\ \longrightarrow \\ {n\to\infty}} 1. 
\end{equation}

(iii) If $1< \kappa<2$, then for any $\varepsilon>0$,   \begin{eqnarray}\label{t:kappa<2}
 \limsup_{n\to\infty} n^{- \frac{\kappa-1}{\kappa} + \varepsilon} \inf_{x\in {\mathbb F}(n)} | x| &=&\infty, \qquad \p^*\mbox{-a.s.},
\\ \nonumber\\ 
 \liminf_{n\to\infty} \sup_{x\in {\mathbb F}(n)} | x| &<&\infty, \qquad \p^*\mbox{-a.s.} \label{t:kappa<2b}
 \end{eqnarray}
 \end{theorem}

  \begin{remark} The assumption \eqref{non-lattice}  is used in the papers \cite{afanasyev, loic16}, see the proof of Lemma \ref{l:maxYexcursion} and Fact \ref{f:loic}.  Should there be a lattice version of their estimates, then   \eqref{non-lattice} could be removed.
    \end{remark}
  
 As mentioned before,  in the slow-movement regime (which corresponds informally to $\kappa=1$), we have proved in~\cite{yzfavtree} that  the set of  favorite sites is tight. %%, but we do not know its almost sure behaviors. 
  
  In the case $1< \kappa <2$,  $\p^*$-a.s., $ \max_{0\le i \le  n}|X_i| = n^{\frac{\kappa-1}{\kappa}+ o(1)}$ (see~\eqref{maxlimsup}), so~\eqref{t:kappa<2} says that up to $n^{o(1)}$, the favorite sites could   reach as far as the upper limits of the walk itself.

It is worthy noticing    the phase transition at $\kappa=2$. The almost sure oscillation  in the case $1 < \kappa < 2$ seems rather surprising, because a priori, we cannot expect a localization of a null-recurrent walk on the tree, in contrast with the one-dimensional random walk in random environment on $\z$. 
%%Moreover, when $1 < \kappa< 2$, the set of the favorite sites $ {\mathbb F}(n) $ is not tight   under $\p^*$. \\

 Let us present now two consequences of  Theorem \ref{t:main} on the structure of the set of favorite sites. The first  one  deals with the cardinality of the set of favorite sites: we show that when    $\kappa \in (2, \infty]$, it is eventually bounded by $3$ and this bound is   optimal.

\begin{corollary}\label{cor:cardinal}
If $\kappa \in (2, \infty]$, then $\p$-almost surely, 
\begin{equation}
\#\mathbb{F}(n)\le 3 \textrm{ for all $n$ large enough.} \label{card3}
\end{equation}
 Moreover,   almost surely there will be infinitely many $n$ such that $\#\mathbb{F}(n)=\min(3, \#\mathscr{M})$.
\end{corollary}

The second corollary studies the vertices which are favorite sites infinitely many times.    We show that for any $\kappa \in (1, \infty]$, such vertices do exist and must be in  $\mathscr{M}$. Moreover,   if $\kappa\in [2, \infty]$, then all vertices in $\mathscr{M}$ will be favorite sites infinitely many times.

\begin{corollary}\label{cor:limsup}
If $\kappa \in (1, \infty]$, then $\p$-almost surely
\begin{equation*}
\limsup_{n\to\infty} \mathbb{F}(n) \subset \mathscr{M}. 
\end{equation*}
Moreover, 
\begin{itemize}
\item[$\bullet$] If $\kappa\in(1, 2)$, then a.s.\ there exists $x\in\mathscr{M}$ such that $x\in\limsup_{n\to\infty} \mathbb{F}(n)$.
\item[$\bullet$] If $\kappa \in [2, \infty]$, then $\limsup_{n\to\infty} \mathbb{F}(n)=\mathscr{M}$ a.s.
\end{itemize}
\end{corollary}

Let us now describe the strategy of the  proof of Theorem~\ref{t:main}. The main step will be the exploration of the  Markov property on the space variable $x$ of the local times process.  It will be more convenient to consider  the edge local times $(\overline L_n(x))_{n\ge 1, x \in \T}$ defined as follows:   \begin{equation}\label{edgelocaltime} \overline L_n(x):= \sum_{i=1}^n 1_{\{ X_{i-1}={\buildrel\leftarrow \over x}, X_i=x\}}, \qquad x\in \T , \, n\ge 1. \end{equation}

% \noindent Recall that $\overline T_1:= \inf\{n\ge 1: X_{n-1}={\buildrel \leftarrow \over \varnothing}, X_n=\varnothing\}$. 
  
 \noindent  
We define a sequence of stopping times $( \overline T_n)_{n\ge1}$  by induction: for any $n\ge1$, \begin{equation}\label{overtj}
\overline T_n:= \inf\{k> \overline T_{n-1}: X_{k-1}= {\buildrel \leftarrow \over \varnothing}, \, X_k=\varnothing\},\end{equation}

\noindent with $\overline T_0:=0$. By definition, $\overline T_n - 1$ is exactly the $n$\textsuperscript{th} return time  to ${\buildrel \leftarrow \over \varnothing}$ of the walk $(X_n)_{n\ge 0}$. %, and $L_{\overline T_n}(\varnothing) \ge n$. 

The key ingredient in the proof of  Theorem~\ref{t:main} (cases $1 < \kappa \le 2$) is the tail distribution of the maximum of  (edge) local times considered at $\overline T_1$. 

Denote by $f(x) \sim g(x)$ as $x\to x_0$ when $\lim_{x\to x_0} f(x)/g(x) =1$ and $f(x) \asymp g(x)$ if $0< \liminf_{x\to x_0} f(x)/g(x) \le \limsup_{x\to x_0} f(x)/g(x) < \infty$.

 \begin{theorem}\label{t:tailmaxlocaltime}  Let $\kappa\in (1, \infty)$. Assume~\eqref{hyp1},~\eqref{kappa},~\eqref{non-lattice} and~\eqref{hyp2}.   As $r \to \infty$, we have 
 $$
 \p\Big( \max_{x\in \T} \overline L_{\overline T_1}(x) \ge r \Big) 
 \,\asymp\, 
 \begin{cases}
 r^{-1} , \qquad &\mbox{if } 1< \kappa<2,\\
 r^{-1} (\log r)^{-1/2}, \qquad &\mbox{if } \kappa=2,\\
 r^{-\kappa/2}, \qquad &\mbox{if } 2< \kappa<\infty.\\
 \end{cases}
 $$
 The same results hold when we replace $ \max_{x\in \T} \overline L_{\overline T_1}(x)$ by $ \max_{x\in \T}  L_{\overline T_1}(x)$. 
 \end{theorem}

\begin{remark} We mention that in the case $1< \kappa< 2$, $ \max_{x\in \T} \overline L_{\overline T_1}(x)$ has a Cauchy-type tail (independent of the value of $\kappa$), such phenomenon is in the same flavor as in Bertoin \cite{bertoin11}, Corollary 1. 
\end{remark}

To see how the asymptotic behaviors of the favorite sites ensue from Theorem~\ref{t:tailmaxlocaltime}, we introduce a set $\ZZ_k$ of vertices in $\T$: \begin{equation}\label{zzk}
\ZZ_k
:=
\Big\{x \in \T: \overline L_{\overline T_k}(x) =1, \min_{\varnothing < y < x } \overline L_{\overline T_k}(y) \ge 2 \Big\}, \qquad k\ge1, \end{equation}

\noindent 
which is the set of the vertices being the first of their ancestry line to be of edge local time $1$. This set is represented on Figure~\ref{f:Zk}. For any fixed $k\ge 1$, by the strong Markov property, we get the following identity in law under the annealed probability measure $\p$:
\begin{equation}\label{lawmax}
  \max_{x\in \T} \overline L_{\overline T_k}(x) 
 \, \law\, 
  \max\Big( \max_{x\le \ZZ_k} \overline L_{\overline T_k}(x) , \, \max_{1\le i \le \#\ZZ_k} \overline L^{*, i}\Big), \end{equation}
  
 \noindent where $x\le \ZZ_k$ means that either $x \in \ZZ_k$ or for any $\varnothing < y \le x, \overline L_{\overline T_k}(y) \ge 2$,  and $(\overline L^{*, i})_{i\ge1}$ are i.i.d.\ copies of $ \max_{x\in \T} \overline L_{\overline T_1}(x)$, independent of $( \max_{x\le \ZZ_k} \overline L_{\overline T_k}(x) , \#\ZZ_k)$.  A similar identity in law holds for the site local times $ L_{\overline T_k}(x)$ instead of the edge local time $\overline L_{\overline T_k}(x)$.

Consider for instance the case $1< \kappa <2$. It is known (see the forthcoming Fact~\ref{f:loic}) that $\#\ZZ_k$ is of order $k$ when $k\to \infty$, and it is not very hard to see that $\max_{x\le \ZZ_k} \overline L_{\overline T_k}(x)$ is also of order $k$. By Theorem~\ref{t:tailmaxlocaltime}, $\overline L^{*, i}$ has the Cauchy-type tail, then an application of the extreme value theory based on~\eqref{lawmax} yields that along some subsequence $k\to \infty$, $\frac1k \max_{x\in \T} \overline L_{\overline T_k}(x) \to \infty$ (Propositon~\ref{p:kappa<2}). This implies  the almost sure unboundedness of the favorite sites.  On the other hand,  the favorite sites are either bounded or escape to infinity at a certain polynomial rate (see~\eqref{favempty}). Combining the two facts we get the upper limits in~\eqref{t:kappa<2}. The lower limits in~\eqref{t:kappa<2b} can be obtained in a similar way. 

\begin{figure}[H]
\begin{tikzpicture}[line cap=round,line join=round,>=triangle 45,x=2.5cm,y=2.5cm]
\clip(2.2,-0.1) rectangle (9.0,2.5);
\draw (5.0,0.0)-- (5.0,0.4);
\draw (5.0,0.4)-- (4.6,0.8);
\draw (5.0,0.4)-- (5.0,0.8);
\draw (5.0,0.4)-- (5.4,0.8);
\draw (4.6,0.8)-- (4.,1.2);
\draw (4.6,0.8)-- (4.6,1.2);
\draw (5.0,0.8)-- (5.0,1.2);
\draw (5.0,0.8)-- (5.4,1.2);
\draw (5.4,0.8)-- (5.8,1.2);
\draw (4.,1.2)-- (3.2,1.6);
\draw (4.,1.2)-- (4.0,1.6);
\draw (4.6,1.2)-- (4.2,1.6);
\draw (4.6,1.2)-- (4.6,1.6);
\draw (5.0,1.2)-- (5.2,1.6);
\draw (5.4,1.2)-- (5.4,1.6);
\draw (5.0,1.2)-- (5.,1.6);
\draw (5.8,1.2)-- (6.2,1.6);
\draw (5.8,1.2)-- (5.8,1.6);
\draw (3.2,1.6)-- (2.8,2.0);
\draw (3.2,1.6)-- (3.2,2.0);
\draw (3.2,1.6)-- (3.6,2.0);
\draw (4.2,1.6)-- (3.8,2.0);
\draw (4.2,1.6)-- (4.2,2.0);
\draw (5.,1.6)-- (4.8,2.0);
\draw (5.2,1.6)-- (5.2,2.0);
\draw (5.8,1.6)-- (5.6,2.0);
\draw (5.8,1.6)-- (6.0,2.0);
\draw (6.2,1.6)-- (6.6,2.0);
\draw (6.2,1.6)-- (6.2,2.0);
\draw (3.8,2.0)-- (3.8,2.4);
\draw (3.8,2.0)-- (4.2,2.4);
\draw (4.8,2.0)-- (4.6,2.4);
\draw (4.8,2.0)-- (5.0,2.4);
\draw (5.6,2.0)-- (5.4,2.4);
\draw (5.6,2.0)-- (5.6,2.4);
\draw (5.6,2.0)-- (5.8,2.4);
\draw (6.2,2.0)-- (6.2,2.4);
\draw (6.2,2.0)-- (6.6,2.4);
\draw (6.6,2.0)-- (7.0,2.4);
\draw (3.2,2.0)-- (3.2,2.4);
\draw (3.6,2.0)-- (3.6,2.4);
\draw (2.8,2.0)-- (2.4,2.4);
\draw (2.8,2.0)-- (2.8,2.4);
\draw (6.2,1.6)-- (5.8,1.2);
\draw (5.8,1.2)-- (6.6,1.6);
\draw (6.6,1.6)-- (7.2,2.0);
\draw (6.6,1.6)-- (6.8,2.0);
\draw (7.2,2.0)-- (7.6,2.4);
\draw (7.2,2.0)-- (7.2,2.4);
\draw [line width=0.4pt,dash pattern=on 1pt off 1pt,color=red] (3.4,1.2)-- (4.,1.2)-- (4.2,1.6)-- (4.4,1.8)-- (4.6000000000000005,1.8)-- (5.,1.6)-- (5.,2.0)-- (5.2,2.2)-- (5.4,2.0)-- (5.6000000000000005,1.6)-- (5.8,1.2)-- (6.4,1.2);
\draw (5.0,0.05) node[anchor=east] {${\buildrel \leftarrow \over \varnothing}$};
\draw (5.0,0.38) node[anchor=east] {$\varnothing$};
\draw[->] (5.17,0.0)-- (5.17,0.4);
\draw (5.17,0.2) node[anchor=west] {$k$ times};

\draw [fill=red] (6.4,0.3) circle (2.5pt);
\draw (6.6,0.3) node[anchor=west] {Vertices $x$ s.t. $\overline{L}_{\overline{T}_k}(x)=1$};
\draw [fill=blue] (6.4,0) circle (2.5pt) ;
\draw (6.6,0) node[anchor=west] {Vertices $x$ s.t. $\overline{L}_{\overline{T}_k}(x)\neq 1$};
\draw [fill=red] (6.4,0.6) circle (2.5pt);
\draw [line width=0.4pt,dash pattern=on 1pt off 1pt,color=red] (6.2,0.6)-- (6.6,0.6);
\draw (6.6,0.6) node[anchor=west] {Vertices of $\mathscr{Z}_k$};
\begin{scriptsize}
\draw [fill=blue] (5.0,0.0) circle (2.5pt);
\draw [fill=blue] (5.0,0.4) circle (2.5pt);
\draw [fill=blue] (4.6,0.8) circle (2.5pt);
\draw [fill=blue] (5.0,0.8) circle (2.5pt);
\draw [fill=blue] (5.4,0.8) circle (2.5pt);
\draw [fill=red] (4.,1.2) circle (2.5pt);
\draw [fill=blue] (4.6,1.2) circle (2.5pt);
\draw [fill=blue] (5.0,1.2) circle (2.5pt);
\draw [fill=blue] (5.4,1.2) circle (2.5pt);
\draw [fill=red] (5.8,1.2) circle (2.5pt);
\draw [fill=blue] (3.2,1.6) circle (2.5pt);
\draw [fill=blue] (4.0,1.6) circle (2.5pt);
\draw [fill=red] (4.2,1.6) circle (2.5pt);
\draw [fill=blue] (4.6,1.6) circle (2.5pt);
\draw [fill=blue] (5.2,1.6) circle (2.5pt);
\draw [fill=blue] (5.4,1.6) circle (2.5pt);
\draw [fill=red] (5.,1.6) circle (2.5pt);
\draw [fill=blue] (6.2,1.6) circle (2.5pt);
\draw [fill=blue] (5.8,1.6) circle (2.5pt);
\draw [fill=blue] (2.8,2.0) circle (2.5pt);
\draw [fill=red] (3.2,2.0) circle (2.5pt);
\draw [fill=blue] (3.6,2.0) circle (2.5pt);
\draw [fill=blue] (3.8,2.0) circle (2.5pt);
\draw [fill=blue] (4.2,2.0) circle (2.5pt);
\draw [fill=red] (4.8,2.0) circle (2.5pt);
\draw [fill=blue] (5.2,2.0) circle (2.5pt);
\draw [fill=blue] (5.6,2.0) circle (2.5pt);
\draw [fill=red] (6.0,2.0) circle (2.5pt);
\draw [fill=red] (6.6,2.0) circle (2.5pt);
\draw [fill=blue] (6.2,2.0) circle (2.5pt);
\draw [fill=blue] (3.8,2.4) circle (2.5pt);
\draw [fill=blue] (4.2,2.4) circle (2.5pt);
\draw [fill=blue] (4.6,2.4) circle (2.5pt);
\draw [fill=red] (5.0,2.4) circle (2.5pt);
\draw [fill=blue] (5.4,2.4) circle (2.5pt);
\draw [fill=blue] (5.6,2.4) circle (2.5pt);
\draw [fill=blue] (5.8,2.4) circle (2.5pt);
\draw [fill=blue] (6.2,2.4) circle (2.5pt);
\draw [fill=blue] (6.6,2.4) circle (2.5pt);
\draw [fill=blue] (7.0,2.4) circle (2.5pt);
\draw [fill=blue] (3.2,2.4) circle (2.5pt);
\draw [fill=blue] (3.6,2.4) circle (2.5pt);
\draw [fill=blue] (2.4,2.4) circle (2.5pt);
\draw [fill=red] (2.8,2.4) circle (2.5pt);
\draw [fill=blue] (6.6,1.6) circle (2.5pt);
\draw [fill=red] (7.2,2.0) circle (2.5pt);
\draw [fill=blue] (6.8,2.0) circle (2.5pt);
\draw [fill=blue] (7.6,2.4) circle (2.5pt);
\draw [fill=blue] (7.2,2.4) circle (2.5pt);
\end{scriptsize}
\end{tikzpicture}
\caption{The set $\mathscr{Z}_k$.}
\label{f:Zk}
\end{figure}

  \medskip
  
  The rest of this paper is organized as follows:
   \begin{itemize*}
  \item Section~\ref{s:pre}: we present the main technical tools, in particular the many-to-one formula for the branching random walk and a change of measure for the edge local times;
\item Section~\ref{s:thm1}: we give the proof of the part~\eqref{t:kappa>2} in Theorem~\ref{t:main};
\item Section~\ref{s:thm2}: we prove Theorem~\ref{t:tailmaxlocaltime} by establishing some results on an associated Markov chain which appears naturally in the change of measure for the edge local times;
\item Section~\ref{s:kappale2}: we prove the remaining parts  in Theorem~\ref{t:main}, namely~\eqref{t:kappa=2},~\eqref{t:kappa<2} and~\eqref{t:kappa<2b}  by using Theorem~\ref{t:tailmaxlocaltime}, and we give the proofs of Corollaries~\ref{cor:cardinal} and~\ref{cor:limsup}.
\end{itemize*}

   Throughout this paper, we denote by $c, c', c''$ (eventually with some subscripts) some positive constants whose values can change from one paragraph to another.

  \section{Preliminaries}\label{s:pre}

%%$\phantom{aob}$

This section is divided into three subsections: in the first subsection, we  introduce the potential $V$ of the biased walk $(X_n)_{n\ge 0}$, and recall some known results  on the potential $V$ and on an associated one-dimensional random walk $S$; in the second subsection, we present a change of measure formula for the edge local times and some consequences; in the third (and last) subsection we collect some known  facts on the random walk $(X_n)_{n\ge 0}$ which will be used in the proofs of Theorems~\ref{t:main} and~\ref{t:tailmaxlocaltime}.

%We shall need an estimate on the maximal displacement of the walk up to time $\overline T_n$, which is a straightforward consequence of the tail estimate on the height of an excursion:

\subsection{The potential} 
Let us introduce $V=(V(x))_{x\in \T}$ the {\it random potential} of the biased random walk $(X_n)_{n\ge 0}$, which will completely determine the behavior of $(X_n)$. Define  
%{\tt (xxxx Est-il n\'ecessaire de pr\'eciser que $V({\buildrel \leftarrow \over \varnothing}):=0$~?)} 
\begin{equation}
V(\varnothing):=0\qquad{\textrm and}\qquad
  V(x) 
  := 
  -
  \sum_{y\in \, ]\!] \varnothing,\, x]\!]}
  \log \, A(y),
  \qquad \textrm{for }x\in \T\backslash\{ \varnothing\} \, ,
  \label{V}
\end{equation}

\noindent where $\, ]\!] \varnothing, \, x]\!] := [\![ \varnothing, \, x]\!] \backslash \{ \varnothing\}$, with $[\![ \varnothing, \, x]\!]$ denoting the set of vertices (including $x$ and $\varnothing$) on the unique shortest path connecting $\varnothing$ to $x$. 
%Throughout the paper, we use $x_i$ (for $0\le i\le |x|$) to denote the ancestor of $x$ in the $i$-th generation; in particular, $x_0 = \varnothing$ and $x_{|x|} =x$. As such, the potential $V$ in (\ref{V}) can also be written as
%$$
%V(x) 
%= - \sum_{i=0}^{|x|-1} \log \,
%A(x_{i+1})
%=
%- \sum_{i=0}^{|x|-1} \log \,
%\frac{\omega(x_i, \, x_{i+1})}{\omega(x_i, \, x_{i-1})},
%\qquad x\in \T\backslash\{ \varnothing\} \, .
%\qquad
%(x_{-1} := {\buildrel \leftarrow \over \varnothing})
%$$
The process $(V(x),{x\in \T})$ is a branching random walk, in the usual sense of Biggins~\cite{biggins77}.

 Let us define a symmetrized version of the potential which will naturally appear in the study of local times: 
%$U({\buildrel \leftarrow \over \varnothing}):=0$ and 
\begin{equation*}\label{U}
  U(x) 
  := 
  V(x) - \log (\frac{1}{\omega(x, \, {\buildrel \leftarrow \over x})}) \, ,
  \qquad
  x\in \T \, .
\end{equation*}

 \noindent Note that 
\begin{equation*}
  \ee^{-U(x)}
  = 
  \frac{1}{\omega(x, \, {\buildrel \leftarrow \over x})} \, \ee^{-V(x)} 
  =
  \ee^{-V(x)} + \sum_{y\in \T: \, {\buildrel \leftarrow \over y} =x} \ee^{-V(y)} ,
  \qquad
  x\in \T \, .
\end{equation*}

% The law of large numbers established by  Hammersly~\cite{H74}, Kingman~\cite{K75} and Biggins~\cite{B76} says that $\frac1n\, \min_{|x|=n} V(x) \to c$, $\p^*$-a.s.  where $c:=...?$ is a positive constant. In particular, $\p^*$-a.s., $V(x) \to \infty$ as $|x|\to \infty$; 
Recall that $\P^*$ is the probability $\P$ conditioned on the  non-extinction of the Galton-Watson tree $\T$. The  Biggins-Hammersley-Kingman (\cite{H74, K75, B76}) law of the large numbers implies that under the assumptions~\eqref{hyp1} and~\eqref{kappa}, $\P^*$-almost surely, $\frac1n\, \min_{|x|=n} V(x)$ converges towards some positive constant. The following simple result deals with  the symmetrized potential:
 
 \begin{lemma}\label{l:U} Assume~\eqref{hyp1},~\eqref{kappa} and~\eqref{hyp2}.  As $n\to \infty$, $$
 \min_{|x|=n} U(x) \to \infty, \qquad  \P^*\mbox{-a.s.}$$
 \end{lemma}

{\noindent\it Proof of Lemma~\ref{l:U}.} According to~\eqref{kappa}, we may choose  some constant $t>1$ such that $$\E \Big( \sum_{x\in \T: \, |x|=1} \ee^{- t V(x)} \Big) < 1, \qquad \mbox{and } \quad c:= \E\big(1+ \sum_{x\in \T: \, |x|=1} \ee^{-  V(x)} \big)^t < \infty.$$

 \noindent By the branching property, we have that for any $n\ge1$, \begin{eqnarray*} 
 \E \Big( \sum_{x\in \T: \, |x|=n} \ee^{- t U(x)} \Big)
 &=&
 \E \Big( \sum_{x\in \T: \, |x|=n} \ee^{- t V(x)} \big(1+ \sum_{y:{\buildrel \leftarrow \over y}=x} \ee^{- (V(y)- V(x))}\big)^t\Big)
 \\
 &=&
 c\, \E \Big( \sum_{x\in \T: \, |x|=n} \ee^{- t V(x)} \Big) 
 \\
 &=&
 c\, \Big( \E \big( \sum_{x\in \T: \, |x|=1} \ee^{- t V(x)} \big) \Big)^n, 
\end{eqnarray*}

\noindent whose sum on $n$ converges. Therefore  Borel-Cantelli's lemma yields that $\sum_{x: \, |x|=n} \ee^{- t U(x)} \to 0$, $\P$-a.s., and Lemma~\ref{l:U} follows. \hfill$\Box$

\medskip
Consequently, the set   of minimums  of $U$, defined as below, is finite almost surely:
\begin{equation}\label{MinU}
{\mathscr M}:= \Big\{x \in \T: U(x) = \min_{y \in \T} U(y)\Big\}.
\end{equation}

 \medskip
 
 Now we introduce   the  ``many-to-one formula", which is by now a standard tool in the study of branching random walk, see for instance Lyons~\cite{lyons} for the probabilistic construction and  Shi~\cite{stf} for the complete references. Under assumption~\eqref{hyp1}, there exists a sequence of i.i.d.\ real-valued random variables $(S_i-S_{i-1}, \, i\ge 0)$, with $S_0=0$, such that for any $n\ge 1$ and any Borel function $g: \r^n \to \r_+$,
\begin{equation}
  \E \Big[ \sum_{x\in \T: \, |x|=n} g(V(x_i), \, 1\le i\le n) \Big]
  =
  \E \Big[ \ee^{S_n} \, g(S_i, \, 1\le i\le n) \Big] \, ,
  \label{many-to-one}
\end{equation}

\noindent where, for any vertex $x\in \T$ such that $|x|=n$, $x_i$ ($0\le i\le n$) denotes the ancestor of $x$ in the $i$-th generation.

Observe that by~\eqref{hyp1}, $$ \E \Big( S_1\Big)= \E \Big( \sum_{|u|=1} V(u) \ee^{-V(u)}\Big)= -\E \Big( \sum_{i=1}^\nu \, A_i \log A_i \Big) >0,$$ and   if $1< \kappa< \infty$, then by~\eqref{kappa} and~\eqref{hyp2}, $$ \E\Big( \ee^{-(\kappa-1) S_1}\Big)=1, \qquad \mbox{ and } \qquad \E \ee^{\delta S_1} + \E \ee^{-(\kappa-1+\delta) S_1} < \infty,$$

\noindent  for some positive constant $\delta$. We mention that $\kappa=\infty$ corresponds to the case where $\E\big( \ee^{- t S_1}\big) < 1$ for any $t>0$, which implies that $\P(S_1\ge0)=1$.

% Let us quote at first a random walk result du to Feller (XII.5.13, pp.411): Assume that $\E(|S_1| \ee^{- (\kappa-1) S_1})< \infty$ and $S_1$ is no-lattice distributed. Then $$ \P \Big( \min_{k\ge0} S_k < - \lambda\Big) \, \sim\, c_F\, \ee^{-(\kappa-1) \lambda } , \qquad \lambda \to \infty.$$ Moreover $\E \Big( \sum_{k=0}^{e_1-1} 1_{( \max_{0 < l \le k} (- S_l) \le r)} \Big)$ is the renewal function related to $S$ hence converges to some positive constant $R_\infty$ as $r\to\infty$.

 We summarize some known results on the transient random walk $(S_n)$ in the following fact:

 \begin{fact}\label{f:rw} Assume~\eqref{hyp1},~\eqref{kappa} and~\eqref{hyp2}. 
 
(i) If $1< \kappa< \infty$, then \begin{equation}\label{kesten}
  \P\Big(\sum_{j=0}^\infty \ee^{-S_j } > r \Big) \, \asymp \, r^{-(\kappa-1)} , \qquad r \to \infty. \end{equation}
 
 (ii)  If $\kappa=\infty$, then for any $p>0$, there exists some constant $c_p>0$ such that 
 \begin{equation}\label{kappainfini}
  \P\Big(\sum_{j=0}^\infty \ee^{-S_j } > r \Big)
  \le 
  c_p\, r^{-p}, \qquad \forall \, r>0.
 \end{equation}
 \end{fact}

 We mention that~\eqref{kesten} comes from Kesten~\cite{kesten73} [the non-lattice case] and Grintsevichyus~\cite{Grintsevichyus} [for the lattice case], and~\eqref{kappainfini} follows easily from the triangular inequality: $\| \sum_{j=0}^\infty \ee^{-S_j }\|_p \le \sum_{j=0}^\infty \| \ee^{-S_j}\|_p= \sum_{j=0}^\infty \big(\| \ee^{-S_1}\|_p\big)^j $, where $\|\cdot\|_p$ denotes the $L^p$-norm, and $\| \ee^{-S_1}\|_p<1$ since $\kappa=\infty$.

 \subsection{Change of measure for the edge local times}
 
 In this subsection we introduce a change of measure formula for the edge local times. This formula  describes the law of the  local times process  under the annealed probability $\p$,  and   plays an important role in the proof of Theorem~\ref{t:tailmaxlocaltime}.

 Recall~\eqref{overtj}. 
For any $k\ge1$, let 
\begin{equation}\label{def-Tk}
\T^{(k)}:= \{x\in \T: \overline L_{\overline T_k}(x) \ge 1\},
\end{equation}

\noindent be the subtree formed by the vertices visited at least once by the walk up to time $\overline T_k$. As  proved by A\"\i d\'ekon and de Raph\'elis (\cite{AR15}, Lemma 3.1), the marked tree $( \T^{(k)}, (\overline L_{\overline T_k}(x))_{x\in \T^{(k)}}) $ is a multi-type Galton-Watson tree with  initial type $k$ at the root (the type of a vertex $x\in \T^{(k)}$ is exactly the edge-local time $\overline L_{\overline T_k}(x)$). 
Applying~\cite{KLPP} and~\cite{BK04}, we get the following fact: 

\begin{fact}[change of measure for the edge-local times] \label{f:measurechange} Assume~\eqref{hyp1} and~\eqref{kappa}. For any $k\ge1$, on an eventually enlarged probability space we may define a probability measure $\widehat \p_k$, a distinguished infinite ray called spine $\{\w_i, i\ge0\}$ and the tree $\T^{(k)}$ such that 

(i) for any $n\ge0$, $$
\widehat \p_k\Big( \w_n= x \, \big| \, {\cal F}_n\Big) = \frac{\overline L_{\overline T_k}(x)}{\sum_{|y|=n} \overline L_{\overline T_k}(y)}, \qquad \forall |x|=n, \, \, x \in \T^{(k)},$$
with ${\cal F}_n=\sigma\{\overline L_{\overline T_k}(x), x\in \T^{(k)}: |x| \le n\}$. 
 
(ii) the marginal of $\widehat\p_k$ on the space of trees is absolutely continuous with respect to $\p$:  $$
\frac{\d \widehat \p_k}{\d \p} \Big|_{{\cal F}_n}= \frac1{k}\sum_{|y|=n} \overline L_{\overline T_k}(y), \qquad \forall \, n\ge 0. $$

(iii) under $\widehat\p_k$, conditioned on $\sigma\{ \overline L_{\overline T_k}(x), {\buildrel \leftarrow \over x} \in \{ \w_n, n\ge0\}\}$, the processes $\{\overline L_{\overline T_k}(y), y \in \T^{(k)}_x\}_{{\buildrel \leftarrow \over x} \in \{ \w_n, n\ge0\}}$ are independent and are distributed as $\{\overline L_{\overline T_l}(y), y \in \T^{(l)}_x\}$ with $l:=\overline L_{\overline T_k}(x)$, where $\T^{(k)}_x$ denotes the subtree of $\T^{(k)}$ rooted at $x$. 
 \end{fact}

As a consequence of (i) and (ii), we get a many-to-one formula for the edge-local times. Let $n\ge1$. For any $k\ge1$ and measurable nonnegative function $f: \N^n \to \r_+$, we have
 \begin{equation}\label{ar15}
 \e\Big[ \sum_{|x|=n} 1_{\{\overline L_{\overline T_k}(x)\ge 1\}} f ( \overline L_{\overline T_k}(x_1), \overline L_{\overline T_k}(x_2), ..., \overline L_{\overline T_k}(x_n))\Big]
 =
 \widehat \e_k \Big[\frac{k}{Y_n} f(Y_1, Y_2, ..., Y_n)\Big],
 \end{equation}
 
 \noindent 
 with    \begin{equation}\label{defyn} 
 Y_n:=\overline L_{\overline T_k}(\w_n),  \qquad n\ge0.
 \end{equation}

 The law of the process $(Y_n)_{n\ge 0}$ under $\widehat\p_k$, is described  in Subsection~6.1 of  A\"\i d\'ekon and  de Raph\'elis~\cite{AR15}:

 \begin{fact}[\cite{AR15}] \label{f:ar} For any $k\ge1$, under $\widehat\p_k$, $(Y_n)_{n\ge 0}$ is a positive recurrent Markov chain taking values in $\N \backslash\{0\}$, started at $k$, with transition probabilities given by 
 \begin{equation}\label{pij}
 P_{i, j}:= C_{i+j-1}^{j-1} \, \E \left( \frac{\ee^{-(j-1) S_1}}{(1+\ee^{-S_1})^{i+j}}\right) , \qquad i, j \ge 1,
 \end{equation}
 
 \noindent 
 where the law of $S_1$ is given in~\eqref{many-to-one}. 
 Moreover, the invariant probability measure $(\pi_i)_{i\ge 1}$ of $(Y_n)_{n\ge 0}$  is given as follows: \begin{equation}\label{pi}
 \pi_i:= i \, \E\left[\frac{\big(\sum_{n=1}^\infty \ee^{-S_n}\big)^{i-1}}{\big(1+\sum_{n=1}^\infty \ee^{-S_n}\big)^{i+1}}\right], 
 \qquad i\ge 1. \end{equation}
 \end{fact} 

 \medskip

%%As a consequence,  we have $ \widehat \p_1(Y_n =i) \le \frac{\pi_i}{\pi_1}$ for any $i$ and $n$. 

% By  the triangular inequality and~\eqref{pi2}, we get that if $\kappa=\infty$, then for any $r>0$, there exists some constant  $c=c_r>0$ such that for all $N\ge 1$ and $z> 1$,  
%\begin{equation}\label{kappainfini}
% \widehat \p_1\Big(\max_{1\le n \le N} Y_n \ge z\Big)
 %\le
% c\, \frac{N}{z^r}.
%\end{equation} 

 By using~\eqref{kesten} and~\eqref{kappainfini}, we get    the asymptotic behaviors of  the invariant probability measure $\pi_i$ as $i \to \infty$:   if $1 < \kappa < \infty$, then  \begin{equation}\label{pi1}
 \pi_i \, \asymp\,  i^{-\kappa} , \qquad i \to \infty, 
 \end{equation}
 
 \noindent and if $\kappa=\infty$, then for any $p>1$, \begin{equation}\label{pi2} \sup_{i\ge1} i^p \, \pi_i < \infty. 
 \end{equation}

\subsection{The random walk $(X_n)_{n\ge 0}$}
We collect some known results on the almost sure limits of the random walk $(X_n)_{n\ge 0}$ in the following fact:

\begin{fact} \label{f:1} Assume~\eqref{hyp1},~\eqref{kappa} and~\eqref{hyp2}.  For any $\kappa \in (1, \infty]$,  $\p^*$-almost surely, \begin{eqnarray} \label{maxlimsup2} 
 \lim_{n\to\infty} \frac1{\log n} \log \max_{0\le i\le n} |X_i| &= & 1 - \max( {1\over 2}, { 1 \over \kappa}) , 
 \\  
 \label{aslocaltime}
 \lim_{n\to\infty} \frac1{\log n} \log L_n({\buildrel \leftarrow \over \varnothing})&=& \max(\frac1\kappa, \frac12) ,
 \end{eqnarray}
where $L_n({\buildrel \leftarrow \over \varnothing}):= \sum_{i=1}^n 1_{\{ X_i = {\buildrel \leftarrow \over \varnothing}\}}$. \end{fact}

%For any $\varepsilon>0$, we have \begin{equation} \label{tail}
%\p\Big( \max_{0\le i \le \overline T_1}|X_i| \ge n \Big) 
%\le 
%n^{-\max(\frac1{\kappa-1}, 1)+\varepsilon}.
%\end{equation}
%Consequently for any $\varepsilon>0$, $\p$-almost surely, \begin{equation} \label{limsup}
% \limsup_{n\to \infty} n^{- \min(\kappa-1, 1) - \varepsilon} \max_{0\le i \le \overline T_n}|X_i| =0.
% \end{equation} 

We mention that~\eqref{maxlimsup2} was proved in~\cite{yzptrf} under  more restrictive assumptions (i.e.\ if $\nu$ equals some constant and $(A_1, ..., A_\nu)$ are i.i.d.), but the same argument still holds in the present case. 
%Otherwise, \eqref{maxlimsup2} can be seen as a corollary of Theorem~1 of~\cite{loic16} for $\kappa\in(1;2]$ and as a corollary of Theorem~1.1 of~\cite{AR15} for $\kappa>2$. 
The statement~\eqref{aslocaltime} was implicitly contained in~\cite{andreoletti-debs1, yzptrf}, see~\cite{subdiff} for further studies on the local times.

Let us consider now the large deviations of the local times at a single vertex of the tree. Let $x \in \T\cup \{{\buildrel \leftarrow \over \varnothing}\}$. Define \begin{eqnarray*}
T_x &:=& \inf\{n\ge0: X_n=x\}, 
\\
T_x^+&:=& \inf\{n\ge 1: X_n= x\},
\end{eqnarray*}

\noindent and we denote by $P_{x, \omega}$ the quenched probability under which the random walk $(X_n)_{n\ge 0}$ starts at $x$ (so $P_\omega=P_{\varnothing, \omega}$). Observe that $\overline T_1= T_{{\buildrel \leftarrow \over \varnothing}}+1$, thus for any  $x\in\T \backslash \{\varnothing\}$, we get $$ P_\omega \big(T_x < \overline T_1\big) 
     =
     P_\omega \big(T_x <  T_{{\buildrel \leftarrow \over \varnothing}}\big) ,
     \qquad
      P_{x,\omega} \big( \overline T_1 < T_x^+ \big)
  =
   P_{x,\omega} \big( T_{{\buildrel \leftarrow \over \varnothing}}< T_x^+ \big) .
$$

 The probabilities $ P_\omega \big(T_x <  T_{{\buildrel \leftarrow \over \varnothing}}\big)$ and $P_{x,\omega} \big( T_{{\buildrel \leftarrow \over \varnothing}}< T_x^+ \big)$ only involve  the restriction at $\{ {\buildrel \leftarrow \over \varnothing}\}\cup [\! [\varnothing, \, x]\!]\,$ of the biased walk $(X_n)_{n\ge 0}$, so a standard result for one-dimensional birth and death  chains (\cite{hoel-port-stone}, pp.31,  formulae (59) and (60)) tells us that \begin{eqnarray}
     P_\omega \big(T_x < \overline T_1\big) 
     &=&\frac{1}{\sum_{z\in \, [\![ \varnothing, \, x]\!]} \ee^{V(z)}}   \, , \label{1D-MAMA}
  \\
  P_{x,\omega} \big( \overline T_1 < T_x^+ \big)
   &=& \frac{\ee^{U(x)}}{\sum_{z\in \, [ \! [ \varnothing, \, x]\! ]} \ee^{V(z)}} \,,  \label{1D-MAMA-bis}
\end{eqnarray}
  
 \noindent where $U(x)$ was defined in~\eqref{U}.

\medskip

 For any $x\in\T \backslash \{\varnothing\}$, the law of $L_{\overline T_n} (x)$ under $P_\omega$ is the law of $\sum_{i=1}^n \xi_i$, where $(\xi_i, \, i\ge 1)$ are  i.i.d.\ random variables with common law given as follows: $P_\omega(\xi_1 =0) =1-a$ and $P_\omega(\xi_1 \ge k) = a \, p^{k-1}$, $\forall k\ge 1$, where 
\begin{eqnarray}
  1-p 
 &:=& P_{x,\omega} \{ \overline T_1< T_x^+ \}
  =
  \frac{\ee^{U(x)}}{\sum_{z\in \, [ \! [ \varnothing, \, x]\! ]} \ee^{V(z)}} \, ,
  \label{valeur(p,a)1}
  \\
  a
 &:=& P_\omega \{ T_x < \overline T_1\} 
  =
  \frac{1}{\sum_{z\in \, [ \! [ \varnothing, \, x]\! ]} \ee^{V(z)}} \, .
  \label{valeur(p,a)2}
\end{eqnarray}

We shall use several times the following lemma which gives the tail estimate of $L_{\overline T_n} (x)$ under $P_\omega$.

\begin{lemma}
\label{l:sum_iid}

 Let $0<a<1$ and $0<p<1$. Let $(\xi_i, \, i\ge 1)$ be an i.i.d.\ sequence of random variables with $\P(\xi_1 =0) =1-a$ and $\P(\xi_1 \ge k) = a \, p^{k-1}$, $\forall k\ge 1$.  Let $n, k\ge2$. 
 If  $\frac{a}{1-p} < \frac{k}{8 n} $, then 
   $$
 \P\Big\{ \sum_{i=1}^n \xi_i \ge k \Big\}
  \le
  6\, n\, a\, \ee^{- \frac{(1-p)k }{8}}\, .
 $$
 \end{lemma}

{\noindent\it Proof of Lemma~\ref{l:sum_iid}.} The above estimate   was borrowed from~\cite{yzfavtree} when $k= \lceil \varepsilon n \rceil$ with some $0< \varepsilon<1$, indeed the same proof presented therein  holds for all $k$ satisfying $\frac{a}{1-p} < \frac{k}{8 n} $, without any modification.  \hfill$\Box$

 \medskip
We end this subsection by a   useful relationship between the edge local times $\overline L_{\overline T_1}(\cdot)$ and the (site) local times $L_{\overline T_1}(\cdot)$: \begin{equation}\label{siteedge}
 L_{\overline T_1}(x)
 = \overline L_{\overline T_1}(x) + \sum_{y: {\buildrel \leftarrow \over y}= x } \overline L_{\overline T_1}(y) 
 =\overline L_{\overline T_1}(x) + \Theta(x),
  \qquad x\in \T , 
 \end{equation}

\noindent where for the notational brevity, we write  \begin{equation}\label{Theta} \Theta(x):= \sum_{y: {\buildrel \leftarrow \over y}= x } \overline L_{\overline T_1}(y) , \qquad x\in \T.
\end{equation}

 For $x\in \T$, let $(x^{(1)}, ..., x^{(\nu_x)})$ be as before  the set of the children of $x$. For any $0 \le s_1, ..., s_i, ... \le1$, we have $$
 E_\omega\Big( \prod_{i=1}^{\nu_x} (s_i)^{\overline L_{\overline T_1}(x^{(i)})} \, \big|\, \overline L_{\overline T_1}(x)= k \Big)
 =
 \Big(1+ \sum_{i=1}^{\nu_x} (1- s_i) A(x^{(i)})\Big)^{-k}, \qquad \forall\, k\ge 1.$$

 \noindent Then, for $0\le s \le 1$, \begin{equation} \label{generating-Thetax}
 E_\omega\Big( s^{\Theta(x)}\, \big|\, \overline L_{\overline T_1}(x)= k \Big)
 =
 \Big(1+  (1-s)  \sum_{y: {\buildrel \leftarrow \over y}= x}  A(y) \Big)^{-k} , \end{equation}

 \noindent which means that under $P_\omega$ and conditioned on $\{\overline L_{\overline T_1}(x)= k\}$, \begin{equation}\label{lawThetax}
 \Theta(x)
 \, \law\,
 \sum_{i=1}^k \eta_i,
 \end{equation} where $(\eta_i,i\ge1)$ are  i.i.d.\ geometric variables with $P_\omega(\eta_1\ge n)= \Big(\frac{ \sum_{y: {\buildrel \leftarrow \over y}= x}  A(y) }{1+\sum_{y: {\buildrel \leftarrow \over y}= x}  A(y) }\Big)^n$, for all $n\ge0 $.  The identity~\eqref{lawThetax} will be explored in Section~\ref{s:thm2}.

\section{Proof of Theorem~\ref{t:main}: Case $\kappa\in (2, \infty]$}\label{s:thm1}

The part~\eqref{t:kappa>2} in Theorem~\ref{t:main} follows  from a result which will also be useful for  $\kappa \in (1, 2]$. 

Let $0<\varepsilon< 1$. We define an integer-valued random variable $K_\varepsilon(\omega)$ by \begin{equation} \label{Komega}
 K_\varepsilon(\omega)
 :=
 \begin{cases} \sup\big\{n\ge1: \min_{|x|=n} U(x) < \frac{8}{\varepsilon} \big\}, \qquad & \mbox{on } \{\T=\infty\}, 
 \\
 \max_{x\in \T} |x|, \qquad & \mbox{on } \{\T <\infty\}.
\end{cases}
 \end{equation} 
 
 \noindent  By Lemma~\ref{l:U},  $K_\varepsilon< \infty$, $\P$-a.s.  The following result holds for any $\kappa \in (1, \infty]$:

 \begin{proposition}\label{p:main2} Let $\kappa \in (1, \infty]$. Assume~\eqref{hyp1},~\eqref{kappa} and~\eqref{hyp2}.  Let $0< \varepsilon< 1/2$. We have \begin{equation}\label{limsupltn}
 \limsup_{n\to \infty} \frac1n\, 
 \,  \max_{  K_\varepsilon\le |x| \le n^{q }}
  L_{\overline T_n} (x) 
  < 
  \varepsilon, \qquad \p\mbox{-a.s.},
  \end{equation}
  where $q>0$ denotes an arbitrary constant smaller than $\kappa-1$ if $\kappa \in (1, \infty) $, and $q$ denotes some (fixed) constant strictly larger than $1$ if $\kappa = \infty$. 
  
 Consequently for any $0<b < \max(\frac{q}{2}, \frac{q}{\kappa})$,  we have \begin{equation}\label{favempty} 
\bigcap_{n\ge 1} \bigcup_{j\ge n} \,  \Big\{ \exists \, x\in \T: K_\varepsilon \le |x| \le j^{b} \mbox{ such that } x \in {\mathbb F}(j)\Big\}
 = 
 \emptyset, \qquad \p\mbox{-a.s.}
  \end{equation}
  \end{proposition}

% $$q\equiv q(\kappa, \delta):= \begin{cases} \kappa-1 - \delta, & \mbox{ if $\kappa\in (1, 2]$} , 
 % \\
 % 1+\delta, & \mbox{ if $\kappa\in (2, \infty]$}.
 % \end{cases}$$ 

By admitting Proposition~\ref{p:main2} for the moment,  we immediately get~\eqref{t:kappa>2}:

\medskip

\noindent{\it Proof of~\eqref{t:kappa>2} in Theorem~\ref{t:main}. } Let $\kappa \in (2, \infty]$. Then we may choose $q>1$ so that~\eqref{favempty} holds with some $b>\frac12$. Fix  $0< \varepsilon < \frac12$.  In view of~\eqref{maxlimsup2}, $\p$-a.s.  for all large $n$, $\max_{0\le i \le n} |X_i| \le n^b$, hence 
\begin{equation} \label{limsupke30} \limsup_{j\to \infty} \sup_{x\in {\mathbb F}(j)} | x|\le K_\varepsilon. 
\end{equation} 

Now we notice that  under $P_\omega$,  $(\ee^{-U(x)})_{x \in \T\cup\{{\buildrel\leftarrow \over \varnothing}\}}$ is the invariant measure of the Markov chain $(X_n)$.  By applying the ergodic theorem for additive functionals of a recurrent Markov chain, we get that  under $P_\omega$, for any $x \in \T$, $$ \lim_{n\to \infty} \frac{L_n(x)}{L_n(\varnothing)} = \ee^{-(U(x)-U(\varnothing))}.$$

\noindent  It follows from \eqref{limsupke30}   that $\p$-a.s.  for all large $n$,  $$ {\mathbb F}(n) \, \subset\, \Big\{x \in \T: U(x) = \min_{ y\in \T:  |y|\le K_\varepsilon} U(y)\Big\}. $$

\noindent Finally, we remark that  $ \min_{ y\in \T: |y|\le K_\varepsilon} U(y)=  \min_{ y\in \T } U(y)$ and get~\eqref{t:kappa>2}. \hfill$\Box$

\medskip
It remains to give the proof of Proposition~\ref{p:main2}, whose main ingredient is  contained in the following lemma:

\begin{lemma}\label{l:pre} Assume~\eqref{hyp1},~\eqref{kappa} and~\eqref{hyp2}.  Let $\kappa \in (1, \infty]$. For any small constant $0< \delta< \kappa-1$, there  exists some constant $c=c_\delta>0$ such that for any $k\ge1$ and $r>1$, 
\begin{equation}\label{pre1}
 \E \sum_{|x|=k} \ee^{-U(x)}  \exp\Big( - r\,  \frac{\ee^{U(x)}}{\sum_{z\in \, [ \! [\varnothing, \, x]\! ]} \ee^{V(z)}} \Big) 
 \le
 \begin{cases}
 c\,  r^{-(\kappa-1-\delta)}, \qquad & \mbox{if } \kappa\in (1, \infty), \\
 c\,  r^{-1-\delta}, \qquad & \mbox{if } \kappa = \infty .
 \end{cases}
  \end{equation}
\end{lemma} 
 
 {\noindent\it Proof of Lemma~\ref{l:pre}:}  Let $r$ be large. Denote by $H_r(x):=  \ee^{-U(x)}  \exp\big( - r\,  \frac{\ee^{U(x)}}{\sum_{z\in \, [ \! [\varnothing, \, x]\! ]} \ee^{V(z)}} \big) $. If $ \sum_{z\in \, [ \! [ \varnothing, \, x]\! ]} \ee^{V(z) -U(x)} \le \frac{r}{(\log r)^2} $ then $H_r(x) \le  \ee^{-U(x)} \ee^{-  (\log r )^2}$, hence
 $$
 \E \sum_{|x|=k} H_r(x) 1_{\{\sum_{z\in \, [ \! [ \varnothing, \, x]\! ]} \ee^{V(z) -U(x)} \le \frac{r}{(\log r)^2}\}} 
\le
 \ee^{-  (\log r )^2} \, \E \sum_{|x|=k} \ee^{-U(x)} 
=
 2 \, \ee^{-  (\log r )^2},
 $$
 
 \noindent where we used the fact that $$\E \sum_{|x|=k} \ee^{-U(x)} 
 = \E \sum_{|x|=k} \ee^{-V(x)} + \E \sum_{|y|=k+1} \ee^{-V(y)} 
 =2.$$

 It remains to deal with the vertices $|x|=k$ such that $ \sum_{z\in \, [ \! [ \varnothing, \, x]\! ]} \ee^{V(z) -U(x)} > \frac{r}{(\log r)^2} $. Let $\varepsilon>0$ be small. Recall that $\ee^{-U(x)}= \ee^{-V(x)} \frac1{\omega(x, {\buildrel\leftarrow \over x})}$. If $\sum_{z\in \, [ \! [ \varnothing, \, x]\! ]} \ee^{V(z) -V(x)} \in (r^{\varepsilon(j-1)}, r^{\varepsilon j}]$ for some $1\le j \le \frac1\varepsilon$, then $ \frac1{\omega(x, {\buildrel\leftarrow \over x})} > r^{1- \varepsilon j} (\log r)^{-2}.$ It follows that 
 \begin{eqnarray}
 && \E \sum_{|x|=k} H_r(x) 1_{\{\sum_{z\in \, [ \! [ \varnothing, \, x]\! ]} \ee^{V(z) -U(x)} > \frac{r}{(\log r)^2}\}} 
  \nonumber \\
 &\le&
  \sum_{j=1}^{\lceil 1/\varepsilon\rceil} \E\, \sum_{|x|=k} 1_{\{\sum_{z\in \, [ \! [ \varnothing, \, x]\! ]} \ee^{V(z) -V(x)} \in (r^{\varepsilon(j-1)}, r^{\varepsilon j}] \}} \ee^{-U(x)} 1_{\{ \frac1{\omega(x, {\buildrel\leftarrow \over x})} > r^{1- \varepsilon j} (\log r)^{-2}\}}  \nonumber\\
 && \qquad +  \E\, \sum_{|x|=k} 1_{\{\sum_{z\in \, [ \! [ \varnothing, \, x]\! ]} \ee^{V(z) -V(x)} > r\}} \ee^{-U(x)} 
 \nonumber \\
 &=:& A_{\eqref{en2}} + B_{\eqref{en2}}, \label{en2}
\end{eqnarray}

\noindent with obvious definitions of $ A_{\eqref{en2}} $ and $ B_{\eqref{en2}}$.

Observe that for any $k\ge1$, conditioned on $\{V(z), |z|\le k\}$, $\{1/\omega(x, {\buildrel\leftarrow \over x})\}_{|x|=k}$ are i.i.d.\ and are distributed as $1+ \sum_{|u|=1} A(u)$ whose expectation is equal to $2$. It follows that 
\begin{eqnarray*}
B_{\eqref{en2}}
&=&  \E \sum_{|x|=k} 1_{\{\sum_{z\in \, [ \! [ \varnothing, \, x]\! ]} \ee^{V(z) -V(x)} > r\}} \ee^{-V(x)} \times \E(1+ \sum_{|u|=1} A(u))
\\
&=&
2 \, \P \Big( \sum_{i=1}^k \ee^{S_i - S_k}  > r \Big),
\end{eqnarray*} 

\noindent by the many-to-one formula~\eqref{many-to-one}. Remark that for any $k\ge 1$ and any $r$, $$ \P \big( \sum_{i=1}^k \ee^{S_i - S_k}  > r \big)
=
 \P \big( \sum_{i=0}^{k-1} \ee^{-S_i } > r\big),
 $$
 
 \noindent by the time-reversal of $S$.  We choose and fix $$
 \eta := \kappa-1, \qquad \mbox{ if } \kappa<\infty,$$
 
 \noindent and we take an arbitrary constant $\eta>1$ if $\kappa=\infty$.  Applying~\eqref{kesten} and~\eqref{kappainfini}, we get that  for all large $r$, $\P \big( \sum_{i=0}^\infty \ee^{-S_i } > r\big) \le r^{-\eta}$; hence $$
B_{\eqref{en2}}
\le 
2 \,  r^{- \eta }
,$$

\noindent for all large $r$.

 To deal with $A_{\eqref{en2}}$, we remark by~\eqref{hyp2} that $c_\alpha:= \E(\frac{1}{\omega(\varnothing, {\buildrel\leftarrow \over \varnothing})} )^{\alpha}=\E\Big(1+\sum_{i=1}^{\nu}A_i\Big)^\alpha < \infty$ for some $\alpha>\kappa$ when $1 < \kappa < \infty$ and for some $\alpha >2$ when $\kappa=\infty$. Therefore $\E (\frac{1}{\omega(\varnothing, {\buildrel\leftarrow \over \varnothing})} 1_{\{ \frac{1}{\omega(\varnothing, {\buildrel\leftarrow \over \varnothing})} > t\}} ) \le c_\alpha \, t^{- (\alpha-1)} $ for any $t >0$. It follows that \begin{eqnarray*}
A_{\eqref{en2}}&\le&
c_\alpha \,  \sum_{j=1}^{\lceil 1/\varepsilon\rceil} \E\, \sum_{|x|=k} 1_{\{\sum_{z\in \, [ \! [ \varnothing, \, x]\! ]} \ee^{V(z) -V(x)} \in (r^{\varepsilon(j-1)}, r^{\varepsilon j}]\}} \ee^{-V(x)} r^{- (\alpha-1) (1- \varepsilon j)} (\log r)^{2 (\alpha-1)}
\\
&=&
c_\alpha\, (\log r)^{2(\alpha-1)} \, \sum_{j=1}^{\lceil 1/\varepsilon\rceil} \, r^{-(\alpha-1) (1- \varepsilon j)} \, \P \Big( \sum_{i=1}^k \ee^{S_i - S_k} \in (r^{\varepsilon(j-1)}, r^{\varepsilon j}]\Big). 
\end{eqnarray*}

 \noindent As for $B_{\eqref{en2}}$, we have by the time-reversal and the choice of $\eta$ that  $$
 \P \Big( \sum_{i=1}^k \ee^{S_i - S_k} \in (r^{\varepsilon(j-1)}, r^{\varepsilon j}]\Big) 
 \le \P \Big( \sum_{i=0}^{k-1} \ee^{- S_i} >r^{\varepsilon(j-1)} \Big)
 \le  r^{-\eta \,\varepsilon \, (j-1)}.
 $$

 \noindent  It follows that $$
 A_{\eqref{en2}}
 \le
 c_\alpha\, \lceil 1/\varepsilon\rceil\, (\log r)^{2(\alpha-1)} \,  \max_{1\le j \le \lfloor 1/\varepsilon\rfloor} r^{-(\alpha-1) (1- \varepsilon j) -\eta \,\varepsilon \, (j-1) } \le r^{-\min(\alpha-1, \eta) +O(\varepsilon)},
 $$
 
\noindent where  $O(\varepsilon)$ denotes some quantity bounded by $c\, \varepsilon$ with some positive constant $c$ depending on $\alpha, \eta$. Since $O(\varepsilon)$ can be chosen as small as desired, we assemble the above estimates  on $A_{\eqref{en2}}$ and on $B_{\eqref{en2}}$ and get~\eqref{pre1}. This completes the proof of Lemma~\ref{l:pre}. \hfill$\Box$.

 \medskip

We are now ready to give the proof of Proposition~\ref{p:main2}:

 \medskip
 {\noindent\it Proof of Proposition~\ref{p:main2}.}
If $\kappa \in (1, 2]$, for any $0< q < \kappa-1$,  we choose a small constant $0< \delta < \kappa-1- q$. If $\kappa \in (2, \infty]$, we choose $q:= 1+\frac\delta2$ with a small constant $\delta>0$.  For $| x| \ge K_\varepsilon$, we have $\ee^{-U(x)} \le \frac\varepsilon8$,  therefore we can apply Lemma~\ref{l:sum_iid} and obtain that 
\begin{eqnarray}
 &&P_\omega \Big\{ 
  \max_{ K_\varepsilon \le |x| \le n^q}
  L_{\overline T_n} (x) \ge \lceil \varepsilon n\rceil \Big\}
 \nonumber \\
 &\le&6 
  \sum_{ K_\varepsilon \le |x| \le n^q}
  \frac{n}{\sum_{z\in \, [ \! [ \varnothing, \, x]\! ]} \ee^{V(z)}} 
  \exp\Big( - \frac{\varepsilon n}{8} \frac{\ee^{U(x)}}{\sum_{z\in \, [ \! [\varnothing, \, x]\! ]} \ee^{V(z)}} \Big) 
  \nonumber \\
  &\le&
c_\varepsilon\,  \sum_{ K_\varepsilon \le |x| \le n^q}
   \ee^{-U(x)}\, 
  \exp\Big( - \frac{\varepsilon n}{16} \frac{\ee^{U(x)}}{\sum_{z\in \, [ \! [\varnothing, \, x]\! ]} \ee^{V(z)}} \Big) 
  , \label{kappa>2}
  \end{eqnarray}

\noindent by using the elementary inequality:   $6 \, r\, \ee^{- \frac{\varepsilon}{8} r} \le c_\varepsilon \ee^{- \frac{\varepsilon}{16} r}$ for any $r>0$. 
 Taking the expectation of the right-hand side of~\eqref{kappa>2}, we deduce from~\eqref{pre1} that for all large $n\ge n_0$, $$
\p \Big\{ 
  \max_{ K_\varepsilon \le |x| \le n^q}
  L_{\overline T_n} (x) \ge \lceil \varepsilon n\rceil \Big\}
\le
  n^{-\varrho},
$$

\noindent with some positive constant $\varrho>0$.  Consider $n_i:= \lceil \ee^{\frac{i}{\log i}}\rceil $ for all large $i\ge i_0$. An  application of Borel-Cantelli's lemma yields that  $$ \limsup_{i\to \infty} \frac1{n_i}\, 
 \,  \max_{  K_\varepsilon\le |x| \le (n_i)^{q }}
  L_{\overline T_{n_i}} (x) 
  \le 
  \varepsilon, \qquad \p\mbox{-a.s.},$$

\noindent which in view of the monotonicity and the fact that $\frac{n_i}{n_{i-1}} \to 1$ impliy~\eqref{limsupltn}.

To get~\eqref{favempty}, we remark that $\p$-a.s.\ for all large $j$, if $n$ is the integer such that $\overline T_{n-1} \le j < \overline T_n$, then ${\mathbb F}(j) \cap \{ x: K_\varepsilon \le |x| \le n^q\}= \emptyset$. Indeed, otherwise for any $x $ belonging to the intersection, $ L_{\overline T_n}(x) \ge L_j(x) \ge L_j(\varnothing)\ge n-1$ which would contradict~\eqref{limsupltn}.  By~\eqref{aslocaltime},  $\p$-a.s.\  for $\overline T_{n-1} \le j < \overline T_n$, we have $n \ge j^{\max(1/\kappa, 1/2)+o(1)}$, implying~\eqref{favempty}.  \hfill$\Box$

 \section{Proof of Theorem~\ref{t:tailmaxlocaltime}}\label{s:thm2}

 This section is devoted to the study of the tail of the maximum of (edge) local times. At first we present several estimates on the Markov chain $(Y_n)_{n\ge 0}$ introduced in~\eqref{defyn}.

  Let \begin{equation}\label{defsigma1} \sigma^+_1:= \inf\{n\ge 1: Y_n=1\}
 \end{equation}

 \noindent be the first return time of $(Y_n)_{n\ge 0}$ to $1$. We estimate the maximum of an excursion of $(Y_n)_{n\ge 0}$ in the following lemma:
 
 \begin{lemma}\label{l:maxYexcursion} Assume~\eqref{hyp1},~\eqref{kappa},~\eqref{non-lattice} and~\eqref{hyp2}. For any $\kappa \in (1, \infty)$, there are two numerical positive constants $c_2 > c_1$  such that for any $r \ge 1$, $$ 
 c_1 \, r^{-(\kappa-1)} 
 \le 
 \widehat \p_1\Big(\max_{1\le n \le \sigma^+_1} Y_n \ge r\Big)
 \le
 c_2\, r^{-(\kappa-1)} ,
 $$
 where we recall that under $ \widehat \p_1$, $Y_0=1$. 
 \end{lemma}

 {\noindent\it Proof of Lemma~\ref{l:maxYexcursion}.} The upper bound follows from the asymptotic behaviors of the invariant probability $(\pi_i)_{i\ge 1}$. Let  \begin{equation}\label{lyj}
  \ell_Y(j):= \sum_{n=1}^{\sigma_1^+} 1_{\{Y_n=j\}}, \qquad j\ge 1,
 \end{equation}
 
 \noindent be the local times of the Markov chain $(Y_n)_{n\ge 0}$ up to $\sigma_1^+$.  Then for any $r>1$,  $$ \widehat \p_1 \Big(\max_{1\le n \le \sigma^+_1} Y_n \ge r\Big)
 \le
 \widehat \p_1 \Big(  \sum_{j=r}^\infty \ell_Y(j) \ge 1\Big)
\le
 \sum_{j=r}^\infty \widehat \e_1 (\ell_Y(j)). 
$$
 
 \noindent 
Since $(\pi_i)_{i\ge 1}$  is the invariant measure for the Markov Chain $(Y_n)_{n\ge 0}$, we obtain  that 
 \begin{equation} \label{upplyj}
  \widehat \e_1 \big( \ell_Y(j) \big)
  =
 \frac{\pi_j}{\pi_1}
 \le
 c\, j^{-\kappa}, \qquad \forall \, j\ge 1,
 \end{equation}
 
 \noindent by applying~\eqref{pi1} to get the above inequality. The upper bound follows.
 
 For the lower bound, we use a representation of $(Y_n)_{n\ge 0}$ in terms of a branching process in random environment (BPRE) with immigration: recalling $P_{i, j}$ from~\eqref{pij}. For $0< a < 1$, let $(\xi_{i, n}(a))_{i, n \ge 1}$ be a family of i.i.d.\ geometric random variables such that $ \P(\xi_{1, 1}(a)= n)= (1-a) \, a^n$ for all $n\ge 1$, and independent of $(S_n)_{n\ge1}$. We observe that for any $i, j\ge 1$, $$
 P_{i, j} 
 =
 \P\Big( \xi_{1, 1}({\tt a_1}) + ...+ \xi_{i+1, 1}({\tt a_1})= j-1\Big),
 $$
 
 \noindent with ${\tt a_1}:= \frac{\ee^{-S_1}}{1+ \ee^{-S_1}}$.  Let ${\tt a_n}:= \frac{\ee^{-(S_n- S_{n-1})}}{1+ \ee^{-(S_n- S_{n-1})}}$ for any $n\ge 2$. Define a process $(\widehat Y_n)_{n\ge 0}$ by induction: $\widehat Y_0=0$ and $$ \widehat Y_n:=
 \sum_{k=1}^{\widehat Y_{n-1}+ 2} \xi_{k, n}({\tt a_n}), \qquad n\ge1.$$ 

\noindent Therefore the  law of $(\widehat Y_n)_{n\ge 0}$ (under $\P$) is exactly the law of $(Y_n-1)_{n\ge 0}$ under $ \widehat \p_1$, consequently $$ \widehat \p_1 \Big(\max_{1\le n \le \sigma^+_1} Y_n \ge r\Big)
=
 \P\Big( \max_{1\le n\le \widehat\sigma_0} \widehat Y_n \ge r-1\Big),$$

\noindent where $\widehat\sigma_0:= \inf \{n\ge 1: \widehat Y_n= 0\}$. Remark that in each generation, there are  two immigrants in $(\widehat Y_n)$, hence  $\max_{1\le n\le \widehat\sigma_0} \widehat Y_n $ is stochastically larger than $\max_{n\ge 0} Z_n$, where $Z$ is a BPRE started at $1$ (without immigration and in the same environment $({\tt a_n})_{n\ge 1}$). Conditioning on $({\tt a_n})_{n\ge 1}$, $Z_1$ has the mean $\frac{{\tt a_1}}{1-{\tt a_1}}= \ee^{-S_1}$. Observe that $\E\big( \ee^{-(\kappa-1) S_1}\big)=1$.  Applying Afanasyev~\cite{afanasyev} gives that for some positive constant $c$, 
% \footnote{In fact Afanasyev~\cite{afanasyev} dealt with the non-lattice case and obtained a true equivalence, however the proof presented there (in the equation (55), page 104) can easily be adopted to the lattice case.} 
$$
 \P\Big( \max_{n\ge 0} Z_n > r \Big) \, \sim \, c\,  r^{- (\kappa-1)}, \qquad r \to \infty.$$

\noindent It follows that $  \widehat \p_1 \big(\max_{1\le n \le \sigma^+_1} Y_n \ge r\big) \ge c_1 r^{-(\kappa-1)}$ for all $r\ge 1$. \hfill$\Box$

\medskip
The following lemma gives a uniform estimate on the tail of $ \max_{0 \le n \le \sigma_1^+} Y_n$ when $Y_0$ is an arbitrary integer: 

\begin{lemma}\label{l:Yexit}  Let $\kappa \in (1, \infty)$.  Assume~\eqref{hyp1},~\eqref{kappa} and~\eqref{hyp2}.    For any $0< \gamma< \kappa-1$, there exists some constant $c_\gamma>0$ such that for all $l, r \ge 1$, $$
 \widehat \p_l \Big( \max_{0 \le n \le \sigma_1^+} Y_n \ge r \Big)
=
 \widehat \p_l \Big( \tau_r < \sigma_1^+\Big)
\, \le \, 
c_\gamma\, \Big(\frac{l}{r}\Big)^{\gamma},$$
where as before $ \widehat \p_l$ means that the Markov chain $(Y_n)_{n\ge 0}$ starts at $l$, and $\tau_r:=\min\{n\ge 0: Y_n \ge r\}$ is the first time that $(Y_n)_{n\ge 0}$ exceeds the level $r$. 
\end{lemma}

{\noindent\it Proof of Lemma~\ref{l:Yexit}.} Fix $0< \gamma < \kappa-1$. Consider the function $f(i):= \frac{\Gamma(i+\gamma)}{\Gamma(i)}$ for $i\ge1$. We mention  that $f$ is  increasing and $$ f(i) \sim i^\gamma, \qquad i \to \infty.$$

 We claim the existence of some integer $i_0\ge1$ such that  \begin{equation}\label{pfi}
Pf(i) \le f(i), \qquad \forall \, i\ge i_0,
\end{equation}

\noindent where $P=(P_{i,j})_{i, j\ge1}$ denotes the matrix transition of the Markov chain $(Y_n)_{n\ge 0}$ under the probability measure $\widehat\p_l$.  The proof of~\eqref{pfi} is given in the Appendix.

 Denote by $\sigma_F:= \inf\{n\ge 0: X_n \in F\}$ where $F:= \{1, 2, ..., i_0\}$. Then by using~\eqref{pfi}, an application of the Markov property of $(Y_n)_{n\ge 0}$ says that $f(Y_{n\wedge \sigma_F})$ is a supermartingale. 

 By the optional stopping theorem, we obtain that  for any $r \ge l >i_0$, $$
 f(l)
 \ge
 \widehat \e_l f(Y_{\sigma_F\wedge \tau_r})
 \ge
f(r)\,  \widehat \p_l( \tau_r < \sigma_F) .
 $$

 \noindent It follows that $$  \widehat \p_l( \tau_r < \sigma_1^+)
 \le
 \widehat \p_l( \tau_r < \sigma_F) +  \widehat \p_l( \sigma_F < \tau_r < \sigma_1^+) 
  \le 
 \frac{f(l)}{f(r)} + \widehat \p_l( \sigma_F < \tau_r < \sigma_1^+) . $$

 \noindent To estimate $ \widehat \p_l( \sigma_F < \tau_r < \sigma_1^+) $, we apply the Markov property at $\sigma_F$: $$
  \widehat \p_l( \sigma_F < \tau_r < \sigma_1^+) 
 =
  \widehat \e_l \Big( 1_{\{ \sigma_F < \tau_r \wedge \sigma_1^+\}} \widehat \p_{Y_{\sigma_F}} ( \tau_r < \sigma_1^+)\Big)
  \le
  \max_{2\le i \le i_0} \widehat \p_i( \tau_r < \sigma_1^+).
 $$
  
  \noindent For any $i \in \{2, ..., i_0\}$ (and $r > i_0$), $$
  \widehat \p_i( \tau_r < \sigma_1^+) \le 
  \frac1{P_{1, i}} \, \widehat \p_1 ( \tau_r < \sigma_1^+)
  =
  \frac1{P_{1, i}} \, \widehat \p_1 \Big( \max_{1\le n \le \sigma_1^+} Y_n \ge r\Big) .$$
 
Applying the upper bound in Lemma~\ref{l:maxYexcursion} gives that there exists some constant $c_{i_0} >0$ such that for any $i \in \{1, ..., i_0\}$, $ \widehat \p_i( \tau_r < \sigma_1^+) \le c_{i_0}\, r^{-(\kappa-1)}$ for all $r\ge1$. 

It follows that for all $l, r\ge1$, $  \widehat \p_l( \tau_r < \sigma_1^+)
 \le
 \frac{f(l)}{f(r)} +  c_{i_0}\, r^{-(\kappa-1)}$, which yields  Lemma~\ref{l:Yexit} because $f(i) \asymp i^\gamma$ for all $i \ge 1$. \hfill $\Box$.

\medskip

Let $k\ge1$, recall the definition of $\ZZ_k$ from~\eqref{zzk}. According to the terminology in~\cite{BK04}, $\ZZ_k$ is an optional line for the multi-type Galton-Watson tree. For any $x \in \T^{(k)}$, recall that we denote by $x \le \ZZ_k$ if for all $\varnothing < y < x$, $\overline L_{\overline T_k}(y) \ge 2$. The following lemma controls the distance of the line $\ZZ_k$ from the root $\varnothing$ when $k$ is large:

\begin{lemma}\label{l:zzk}  Let $\kappa \in (1, \infty)$.   Assume~\eqref{hyp1},~\eqref{kappa} and~\eqref{hyp2}.    There is some positive constant $c_3$ such that $\p^*$-almost surely, for all large enough $k$, \begin{equation}\label{l:zzk1}
 \max_{x \in \ZZ_k} |x| \le c_3 \, \log k.
\end{equation}
\noindent Moreover, there exists $(a_k)_{k\ge 0}$ an increasing deterministic sequence, $a_k \to \infty$, such that $\p^*$-almost surely, for all large enough $k$, 
\begin{equation}\label{l:zzk2}
 \min_{x \in \ZZ_k} |x| \ge a_k.
 \end{equation} 

\end{lemma}

\begin{remark}
Under some additional integrability assumption, for instance if there exists some $\delta>0$ such that $\E (\sum_{i=1}^\nu A_i^{-\delta}) < \infty$, then we may take $a_k= c \log k$ for some positive constant $c$ in~\eqref{l:zzk2}.
\end{remark}

{\noindent\it Proof of Lemma~\ref{l:zzk}.} When $\kappa >2$,~\eqref{l:zzk1} follows immediately from  Lemma 7.2 in~\cite{AR15}. We mention that this Lemma 7.2 is also valid when $1< \kappa \le 2$. In fact, with the notations and the equality (5.6) therein, it is enough to remark that $P_{\T}(N^{(1)}_y\ge 1) ^2 \le P_{\T}(N^{(1)}_y\ge 1)^\beta$ with some $1< \beta < \kappa\wedge 2$. Since $P_{\T}(N^{(1)}_y\ge 1) \le \ee^{-V(y)}$ and $\E \sum_{|y|=\ell/2} \ee^{- \beta V(x)} $ decays exponentially fast in $\ell$,~\eqref{l:zzk1} follows from an application of the convergence part of Borel-Cantelli's lemma.

To get~\eqref{l:zzk2}, we first remark an elementary fact: for any    sequence $(\xi_k)_{k\ge 1}$ such that $\xi_k \to \infty$ $\p^*$-almost surely, there exists some increasing deterministic sequence $a_k \to \infty$ such that $\p^*$-almost surely, \begin{equation}\label{xkak} \xi_k \ge a_k,
\end{equation} for all large enough $k$. Indeed, defining $b_0:=1$ and $b_i:= \inf\{n> b_{i-1}: \p^*(\inf_{k\ge n} \xi_k < i) \le i^{-2}\}$ for any $i\ge 1$, the Borel-Cantelli lemma yields that $\p^*$-a.s., for all $\inf_{k\ge b_i}\xi_k \ge i$ for all large $i$.  We define  $a_k:= \inf\{i\ge 0: b_{i+1} \ge k\}$, which satisfies~\eqref{xkak}.

Now, observe that by hypothesis,   $\p^*$-almost surely, every generation $\T$ is finite. Hence by the definition of $\ZZ_k$, $\min_{x \in \ZZ_k} |x| \to_{k \to \infty} \infty$ (as the second equation of Fact~\ref{f:loic} ensures that $\p^*$-a.s., $\ZZ_k$ is non-empty for $k$ large enough). This completes the proof of Lemma~\ref{l:zzk}. \hfill$\Box$

%%%Proof of the lower bound with $\log n$ where we need to assume $\phi(-\delta)< \infty)$ for some $\delta>0$, this condition assumes the finiteness of the velocity in the BRW.
%For the lower bound, we first apply the Biggins-Hammersley-Kingman theorem on the first order of $\max_{|x|=l} V(x)$
%(\cite{H74, K75, B76}): There is some positive and finite constant $c_3>0$ such that \begin{equation}\label{llnbrw} \lim_{l\to \infty} \frac1{l} %\max_{|x|\le l} V(x) = c_3, \qquad \P^*\mbox{-a.s.}
%\end{equation}

%\noindent 
%Observe that for any $l\ge 1$, $$
%P_\omega\Big( \min_{x \in \ZZ_k} |x| < l \Big)
%\le
%P_\omega\Big( \exists |y|=l: \overline L_{\overline T_k}(y) =0\Big)
%\le
%\sum_{|y|=l} \Big( 1-  \frac{1}{\sum_{z\in \, [ \! [ \varnothing, \, y]\! ]} \ee^{V(z)}} \Big)^k, 
%$$

%\noindent by using~\eqref{valeur(p,a)2}.  Let $c'_3> c_3$ and consider $l= \frac{\log k}{2 c'_3}$.  It follows that 
%\begin{eqnarray*}
%&& \p\Big( \min_{x \in \ZZ_k} |x| < l, \max_{|x|\le l} V(x) \le c'_3\, l \Big)
%\\
%&\le&
%\E \Big( 1_{\{ \max_{|x|\le l} V(x) \le c'_3\, l\}} \sum_{|y|=l} \Big( 1-  \frac{1}{\sum_{z\in \, [ \! [ \varnothing, \, y]\! ]} \ee^{V(z)}} %\Big)^k \Big)
%\\
%&\le&
%\E \Big( 1_{\{ \max_{i \le l} S_i \le c'_3 \, l\}}  \ee^{S_l} \Big( 1-  \frac{1}{\sum_{0\le i \le l} \ee^{S_i}} \Big)^k \Big)
%\\
%&\le& \ee^{c'_3 l}\, \ee^{- \frac{k}{l}\, \ee^{-c'_3 l}},
%\end{eqnarray*}
%\noindent where in the above second inequality we have used the many-to-one formula~\eqref{many-to-one}.  Applying the Borel-Cantelli lemma gives the lower bound in  Lemma~\ref{l:zzk}. \hfill$\Box$ 

\medskip
The following result estimates the maximum of edge local times up to the optional line $\ZZ_1$ ($\ZZ_1$ being defined in~\eqref{zzk}): 

\begin{lemma}\label{l:bush}  Let $\kappa \in (1, \infty)$.  Assume~\eqref{hyp1},~\eqref{kappa},~\eqref{non-lattice} and~\eqref{hyp2}.   There exist some positive constants $c_4, c_5$ such that for any $r \ge 1$, 
$$
c_4 \, r^{-\kappa} 
\le
\p \Big( \max_{ x \le \ZZ_1} \overline L_{\overline T_1}(x) \ge r \Big) 
\le 
c_5\, r^{-\kappa}. 
$$
\end{lemma}

{\noindent\it Proof of Lemma~\ref{l:bush}.}  Let $r>1$. By considering the first generation $n$ such that $ \max_{|x|=n} \overline L_{\overline T_1}(x)\ge r$, we get that  \begin{eqnarray*}
\p\Big(\max_{ x \le \ZZ_1} \overline L_{\overline T_1}(x) \ge r \Big) 
&\le& 
\sum_{n=1}^\infty \, \e \Big( \sum_{|x|=n} 1_{\{ \overline L_{\overline T_1}(x) \ge r, \, \max_{1\le i < n} \overline L_{\overline T_1}(x_i)< r, \, \min_{1\le i < n} \overline L_{\overline T_1}(x_i) \ge 2\}}\Big)
\\
&=&
\sum_{n=1}^\infty \, \widehat \e_1 \left(\frac1{Y_n} 1_{\{Y_n \ge r, \, \max_{1\le i < n} Y_i < r, \, \min_{1\le i < n} Y_i \ge 2\}}\right)
,
\end{eqnarray*}

\noindent by using the many-to-one formula~\eqref{ar15} for the edge local times. It follows that \begin{eqnarray*}
\p\Big(\max_{ x \le \ZZ_1} \overline L_{\overline T_1}(x) \ge r \Big) 
&\le& \frac1{r}\, \sum_{n=1}^\infty \, \widehat \p_1\Big(Y_n \ge r, \, \max_{1\le i < n} Y_i < r, \min_{1\le i < n} Y_i \ge 2\Big) 
\\
&=&
\frac1{r}\, \widehat \p_1 \Big( \max_{1\le i < \sigma_1^+} Y_i \ge r\Big),
\end{eqnarray*}

\noindent by using the first return time $\sigma_1^+:= \min\{i\ge 1: Y_i=1\}$. By the upper bound of Lemma~\ref{l:maxYexcursion}, we get that $$
\p\Big(\max_{ x \le \ZZ_1} \overline L_{\overline T_1}(x) \ge r \Big) 
\le 
c_2\, r^{-\kappa}.$$

 To get the lower bound, we introduce  \begin{eqnarray*}
 Z&:=&
 \sum_{x \le \ZZ_1} 1_{\{\overline L_{\overline T_1}(x) \ge r, \, \max_{\varnothing < y < x} \overline L_{\overline T_1}(y)< r\}}
 \\
 &=&
 \sum_{n=1}^\infty \sum_{x: |x|=n} 1_{\{ \overline L_{\overline T_1}(x) \ge r, \, \max_{1\le i < n} \overline L_{\overline T_1}(x_i)< r, \, \min_{1\le i < n} \overline L_{\overline T_1}(x_i) \ge 2\}}, 
 \end{eqnarray*}
 
\noindent which is the cardinal of the set of vertices which are the first of their ancestry line to have their edge local time to overshoot $r$. Remark that $\{ Z \ge 1\} \subset \{ \max_{ x \le \ZZ_1} \overline L_{\overline T_1}(x) \ge r\}$. We choose and fix $p>1$ such that \begin{equation} \label{def-p}
\begin{cases}
 \kappa < p < \min (\alpha, 2) , \qquad & \mbox{if } 1< \kappa <2, \\
 p := 2, \qquad & \mbox{if } \kappa \ge 2, 
\end{cases}
\end{equation}

\noindent where $\alpha>\kappa$ is the constant in the assumption~\eqref{hyp2}. By the Paley-Zygmund inequality,  \begin{equation} \label{paley}
 \p\Big( Z \ge 1 \Big) \ge \frac{(\e (Z))^{p/(p-1)}}{(\e(Z^p))^{1/(p-1)}} . \end{equation}
 
First we estimate  $\e(Z)$. It follows from the many-to-one formula~\eqref{ar15} that   $$ \e (Z)
 =
 \sum_{n=1}^\infty  \widehat\e_1 \Big( \frac{1}{Y_n} 1_{\{\tau_r=n , \, \tau_r < \sigma^+_1\}}\Big)
 = \widehat\e_1 \Big( \frac{1}{Y_{\tau_r}} 1_{\{\tau_r < \sigma^+_1\}}\Big),$$

\noindent where as in Lemma~\ref{l:Yexit},  $\tau_r =\min\{n\ge 0: Y_n \ge r\}$ denotes the first time that $(Y_n)_{n\ge 0}$ exceeds $r$. Applying the upper bound in Lemma~\ref{l:maxYexcursion} gives $$ \e(Z)
\le
c_2 \, r^{- \kappa}.$$

On the other hand, by applying the lower bound in Lemma~\ref{l:maxYexcursion}, we see that for some constant $c$ large enough and for any $r\ge1$, $$ \widehat\p_1 \Big( \tau_r < \sigma_1^+ , Y_{\tau_r} \le cr \Big) 
=
 \widehat\p_1 \Big( \tau_r < \sigma_1^+  \Big) - \widehat\p_1 \Big( \tau_{ c r} < \sigma_1^+  \Big) 
\ge
c_1 r^{-\kappa} - c_2 (c r)^{-\kappa}
\ge
c_6 r^{-\kappa}.$$

Hence $\e(Z) \ge (cr)^{-1} \widehat\p_1\big( \tau_r < \sigma_1^+ , Y_{\tau_r} \le cr \big) \ge \frac{c_6}{c} \, r^{-\kappa}$. Consequently we have proved that \begin{equation}\label{moment1Z}
c_7 \, r^{-\kappa} \le \e(Z) \le c_2\, r^{-\kappa}, \qquad \forall r\ge1.
\end{equation}

Now we estimate the $p$-th moment of $Z$ by using the change of measure formula in Fact~\ref{f:measurechange}. At first, we recall the definition of $\T^{(1)}$ in \eqref{def-Tk} and    introduce the following optional line (in the sense of Biggins and Kyprianou~\cite{BK04}): $$ \L_r:=\Big\{ x \in \T^{(1)}: \overline L_{\overline T_1}(x) \ge r, \, \max_{\varnothing < y < x} \overline L_{\overline T_1}(y)< r\Big\}, \qquad r>1,$$

\noindent   the set of vertices which are the first on  the  ancestry line to have an  edge local time   exceeding  $r$. Recall the definition of $\ZZ_1$ from~\eqref{zzk}. We remark that for any $r>1$, $$ Z= \sum_{x \in \L_r} 1_{\{ x \le \ZZ_1\}}.$$

A standard argument in the studies of branching random walk, see for instance Biggins and Kyprianou~\cite{BK04}, shows that we may replace the set $\{|x|=n\}$ by an optional line $\L_r$, and the corresponding change of measure formula in Fact~\ref{f:measurechange} still holds (with obvious modifications). Write $x\sim \w_i$ iff ${\buildrel \leftarrow \over x}= \w_{i-1}$ and $x\neq \w_i$ for any $i\ge1$. Then we have \begin{eqnarray}
\e \Big( Z^p \Big)
&=&
\widehat\e_1\Big[ \frac{1}{Y_{\tau_r}} 1_{\{\w_{\tau_r} \le \ZZ_1\}} \, \Big(\sum_{x \in \L_r} 1_{\{ x \le \ZZ_1\}} \Big)^{p-1} \Big]
\nonumber \\
&=&
\widehat\e_1\Big[ \frac{1}{Y_{\tau_r}} 1_{\{\tau_r < \sigma^+_1\}} \, \Big( 1+ \sum_{i=1}^{\tau_r} \sum_{x\sim \w_i} \, Z^{(x)} \Big)^{p-1}\Big], \label{z2a}
\end{eqnarray}

\noindent where for any $x \sim \w_i$, $$Z^{(x)}:= \sum_{ y \in \L_r, y\ge x} 1_{\{ y \le \ZZ_1\}}.$$

Let ${\cal Y}: =\sigma\{Y_i, \w_i,  \overline L_{\overline T_1}(x), x \sim \w_i ,  i \le \tau_r\}$ be the $\sigma$-fields generated by the spine up to $\tau_r$. By the choice of $p$ in~\eqref{def-p}, $0< p-1 \le 1$, and it follows that  \begin{equation}\label{zpy} \widehat\e_1\Big[ \Big( 1+ \sum_{i=1}^{\tau_r} \sum_{x\sim \w_i} \, Z^{(x)} \Big)^{p-1} \, \big| {\cal Y} \Big] 
\le 
 1+ \sum_{i=1}^{\tau_r} \Big( \sum_{x\sim \w_i} \, \widehat \e_1  \big[Z^{(x)} \, |\, {\cal Y} \big] \Big)^{p-1}. 
 \end{equation}
Notice that by the branching property outside the spine $(\w_i)$, on the event $\{ \overline L_{\overline T_1}(x) = l \}$, $$ \widehat \e_1 \big[Z^{(x)} \, |\, {\cal Y} \big]
=
\e(Z_l),$$

\noindent where $Z_0:=0$ and $$Z_l:=
 \sum_{x \in \T^{(l)}} 1_{\{ \overline L_{\overline T_l}(x) \ge r, \, \max_{\varnothing < y < x} \overline L_{\overline T_l}(y)< r, \, \min_{\varnothing < y < x} \overline L_{\overline T_l}(y) \ge 2\}}, \qquad l\ge 1. $$

 By applying~\eqref{ar15} and then Lemma~\ref{l:Yexit} , we get that for any $ \beta \in [1, \kappa)$, there exists some constant $c_\beta>1$ such that for all $r \ge l\ge1 $, \begin{equation}\label{zelbeta} \e(Z_l) 
=
 \widehat\e_l \Big( \frac{l}{Y_{\tau_r}} 1_{\{\tau_r < \sigma^+_1\}}\Big)
\le
c_\beta \, \Big(\frac{l}{r}\Big)^{\beta}.
\end{equation}

\noindent 
Going back to~\eqref{z2a}, we deduce from~\eqref{zpy} and~\eqref{zelbeta} that  for any $1\le \beta < \kappa$, \begin{eqnarray}
\e(Z^p) 
&\le&
 \frac1r \widehat \p_1\Big(\tau_r < \sigma^+_1\Big)+ (c_\beta)^{p-1}\, \widehat\e_1\Big[ \frac{1}{Y_{\tau_r}} 1_{\{\tau_r < \sigma^+_1\}} \, \sum_{i=1}^{\tau_r} \Big( \sum_{x\sim \w_i} \, \big(\frac{\overline L_{\overline T_1}(x)}{r}\big)^\beta \Big)^{p-1} \Big]
\nonumber \\
&\le&
c_2 \, r^{-\kappa} + c_{\beta, p} \, r^{-1- \beta (p-1)}\, \widehat\e_1\Big[ 1_{\{\tau_r < \sigma^+_1\}} \, \sum_{i=1}^{\tau_r} \Big( \sum_{x\sim \w_i} \, \big(\overline L_{\overline T_1}(x)\big)^\beta\Big)^{p-1}\Big]
\nonumber \\
&\le&
c_2 \, r^{-\kappa} + c_{\beta, p} \, r^{-1- \beta (p-1)}\, \widehat\e_1\Big[ 1_{\{\tau_r < \sigma^+_1\}} \, \sum_{i=1}^{\tau_r} \Big( \sum_{x\sim \w_i} \, \overline L_{\overline T_1}(x) \Big)^{\beta(p-1)}\Big], \label{zp1}
\end{eqnarray}

\noindent
where we have used the upper bound of Lemma~\ref{l:maxYexcursion}  for the  second inequality.

 Now we choose (and fix) a constant $\beta$ such that \begin{equation} \label{def-beta}
\begin{cases}
\beta:=1 , \qquad & \mbox{if } 1< \kappa <2, \\
 \kappa-1 < \beta < \min(\alpha-1, \kappa) , \qquad & \mbox{if } \kappa \ge 2, 
\end{cases}
\end{equation}
Let us admit for the moment the existence of some positive constant $c $ such that \begin{equation} \label{ebetal} 
\widehat\e_1\Big[ \Big( \sum_{x\sim \w_i} \, \overline L_{\overline T_1}(x)\Big)^{\beta (p-1)} \, \big| \, \overline L_{\overline T_1}(\w_{i-1})= l\Big]\le c\, l^{\beta (p-1)}, \qquad \forall l\ge1.
\end{equation}

 \noindent Therefore we deduce from~\eqref{zp1} that $$
 \e(Z^p) 
 \le
c_2 \, r^{-\kappa} +  c'\, r^{-1- \beta (p-1)}\, \widehat\e_1\Big( 1_{\{\tau_r < \sigma^+_1\}} \, \sum_{i=1}^{\tau_r} (Y_{i-1})^{\beta(p-1)}\Big).$$

 \noindent By using the local time process $\ell_Y(\cdot)$ of $(Y_n)_{n\ge 0}$ defined in~\eqref{lyj}, we get  that (recalling that $\pi_j \asymp j^{-\kappa}$ and noticing that $\beta(p-1) > \kappa-1$, by the choices of $p$  in \eqref{def-p} and of $\beta$ in  \eqref{def-beta}) $$
 \widehat\e_1\Big( 1_{\{\tau_r < \sigma^+_1\}} \, \sum_{i=1}^{\tau_r} (Y_{i-1})^{\beta (p-1)}\Big)
 \le
 \sum_{j=1}^r \, j^{\beta (p-1)}\, \widehat\e_1(\ell_Y(j))
 =
 \sum_{j=1}^r \, j^{\beta (p-1)}\, \frac{\pi_j}{\pi_1}
 \le
 c\, r^{\beta(p-1)-\kappa+1}.
 $$
  
 \noindent 
 Consequently $$\e(Z^p) \le c''\, r^{-\kappa},$$

 \noindent which in view of~\eqref{paley} implies that $$
 \p\Big( \max_{ x \le \ZZ_1} \overline L_{\overline T_1}(x) \ge r \Big) 
 \ge 
 \p(Z\ge 1) 
 \ge
  c_8\, r^{-\kappa}, \qquad r\ge1.$$
 
 \noindent This proves the lower bound of  Lemma~\ref{l:bush}.

 It remains to check~\eqref{ebetal}. According to Fact~\ref{f:measurechange}, the expectation term in the left-hand-side of~\eqref{ebetal} is equal to \begin{eqnarray*}
 \widehat\e_l\Big( \sum_{|x|=1: x \sim \w_1} \, \overline L_{\overline T_l}(x) \Big)^{\beta (p-1)}
 &=&
 \e \Big[ \frac1{l}\, \sum_{|x|=1} \overline L_{\overline T_l}(x) \, \Big( \sum_{y: |y|=1, y \neq x} \overline L_{\overline T_l}(y)\Big)^{\beta (p-1)} \Big] 
 \\
 &\le&
  \frac1{l}\, \e \Big( \sum_{|x|=1} \overline L_{\overline T_l}(x) \Big)^{1+\beta (p-1)}.
 \end{eqnarray*}

 \noindent We shall use an application of Hölder's inequality: $(\sum_{i=1}^l t_i)^{1+\beta (p-1)} \le l^{\beta (p-1)}\, \sum_{i=1}^l t_i^{1+\beta (p-1)}$ for any $t_i\ge0$. Observe that under $\p$, $\{ \overline L_{\overline T_l}(\cdot) - \overline L_{\overline T_{l-1}}(\cdot)\}_{l\ge1}$ are identically distributed (but not independent). It follows that  $$ \frac1{l}\, \e \Big( \sum_{|x|=1} \overline L_{\overline T_l}(x) \Big)^{1+\beta (p-1)}
 \le
 l^{\beta (p-1)} \, \e \Big( \sum_{|x|=1} \overline L_{\overline T_1}(x) \Big)^{1+\beta (p-1)} = l^{\beta (p-1)} \e \Big( \Theta(\varnothing)\Big)^{1+\beta (p-1)}, $$

 \noindent by using the notation in~\eqref{Theta}. Then to check~\eqref{ebetal} it suffices to prove that $ \e \big( \Theta(\varnothing)\big)^{1+\beta (p-1)} < \infty.$ By~\eqref{lawThetax} (with $x=\varnothing$ and $k=1$ there), we get that $$ \e \big( \Theta(\varnothing)\big)^{1+\beta (p-1)} 
 =
 \E \sum_{n=1}^\infty n^{1+\beta (p-1)} \, \big(\frac{\sum_{|y|=1} A(y)}{1+\sum_{|y|=1} A(y)}\big)^{n} \, \frac1{1+\sum_{|y|=1} A(y)} .
 $$
 
 \noindent Elementary computations say that $\sum_{n=1}^\infty n^{1+\beta (p-1)} \, \big(\frac{t}{1+t}\big)^{n}  \le c_\beta \, (1+t)^{2+\beta (p-1)} $ for any $t>0$, it follows that $$ \e \big( \Theta(\varnothing)\big)^{1+\beta (p-1)} 
 \le
 c_\beta \, \E \Big( 1+ \sum_{|y|=1} A(y) \Big)^{1+\beta (p-1)},$$

 \noindent which is finite according to the assumption~\eqref{hyp2} (recalling that $1+ \beta (p-1)< \alpha$). This implies~\eqref{ebetal} and completes the proof of Lemma~\ref{l:bush}.  \hfill$\Box$

 \medskip
  Let us consider now the maximum of the (site) local times $ L_{\overline T_1}(x)$ instead of the edge local times $\overline L_{\overline T_1}(x)$:

  \begin{lemma}\label{l:bush2} Assume~\eqref{hyp1},~\eqref{kappa},~\eqref{non-lattice} and~\eqref{hyp2}.  There exist some positive constants $c_9, c_{10}$ such that for any $r \ge 1$, 
$$
c_9 \, r^{-\kappa} 
\le
\p \Big( \max_{ x \le \ZZ_1}  L_{\overline T_1}(x) \ge r \Big) 
\le 
c_{10} \, r^{-\kappa}. 
$$
\end{lemma}

{\noindent\it Proof of Lemma~\ref{l:bush2}.} In view of the lower bound in  Lemma~\ref{l:bush}, it is enough to prove the upper bound in Lemma~\ref{l:bush2}.

Recall~\eqref{siteedge} and~\eqref{Theta}.  
By applying Lemma~\ref{l:sum_iid} to~\eqref{lawThetax}, we get  that for any $r>1$, $k\ge1$, if $ \sum_{y: {\buildrel \leftarrow \over y}= x}  A(y) \le \frac18 \frac{r}{k}$, then \begin{equation}\label{Thetaxzk}
 P_\omega\Big( \Theta(x) \ge r \big| \overline L_{\overline T_1}(x) = k\Big)
 \le
 6 \, k \, \ee^{- \frac{r}{8 (1+ \sum_{y: {\buildrel \leftarrow \over y}= x}  A(y) ) }}.
 \end{equation} 
 
  \noindent 
 It follows that \begin{eqnarray}
&&P_\omega \Big( \max_{ x \le \ZZ_1} \Theta(x) \ge r\Big)
 \nonumber \\
&\le&
\sum_{n=0}^\infty \sum_{|x|=n} E_\omega \Big[ 1_{\{x \le \ZZ_1\}} \Big( 6 \, \overline L_{\overline T_1}(x) \, \ee^{- \frac{r}{8 (1+ \sum_{y: {\buildrel \leftarrow \over y}= x}  A(y) ) }} + 1_{\{ \overline L_{\overline T_1}(x) > \frac{r}{8\sum_{y: {\buildrel \leftarrow \over y}= x}  A(y)}\}} \Big) \Big] . \label{theta3}
\end{eqnarray}

\noindent Observe that the event $\{x \le \ZZ_1\}$ and the variable $\overline L_{\overline T_1}(x)$ only involve  those environments which are independent of $\sum_{y: {\buildrel \leftarrow \over y}= x}  A(y)$.  Moreover $\sum_{y: {\buildrel \leftarrow \over y}= x}  A(y)$ is distributed as $\sum_{i=1}^\nu A_i$. By taking the expectation with respect to the environment, we deduce from~\eqref{theta3} that \begin{eqnarray*}
&& \p \Big( \max_{ x \le \ZZ_1} \Theta(x) \ge r\Big)
 \nonumber \\
&\le&
\int \P(\sum_{i=1}^\nu A_i \in \d t) \, \sum_{n=0}^\infty \e \Big[ \sum_{|x|=n} 1_{\{x \le \ZZ_1\}} \Big( 6 \, \overline L_{\overline T_1}({\buildrel\leftarrow \over x}) \, \ee^{- \frac{r}{8 (1+ t ) }} + 1_{\{ \overline L_{\overline T_1}({\buildrel\leftarrow \over x}) > \frac{r}{8t}\}} \Big)\Big] 
\\
&=&
\int \P(\sum_{i=1}^\nu A_i \in \d t) \, \sum_{n=0}^\infty \widehat \e_1 \Big[ \frac{1_{\{ n \le \sigma_1^+\}}}{Y_n} \Big( 6 Y_n \, \ee^{- \frac{r}{8 (1+ t ) }} + 1_{\{ Y_n > \frac{r}{8t}\}} \Big)\Big] ,
\end{eqnarray*}

\noindent by using the many-to-one formula~\eqref{ar15} and the notation $\sigma_1^+$ in~\eqref{defsigma1}. Since the Markov Chain $(Y_n)_{n\ge 0}$ is positive recurrent, $\widehat\e_1(\sigma_1^+) < \infty$. This together with~\eqref{lyj} and~\eqref{upplyj} yield that \begin{eqnarray*}
 \sum_{n=0}^\infty \widehat \e_1 \Big[ \frac{1_{\{ n \le \sigma_1\}}}{Y_n} \Big( 6 Y_n \, \ee^{- \frac{r}{8 (1+ t ) }} + 1_{\{ Y_n > \frac{r}{8t}\}} \Big)\Big]
&\le&
c\, \, \ee^{- \frac{r}{8 (1+ t ) }} + \sum_{j> r/(8t)} \frac1{j} \widehat \e_1(\ell_Y(j)) 
\\
&\le&
c\, \, \ee^{- \frac{r}{8 (1+ t ) }} + c \, r^{-\kappa} \, t^\kappa.
\end{eqnarray*}

\noindent Hence $$ \p \Big( \max_{ x \le \ZZ_1} \Theta(x) \ge r\Big)
\le
c \, \E \ee^{- \frac{r}{8 (1+ \sum_{i=1}^\nu A_i) }} + c \, r^{-\kappa} \, \E \Big( \sum_{i=1}^\nu A_i\Big)^\kappa. $$

\noindent 
If we introduce an independent standard exponential variable ${\bf e}$, then $ \E \ee^{- \frac{r}{8 (1+ \sum_{i=1}^\nu A_i) }} = \P\big( (1+ \sum_{i=1}^\nu A_i) > \frac{r}{8{\bf e}}\big) \le c'\, \E \big(\frac{r}{8{\bf e}}\big)^{-\kappa}= c'\, r^{-\kappa}\, \E(8 {\bf e})^\kappa$, with $c':= \E (1+ \sum_{i=1}^\nu A_i)^\kappa< \infty$ by assumption. Then we get that $\p \big( \max_{ x \le \ZZ_1} \Theta(x) \ge r\big) \le c'' \, r^{-\kappa}$, which in view of~\eqref{siteedge} and  Lemma~\ref{l:bush} complete the proof of Lemma~\ref{l:bush2}.  \hfill$\Box$

\medskip

 Recall the definition of $\ZZ_1$ in~\eqref{zzk}. We are interested in the number of vertices in this optional line $\ZZ_1$.

 \begin{fact}[\cite{AR15},\cite{loic16}] \label{f:loic} Assume~\eqref{hyp1},~\eqref{kappa},~\eqref{non-lattice} and~\eqref{hyp2}.  For any $\kappa \in (1, 2]$, we have 
  \begin{equation}
 \label{loic}  \p\Big( \#\ZZ_1 > r \Big) \, \asymp\, r^{- \kappa}, \qquad \forall \, r\ge1, 
 \end{equation}
 and if $\kappa>2$, then
 \begin{equation}\label{loic2}
 \e\Big[\big(\#\ZZ_1\big)^2\Big]<\infty. 
 \end{equation}
  Moreover if $\kappa\in(1, 2]$, as $l\to \infty$, $$ \frac{\#\ZZ_l}{l}\quad {\buildrel \mbox{\tiny $L^p(\p^*)$} \over \longrightarrow}\, W_\infty, $$
  for any $1\le p < \kappa$, where $W_\infty$ denotes  the limit of the additive martingale $W_n:= \sum_{|x|=n} \ee^{-V(x)}$, which is positive $\p^*$-a.s.
 \end{fact}

 Now we are ready to prove Theorem~\ref{t:tailmaxlocaltime}:

 {\noindent\it Proof of Theorem~\ref{t:tailmaxlocaltime}.} We prove first the estimates on the edge local times $ \max_{x\in \T} \overline L_{\overline T_1}(x)$. By the strong Markov property, we get the following identity in law under the annealed probability measure $\p$:
 $$
  \max_{x\in \T} \overline L_{\overline T_1}(x) 
  \law
  \max\Big( \max_{x\le \ZZ_1} \overline L_{\overline T_1}(x) , \, \max_{1\le i \le \#\ZZ_1} \overline L^{*, i}\Big),$$
  
 \noindent where $\overline L^{*, i}, i\ge1$, are i.i.d.\ copies of $ \max_{x\in \T} \overline L_{\overline T_1}(x)$, independent of $(\max_{x\le \ZZ_1} \overline L_{\overline T_1}(x), \#\ZZ_1)$. Write $a_r:= \p(\max_{x\in \T} \overline L_{\overline T_1}(x) > r)$ for all $r\ge0$ and $\ell^*:=  \max_{x\le \ZZ_1} \overline L_{\overline T_1}(x)$ for notational brevity. We get from the above identity in law that for all $r>0$, $$ 
 a_r= 1- \e \Big( 1_{\{\ell^*\le r\}} (1- a_r)^{\#\ZZ_1}\Big)= 1- \e \Big((1-a_r)^{\#\ZZ_1}\Big) + \e \Big( 1_{\{\ell^*> r\}} (1- a_r)^{\#\ZZ_1}\Big).$$
 
 \noindent 
Define a function $f$ by $$ f(\varepsilon):= \varepsilon-1+  \e \Big((1-\varepsilon)^{\#\ZZ_1}\Big), \qquad 0< \varepsilon <1.$$

\noindent Then we get 
 \begin{equation}\label{Fr1}
   \e \Big( 1_{\{\ell^*> r\}} (1-  a_r )^{\#\ZZ_1}\Big) = f(a_r). \end{equation}

 \noindent Note that $1-  \e \big((1-\varepsilon)^{\#\ZZ_1}\big) =   \varepsilon \sum_{j=1}^\infty (1-\varepsilon)^{j-1}\p\big( \#\ZZ_1\ge j\big)$.  Note that  $\sum_{j=1}^\infty  \p\big( \#\ZZ_1\ge j\big)=\e(\#\ZZ_1)=1$. In fact, by the definition of $\ZZ_1$ in  \eqref{zzk} and by applying  the many-to-one formula \eqref{ar15},  
 \begin{eqnarray*} \e(\#\ZZ_1)
 &=& \sum_{n=1}^\infty \e\Big(\sum_{|x|=n}1_{\{ \overline L_{\overline T_1}(x_1) \ge 2, ..., \overline L_{\overline T_1}(x_{n-1})\ge 2, \overline L_{\overline T_1}(x_n)=1\}}\Big)
 \\
 &=&
   \sum_{n=1}^\infty  \widehat \e_1(\frac{1}{Y_n}  1_{\{ Y_1\ge 2, ..., Y_{n-1} \ge 2, Y_n =1\}})
   \\
   &=&  \sum_{n=1}^\infty  \widehat \p_1(  Y_1\ge 2, ..., Y_{n-1} \ge 2, Y_n =1)
 =1.
   \end{eqnarray*}
 
\noindent It follows          $$ f(\varepsilon) = \varepsilon   \, \sum_{j=1}^\infty (1-(1- \varepsilon)^{j-1}) \p\big( \#\ZZ_1\ge j\big), \qquad 0< \varepsilon< 1.$$

 \noindent  Based on~\eqref{loic} and~\eqref{loic2}, it is elementary to check that\footnote{The cases $1 < \kappa \le 2$ follow from an application of Tauberian theorem. Let us   give the details when  $\kappa>2$: Using the fact that  $1-(1- \varepsilon)^{j-1} \le \varepsilon (j-1)$, we get that $f(\varepsilon) \le \varepsilon^2 \e \sum_{j=1}^\infty (j-1) 1_{\{ j \le\#\ZZ_1\}} \le \varepsilon^2 \, \e[(\#\ZZ_1)^2]$,  implying  the upper bound on $f(\varepsilon)$. For the lower bound, we use the inequality that $1-(1- \varepsilon)^{j-1} \ge  \varepsilon (j-1)  (1- \varepsilon)^{j-1}$ and get that $\frac{f(\varepsilon)}{\varepsilon^2}  \ge    \e \sum_{j=1}^\infty (j-1) (1-\varepsilon)^{j-1}1_{\{ j \le\#\ZZ_1\}} \to   \e \sum_{j=1}^\infty (j-1)  1_{\{ j \le\#\ZZ_1\}}$ as $\varepsilon\to0$ by the monotone convergence,  and the lower bound on $f(\varepsilon)$ follows.}  as $\varepsilon \to 0$, \begin{equation}\label{Fr}
f(\varepsilon)\,
 \asymp \, \begin{cases}
 \varepsilon^\kappa , \qquad &\mbox{if } 1< \kappa<2,\\
\varepsilon^2 \log\frac1{\varepsilon}, \qquad &\mbox{if } \kappa=2,\\
 \varepsilon^2, \qquad &\mbox{if } 2< \kappa<\infty.\\
 \end{cases}
\end{equation}

For the upper bound of $ \p(\max_{x\in \T} \overline L_{\overline T_1}(x) > r)\equiv a_r$, we use  the equality~\eqref{Fr1} to see that $f(a_r) \le \p\big( \ell^* > r\big) \le c_5 \, r^{-\kappa}$ by
 Lemma~\ref{l:bush},  and the corresponding upper bound for $ a_r$ follows immediately from  \eqref{Fr}.

To get the lower bound of $ a_r$, we remark that for $r$ large enough (such that $a_r\le 1/2$), $$ f(a_r) = \e \big( 1_{\{\ell^*> r\}} (1-  a_r )^{\#\ZZ_1}\big)  \ge c\, \p\Big(\ell^* >r, \#\ZZ_1 \le \frac{1}{ a_r}\Big),$$

\noindent with $c:= \inf_{0<\varepsilon<\frac12} (1-\varepsilon)^{1/\varepsilon}>0$. It follows that $$  \p\Big( \ell^*> r \Big) 
\le
\frac1{c}  f(a_r)  + \p\Big( \#\ZZ_1 > \frac{1}{a_r}\Big)
\le 
\begin{cases}
 c'\, (a_r)^\kappa , \qquad &\mbox{if } 1< \kappa<2,\\
c'\,  (a_r)^2 \log\frac1{a_r}, \qquad &\mbox{if } \kappa=2,\\
c'\,   (a_r)^2, \qquad &\mbox{if } 2< \kappa<\infty,\\
 \end{cases} $$

\noindent  by using \eqref{Fr},  \eqref{loic} and  \eqref{loic2}.   By    Lemma~\ref{l:bush}, $   \p\big( \ell^*> r \big) \ge c_1 r^{-\kappa}$, which  gives the lower bound of $ a_r$. 

To deal with the local times instead of the edge-local times  $ \max_{x\in \T}  L_{\overline T_1}(x)$, we remark that again by the Markov property, under the annealed probability measure $\p$, 
 $$
  \max_{x\in \T}  L_{\overline T_1}(x) 
  \law
  \max\Big( \max_{x\le \ZZ_1}  L_{\overline T_1}(x) , \, \max_{1\le i \le \#\ZZ_1} \widehat L^{*, i}\Big),$$
  
 \noindent where $\widehat L^{*, i}, i\ge1$, are i.i.d.\ copies of $ \max_{x\in \T}  L_{\overline T_1}(x)$, independent of $(\max_{x\le \ZZ_1}  L_{\overline T_1}(x), \#\ZZ_1)$.

The rest of the proof goes exactly in the same way as that for the edge-local times, by applying Lemma~\ref{l:bush2} instead of Lemma~\ref{l:bush}.   \hfill$\Box$

\medskip

\section{Proof of Theorem~\ref{t:main}: Case $1< \kappa \le 2$}\label{s:kappale2}

%The case when $\kappa\in (1, 2]$ is more delicate. The main step is the integrability of the maximum of the (site) local times considered at $\overline T_1$. At first we study  the maximum of the edge local times by using the many-to-one formula~\eqref{ar15}. 
 
 \subsection{Proofs of~\eqref{t:kappa<2} and~\eqref{t:kappa<2b}, case $1< \kappa <2$}

\begin{proposition}\label{p:kappa<2} Let $1< \kappa< 2$. Assume~\eqref{hyp1},~\eqref{kappa},~\eqref{non-lattice} and~\eqref{hyp2}. Then $\p^*$-almost surely, $$ \limsup_{n \to \infty} \frac1n\, \max_{x \in \T} \overline L_{\overline T_n}(x) =\infty.$$
\end{proposition}

We shall use several times the following  fact in the proof of Proposition~\ref{p:kappa<2}:

\begin{fact}[\cite{shiryayev}, Chapter VII, Corollary 1] \label{f:xi} Let $(\xi_i)_{i\ge0}$ be a sequence of random variables adapted to some filtration $({\cal G}_i)_{i\ge0}$. Suppose that $\p$-almost surely,   $0\le \xi_i \le 1$ for any $i\ge 1$.  Then  $$
 \Big\{ \sum_{i=1}^\infty \xi_i = \infty \Big\}
=
\Big\{ \sum_{i=1}^\infty \e\big( \xi_i \big | {\cal G}_{i-1}\big)= \infty \Big\}, \qquad \p\mbox{-a.s.}$$
\end{fact}

{\noindent\it Proof of Proposition~\ref{p:kappa<2}.}  
Let $K $ be a large constant and $\varepsilon >0$ be small. Choose an increasing subsequence $n_j \to \infty$ as $j \to \infty$ such that 
\begin{equation}
a_{n_j} > n_{j-1},
\label{defnj}
\end{equation}

\noindent where the (deterministic) sequence $(a_k)$ is as in~\eqref{l:zzk2}. Recall the definition of $\ZZ_k$ for $k\ge 1$ from~\eqref{zzk}. Define for any $j\ge1$, %$$ B_j:= \Big\{ \max_{x \in \T} L^x_{\overline T_{n_j}} \ge K \, n_j,\, \max_{0\le i \le \overline T_{n_{j-1}}} | X_i| < \min_{x \in \ZZ_{n_j}} |x| , \, \# \ZZ_{n_j} \ge \varepsilon n_j \Big\} ,$$

 $$ B_j:= 
\Big\{ \max_{x \in \T} \overline L_{\overline T_{n_j}}(x) \ge K \, n_j,\, \max_{0\le i \le \overline T_{n_{j-1}}} | X_i| < \min_{x \in \ZZ_{n_j}} |x| , \, \# \ZZ_{n_j} \ge \varepsilon n_j 
\Big\} .$$

\definecolor{aqaqaq}{rgb}{0.6274509803921569,0.6274509803921569,0.6274509803921569}
\definecolor{zzttqq}{rgb}{0.6,0.2,0.0}
\begin{figure}[H]
\definecolor{cqcqcq}{rgb}{0.7529411764705882,0.7529411764705882,0.7529411764705882}
\definecolor{zzttqq}{rgb}{0.6,0.2,0.0}
\begin{tikzpicture}[line cap=round,line join=round,>=triangle 45,x=2cm,y=2cm]
\clip(6,-0.1) rectangle (14.0,3.8);
\fill[line width=0.0pt,color=zzttqq,fill=zzttqq,fill opacity=0.1] (8.9,0.0) -- (8.5,0.8) -- (8.2,0.9) -- (7.7,1.4) -- (7.1,1.6) -- (6.9,1.8) -- (6.9,2.2) -- (6.7,2.6) -- (6.3,3.0) -- (6.4,3.1) -- (6.5,3.0) -- (6.7,2.9) -- (6.9,2.8) -- (6.9,2.6) -- (7.0,2.4) -- (7.1,2.3) -- (7.2,2.5) -- (7.3,2.5) -- (7.3,2.6) -- (7.4,2.7) -- (7.5,2.6) -- (7.3,3.0) -- (7.1,3.2) -- (7.0,3.3) -- (8.7,3.3) -- (8.7,3.2) -- (8.9,3.0) -- (8.7,2.8) -- (8.9,2.4) -- (9.3,2.2) -- (9.3,1.6) -- (9.5,1.6) -- (9.5,1.9) -- (9.7,1.9) -- (9.7,2.0) -- (9.9,2.1) -- (9.9,2.2) -- (10.0,2.3) -- (10.1,2.2) -- (10.1,2.0) -- (10.3,1.8) -- (10.2,1.5) -- (10.3,1.3) -- (10.4,1.5) -- (10.5,1.5) -- (10.5,1.8) -- (10.3,2.0) -- (10.4,2.1) -- (10.5,2.0) -- (10.6,2.0) -- (10.7,2.2) -- (10.500000000000004,2.6000000000000085) -- (10.6,2.7) -- (10.7,2.8) -- (10.8,2.9) -- (10.7,3.0) -- (10.8,3.1) -- (11.1,2.8) -- (11.1,2.4) -- (11.3,2.2) -- (11.5,2.4) -- (11.5,2.6) -- (11.6,2.7) -- (11.7,2.6) -- (11.9,2.8) -- (11.9,3.0) -- (12.1,3.2) -- (12.1,3.3) -- (12.3,3.3) -- (12.3,3.1) -- (12.1,2.8) -- (11.7,2.4) -- (11.9,2.2) -- (11.3,1.6) -- (10.7,1.4) -- (10.5,1.2) -- (9.9,0.8) -- (9.1,0.0) -- (9.1,-0.1) -- (8.9,-0.1) -- cycle;
\fill[line width=0.0pt,color=zzttqq,fill=zzttqq,fill opacity=0.25] (10.7,0.1) -- (10.8,0.1) -- (10.8,0.0) -- (10.7,0.0) -- cycle;
\fill[line width=0.0pt,color=zzttqq,fill=zzttqq,fill opacity=0.1] (10.7,0.3) -- (10.8,0.3) -- (10.8,0.2) -- (10.7,0.2) -- cycle;
\fill[line width=0.0pt,color=zzttqq,fill=zzttqq,fill opacity=0.25] (9.1,-0.1) -- (8.9,-0.1) -- (8.9,0.0) -- (8.5,0.8) -- (8.7,0.8) -- (8.7,0.9) -- (8.95,0.9) -- (8.95,1.05) -- (9.05,1.05) -- (9.05,0.9) -- (9.2,0.9) -- (9.45,0.9) -- (9.45,0.85) -- (9.65,0.85) -- (9.65,0.75) -- (9.55,0.65) -- (9.65,0.55) -- (9.1,0.0) -- cycle;
\draw [line width=0.2pt] (9.0,0.0)-- (8.8,0.4);
\draw [line width=0.2pt] (9.0,0.0)-- (9.4,0.4);
\draw [line width=0.2pt] (9.4,0.4)-- (9.2,1.0);
\draw [line width=0.2pt] (9.4,0.4)-- (9.8,0.8);
\draw [line width=0.2pt] (8.8,0.4)-- (9.0,0.8);
\draw [line width=0.2pt] (8.8,0.4)-- (8.6,0.8);
\draw [line width=0.2pt] (8.6,0.8)-- (8.2,1.0);
\draw [line width=0.2pt] (8.6,0.8)-- (8.8,1.2);
\draw [line width=0.2pt] (9.2,1.0)-- (9.4,1.4);
\draw [line width=0.2pt] (9.8,0.8)-- (9.6,1.2);
\draw [line width=0.2pt] (9.8,0.8)-- (10.0,1.2);
\draw [line width=0.2pt] (9.8,0.8)-- (10.4,1.2);
\draw [line width=0.2pt] (8.6,0.8)-- (8.4,1.2);
\draw [line width=0.2pt] (8.2,1.0)-- (7.8,1.4);
\draw [line width=0.2pt] (8.8,1.2)-- (9.2,1.6);
\draw [line width=0.2pt] (8.4,1.2)-- (8.8,1.6);
\draw [line width=0.2pt] (8.4,1.2)-- (8.2,1.6);
\draw [line width=0.2pt] (7.8,1.4)-- (8.0,1.8);
\draw [line width=0.2pt] (7.8,1.4)-- (7.4,1.8);
\draw [line width=0.2pt] (9.4,1.4)-- (10.0,1.6);
\draw [line width=0.2pt] (9.4,1.4)-- (9.4,1.8);
\draw [line width=0.2pt] (10.0,1.2)-- (10.4,1.6);
\draw [line width=0.2pt] (10.4,1.2)-- (10.6,1.4);
\draw [line width=0.2pt] (10.0,1.6)-- (9.8,2.0);
\draw [line width=0.2pt] (10.4,1.6)-- (10.2,2.0);
\draw [line width=0.2pt] (10.0,1.6)-- (10.2,1.8);
\draw [line width=0.2pt] (8.8,1.6)-- (8.6,1.8);
\draw [line width=0.2pt] (8.8,1.6)-- (9.0,1.8);
\draw [line width=0.2pt] (8.6,1.8)-- (8.4,2.0);
\draw [line width=0.2pt] (9.0,1.8)-- (9.2,2.0);
\draw [line width=0.2pt] (9.0,1.8)-- (8.8,2.0);
\draw [line width=0.2pt] (8.6,1.8)-- (8.6,2.2);
\draw [line width=0.2pt] (8.0,1.8)-- (7.8,2.2);
\draw [line width=0.2pt] (7.4,1.8)-- (7.6,2.0);
\draw [line width=0.2pt] (7.4,1.8)-- (7.2,2.0);
\draw [line width=0.2pt] (8.2,1.6)-- (8.2,2.0);
\draw [line width=0.2pt] (10.6,1.4)-- (11.0,1.8);
\draw [line width=0.2pt] (10.6,1.4)-- (10.6,1.8);
\draw [line width=0.2pt] (10.6,1.4)-- (11.2,1.6);
\draw [line width=0.2pt] (7.8,1.4)-- (7.2,1.6);
\draw [line width=0.2pt] (7.2,1.6)-- (7.0,1.8);
\draw [line width=0.2pt] (7.2,1.6)-- (6.6,1.8);
\draw [line width=0.2pt] (6.6,1.8)-- (6.8,2.2);
\draw [line width=0.2pt] (6.6,1.8)-- (6.6,2.2);
\draw [line width=0.2pt] (7.0,1.8)-- (7.0,2.2);
\draw [line width=0.2pt] (7.8,2.2)-- (7.6,2.6);
\draw [line width=0.2pt] (7.8,2.2)-- (8.0,2.6);
\draw [line width=0.2pt] (8.6,2.2)-- (8.4,2.6);
\draw [line width=0.2pt] (8.2,1.6)-- (8.4,1.8);
\draw [line width=0.2pt] (8.2,1.0)-- (8.2,1.2);
\draw [line width=0.2pt] (8.6,0.8)-- (8.8,1.0);
\draw [line width=0.2pt] (9.0,0.8)-- (9.0,1.0);
\draw [line width=0.2pt] (9.4,0.4)-- (9.4,0.6);
\draw [line width=0.2pt] (9.4,0.6)-- (9.6,0.8);
\draw [line width=0.2pt] (9.8,0.8)-- (9.8,1.0);
\draw [line width=0.2pt] (9.6,0.8)-- (9.4,1.0);
\draw [line width=0.2pt] (9.4,0.6)-- (9.4,0.8);
\draw [line width=0.2pt] (8.8,0.4)-- (9.0,0.6);
\draw [line width=0.2pt] (9.0,0.0)-- (9.0,0.4);
\draw [line width=0.2pt] (9.0,0.4)-- (9.2,0.6);
\draw [line width=0.2pt] (9.0,0.6)-- (9.2,0.8);
\draw [line width=0.2pt] (8.8,0.4)-- (8.8,0.6);
\draw [line width=0.2pt] (9.0,1.0)-- (9.2,1.2);
\draw [line width=0.2pt] (9.0,1.0)-- (9.0,1.2);
\draw [line width=0.2pt] (8.8,1.2)-- (8.8,1.4);
\draw [line width=0.2pt] (8.2,1.2)-- (8.0,1.4);
\draw [line width=0.2pt] (8.2,1.2)-- (8.2,1.4);
\draw [line width=0.2pt] (9.8,1.0)-- (9.8,1.2);
\draw [line width=0.2pt] (10.0,1.2)-- (9.8,1.4);
\draw [line width=0.2pt] (9.6,1.2)-- (9.6,1.4);
\draw [line width=0.2pt] (9.2,1.2)-- (9.2,1.4);
\draw [line width=0.2pt] (8.8,1.4)-- (9.0,1.6);
\draw [line width=0.2pt] (8.4,1.2)-- (8.4,1.4);
\draw [line width=0.2pt] (8.4,1.4)-- (8.6,1.6);
\draw [line width=0.2pt] (8.4,1.4)-- (8.4,1.6);
\draw [line width=0.2pt] (8.0,1.4)-- (8.0,1.6);
\draw [line width=0.2pt] (8.2,2.0)-- (8.0,2.2);
\draw [line width=0.2pt] (8.2,2.0)-- (8.2,2.2);
\draw [line width=0.2pt] (8.4,2.0)-- (8.4,2.2);
\draw [line width=0.2pt] (8.8,2.0)-- (8.8,2.2);
\draw [line width=0.2pt] (9.2,2.0)-- (9.0,2.2);
\draw [line width=0.2pt] (9.2,2.0)-- (9.2,2.2);
\draw [line width=0.2pt] (9.4,1.8)-- (9.4,2.0);
\draw [line width=0.2pt] (9.4,2.0)-- (9.4,2.2);
\draw [line width=0.2pt] (9.4,2.0)-- (9.8,2.2);
\draw [line width=0.2pt] (10.2,1.8)-- (10.0,2.0);
\draw [line width=0.2pt] (10.2,2.0)-- (10.4,2.2);
\draw [line width=0.2pt] (10.6,1.8)-- (10.8,2.2);
\draw [line width=0.2pt] (10.8,2.2)-- (10.6,2.6);
\draw [line width=0.2pt] (10.8,2.2)-- (11.0,2.4);
\draw [line width=0.2pt] (10.4,2.2)-- (10.0,2.4);
\draw [line width=0.2pt] (10.4,2.2)-- (10.4,2.4);
\draw [line width=0.2pt] (9.4,2.2)-- (9.0,2.4);
\draw [line width=0.2pt] (9.4,2.2)-- (9.6,2.4);
\draw [line width=0.2pt] (8.6,2.2)-- (8.8,2.4);
\draw [line width=0.2pt] (8.8,2.4)-- (8.6,2.8);
\draw [line width=0.2pt] (9.0,2.4)-- (9.2,2.6);
\draw [line width=0.2pt] (7.6,2.0)-- (7.4,2.2);
\draw [line width=0.2pt] (7.6,2.0)-- (7.6,2.2);
\draw [line width=0.2pt] (7.4,2.2)-- (7.2,2.4);
\draw [line width=0.2pt] (7.6,2.2)-- (7.4,2.4);
\draw [line width=0.2pt] (7.2,2.0)-- (7.2,2.2);
\draw [line width=0.2pt] (7.0,2.2)-- (6.8,2.6);
\draw [line width=0.2pt] (7.0,2.2)-- (7.2,2.6);
\draw [line width=0.2pt] (7.2,2.6)-- (7.4,2.8);
\draw [line width=0.2pt] (7.2,2.6)-- (7.0,2.8);
\draw [line width=0.2pt] (6.8,2.6)-- (6.6,2.8);
\draw [line width=0.2pt] (6.8,2.6)-- (6.8,2.8);
\draw [line width=0.2pt] (6.6,2.2)-- (6.4,2.4);
\draw [line width=0.2pt] (6.8,2.2)-- (6.6,2.4);
\draw [line width=0.2pt] (10.0,2.0)-- (10.0,2.2);
\draw [line width=0.2pt] (10.6,1.8)-- (10.8,2.0);
\draw [line width=0.2pt] (10.6,1.8)-- (10.4,2.0);
\draw [line width=0.2pt] (11.0,1.8)-- (11.0,2.0);
\draw [line width=0.2pt] (11.0,1.8)-- (11.2,2.0);
\draw [line width=0.2pt] (11.2,1.6)-- (11.6,2.0);
\draw [line width=0.2pt] (11.6,2.0)-- (11.4,2.2);
\draw [line width=0.2pt] (11.6,2.0)-- (11.8,2.2);
\draw [line width=0.2pt] (11.2,2.0)-- (11.2,2.2);
\draw [line width=0.2pt] (11.0,2.0)-- (11.0,2.2);
\draw [line width=0.2pt] (11.4,2.2)-- (11.2,2.4);
\draw [line width=0.2pt] (11.4,2.2)-- (11.6,2.4);
\draw [line width=0.2pt] (9.6,2.4)-- (9.4,2.8);
\draw [line width=0.2pt] (9.6,2.4)-- (9.8,2.8);
\draw [line width=0.2pt] (10.4,2.4)-- (10.2,2.8);
\draw [line width=0.2pt] (10.4,2.4)-- (10.6,2.8);
\draw [line width=0.2pt] (11.0,2.4)-- (10.8,2.6);
\draw [line width=0.2pt] (11.0,2.4)-- (11.0,2.8);
\draw [line width=0.2pt] (9.0,2.4)-- (8.8,2.8);
\draw [line width=0.2pt] (8.4,2.6)-- (8.0,3.0);
\draw [line width=0.2pt] (8.0,2.6)-- (7.8,2.8);
\draw [line width=0.2pt] (7.6,2.6)-- (7.8,3.0);
\draw [line width=0.2pt] (7.6,2.6)-- (7.4,3.0);
\draw [line width=0.2pt] (7.0,2.8)-- (7.0,3.2);
\draw [line width=0.2pt] (7.0,2.8)-- (6.6,3.2);
\draw [line width=0.2pt] (9.2,2.6)-- (9.0,3.0);
\draw [line width=0.2pt] (9.2,2.6)-- (9.2,3.0);
\draw [line width=0.2pt] (8.6,2.8)-- (8.4,3.0);
\draw [line width=0.2pt] (8.6,2.8)-- (8.8,3.0);
\draw [line width=0.2pt] (10.0,2.4)-- (9.8,2.6);
\draw [line width=0.2pt] (10.0,2.4)-- (10.0,2.8);
\draw [line width=0.2pt] (11.2,2.4)-- (11.6,2.8);
\draw [line width=0.2pt] (11.2,2.4)-- (11.2,2.8);
\draw [line width=0.2pt] (11.6,2.4)-- (11.8,2.6);
\draw [line width=0.2pt] (11.6,2.4)-- (11.6,2.6);
\draw [line width=0.2pt] (11.2,2.8)-- (11.0,3.0);
\draw [line width=0.2pt] (11.0,2.8)-- (10.8,3.0);
\draw [line width=0.2pt] (10.2,2.8)-- (10.4,3.2);
\draw [line width=0.2pt] (10.6,2.8)-- (10.4,3.0);
\draw [line width=0.2pt] (10.6,2.8)-- (10.6,3.0);
\draw [line width=0.2pt] (9.4,2.8)-- (9.6,3.2);
\draw [line width=0.2pt] (9.4,2.8)-- (9.4,3.2);
\draw [line width=0.2pt] (9.0,3.0)-- (8.8,3.2);
\draw [line width=0.2pt] (8.8,3.0)-- (8.6,3.2);
\draw [line width=0.2pt] (8.0,3.0)-- (8.2,3.2);
\draw [line width=0.2pt] (7.8,3.0)-- (7.6,3.2);
\draw [line width=0.2pt] (7.8,3.0)-- (8.0,3.2);
\draw [line width=0.2pt] (7.4,3.0)-- (7.2,3.2);
\draw [line width=0.2pt] (8.4,3.0)-- (8.4,3.2);
\draw [line width=0.2pt] (9.2,3.0)-- (9.0,3.2);
\draw [line width=0.2pt] (9.2,3.0)-- (9.2,3.2);
\draw [line width=0.2pt] (9.8,2.8)-- (9.8,3.0);
\draw [line width=0.2pt] (9.8,3.0)-- (10.2,3.2);
\draw [line width=0.2pt] (9.8,3.0)-- (9.8,3.2);
\draw [line width=0.2pt] (10.6,3.0)-- (10.8,3.2);
\draw [line width=0.2pt] (11.0,3.0)-- (11.2,3.2);
\draw [line width=0.2pt] (11.6,2.8)-- (11.4,3.0);
\draw [line width=0.2pt] (11.6,2.8)-- (11.8,3.0);
\draw [line width=0.2pt] (11.4,3.0)-- (11.6,3.2);
\draw [line width=0.2pt] (11.8,3.0)-- (12.0,3.2);
\draw [line width=0.2pt] (11.8,2.6)-- (12.0,2.8);
\draw [line width=0.2pt] (12.0,2.8)-- (12.0,3.0);
\draw [line width=0.2pt] (12.0,2.8)-- (12.2,3.2);
\draw [line width=0.2pt] (6.4,2.4)-- (6.2,2.6);
\draw [line width=0.2pt] (6.4,2.4)-- (6.4,2.8);
\draw [line width=0.2pt] (6.2,2.6)-- (6.2,2.8);
\draw [line width=0.2pt] (6.4,2.8)-- (6.2,3.0);
\draw [line width=0.2pt] (6.2,3.0)-- (6.4,3.2);
\draw [line width=0.2pt] (6.2,3.0)-- (6.2,3.2);
\draw [line width=0.2pt] (6.6,2.8)-- (6.4,3.0);
\draw [line width=0.2pt] (8.0,2.2)-- (8.0,2.4);
\draw [line width=0.2pt] (8.4,2.2)-- (8.2,2.4);
\draw [line width=0.2pt] (8.4,2.2)-- (8.4,2.4);
\draw [line width=0.2pt] (9.4,1.4)-- (9.6,1.6);
\draw [line width=0.2pt] (9.6,1.6)-- (9.6,1.8);
\draw [line width=0.2pt] (9.6,1.6)-- (9.8,1.8);
\draw [line width=0.2pt] (9.2,1.6)-- (9.2,1.8);
\draw [line width=0.2pt] (8.6,0.8)-- (8.6,1.0);
\draw [line width=0.2pt] (8.6,1.0)-- (8.6,1.2);
\draw [line width=0.2pt] (10.0,1.2)-- (10.0,1.4);
\draw [line width=0.2pt] (9.8,0.8)-- (10.0,1.0);
\draw [line width=0.2pt] (10.0,1.0)-- (10.2,1.2);
\draw [line width=0.2pt] (10.4,1.2)-- (10.4,1.4);
\draw [line width=0.2pt] (10.4,1.6)-- (10.4,1.8);
\draw [line width=0.2pt] (7.8,1.4)-- (7.8,1.6);
\draw [line width=0.2pt] (7.8,1.6)-- (7.6,1.8);
\draw [line width=0.2pt] (7.8,1.6)-- (7.8,1.8);
\draw [line width=0.2pt] (7.8,1.8)-- (7.8,2.0);
\draw [line width=0.2pt] (7.6,2.2)-- (7.6,2.4);
\draw [line width=0.2pt] (7.4,2.4)-- (7.4,2.6);
\draw [line width=0.2pt] (6.8,2.2)-- (6.8,2.4);
\draw [line width=0.2pt] (9.4,2.2)-- (9.2,2.4);
\draw [line width=0.2pt] (9.4,2.2)-- (9.4,2.4);
\draw [line width=0.2pt] (9.8,2.2)-- (9.8,2.4);
\draw [line width=0.2pt] (10.0,2.8)-- (10.0,3.0);
\draw [line width=0.2pt] (10.0,2.8)-- (10.2,3.0);
\draw [line width=0.2pt] (10.8,2.6)-- (10.8,2.8);
\draw [line width=0.2pt] (9.4,1.0)-- (9.4,1.2);
\draw [line width=0.2pt] (9.6,0.8)-- (9.6,1.0);
\draw [dash pattern=on 1pt off 1pt,color=blue] (7.0000000000000036,2.1999999999999855)-- (7.200000000000001,1.9999999999999982)-- (7.59999999999999,2.200000000000013)-- (7.8,2.2)-- (8.200000000000001,1.9999999999999956)-- (8.400000000000007,1.9999999999999916)-- (8.399999999999984,2.600000000000019)-- (8.8,2.4000000000000066)-- (10.8,2.200000000000001)-- (11.59999999999998,2.3999999999999813);
\draw [line width=0.2pt,dash pattern=on 1pt off 1pt,color=blue] (11.59999999999998,2.3999999999999813)-- (12.1,2.4);
\draw [line width=0.2pt,dash pattern=on 1pt off 1pt,color=blue] (7.0000000000000036,2.1999999999999855)-- (6.4,2.2);
\draw (10.901263762226757,0.6460310609572653) node[anchor=north west] {Vertices of $\mathscr{Z}_{n_{j}}$};
\draw (10.901263762226757,0.4443601946672512) node[anchor=north west] {Vertices visited before $\overline{T}_{n_j}$};
\draw (10.921430848855758,0.24268932837723717) node[anchor=north west] {Vertices visited before $\overline{T}_{n_{j-1}}$};
\draw [line width=0.2pt,color=red,<-] (7.4,3.1)-- (7.4,3.5);
\draw [line width=0.2pt,color=red,<-] (7.1,2.3)-- (8.8,3.4);
\draw [line width=0.2pt,color=red,<-] (7.300000000000004,2.2999999999999963)-- (8.9,3.4);
\draw [line width=0.2pt,color=red,<-] (7.7,2.3)-- (9.0,3.4);
\draw [line width=0.2pt,color=red,<-] (7.9,2.3)-- (9.1,3.4);
\draw [line width=0.2pt,color=red,<-] (8.3,2.1)-- (9.3,3.4);
\draw [line width=0.2pt,color=red,<-] (8.5,2.7)-- (9.2,3.4);
\draw [line width=0.2pt,color=red,<-] (8.9,2.5)-- (9.4,3.4);
\draw [line width=0.2pt,color=red,<-] (10.7,2.3)-- (9.5,3.4);
\draw [line width=0.2pt,color=red,<-] (11.5,2.5)-- (9.6,3.4);
\draw (6.009016042829479,3.832430748339487) node[anchor=north west,color=red] {$\exists x, \overline{L}_{\overline{T}_{n_j}}(x)>Kn_j$};
\draw (8.6,3.7) node[anchor=north west,color=red] {$\#\mathscr{Z} _{n_j}>\varepsilon n_j$};
\draw [line width=0.2pt] (9.0,-0.0)-- (9.1,0.2);
\draw [line width=0.2pt] (9.1,0.2)-- (9.1,0.3);
\draw [line width=0.2pt] (9.1,0.2)-- (9.2,0.3);
\draw [line width=0.2pt] (9.2,0.3)-- (9.2,0.4);
\draw [line width=0.2pt] (9.2,0.3)-- (9.3,0.4);
\draw [line width=0.2pt] (9.1,0.3)-- (9.1,0.4);
\draw [line width=0.2pt] (8.999999999999998,0.40000000000000036)-- (9.0,0.5);
\draw [line width=0.2pt] (9.0,-0.0)-- (8.9,0.19999999999999973);
\draw [line width=0.2pt] (8.9,0.19999999999999973)-- (8.9,0.3);
\draw [line width=0.2pt] (8.9,0.3)-- (8.9,0.4);
\draw [line width=0.2pt] (9.0,0.5)-- (9.1,0.6);
\draw [line width=0.2pt] (9.2,0.6)-- (9.2,0.7);
\draw [line width=0.2pt] (9.3,0.4)-- (9.2,0.5);
\draw [line width=0.2pt] (9.3,0.4)-- (9.3,0.5);
\draw [line width=0.2pt] (8.8,0.6)-- (8.7,0.7);
\draw [line width=0.2pt] (8.8,0.6)-- (8.8,0.7);
\draw [line width=0.2pt] (8.8,0.6)-- (8.9,0.7);
\draw [line width=0.2pt] (9.1,0.2)-- (9.05,0.25);
\draw [line width=0.2pt] (9.1,0.3)-- (9.05,0.35);
\draw [line width=0.2pt] (9.1,0.3)-- (9.15,0.35);
\draw [line width=0.2pt] (9.15,0.35)-- (9.15,0.4);
\draw [line width=0.2pt] (9.15,0.4)-- (9.2,0.45);
\draw [line width=0.2pt] (9.15,0.4)-- (9.15,0.45);
\draw [line width=0.2pt] (9.15,0.4)-- (9.1,0.45);
\draw [line width=0.2pt] (8.9,0.4)-- (8.95,0.45);
\draw [line width=0.2pt] (8.9,0.4)-- (8.9,0.45);
\draw [line width=0.2pt] (8.9,0.19999999999999973)-- (8.95,0.25);
\draw [line width=0.2pt] (8.95,0.25)-- (8.95,0.3);
\draw [line width=0.2pt] (8.9,0.3)-- (8.95,0.35);
\draw [line width=0.2pt] (9.3,0.4)-- (9.35,0.45);
\draw [line width=0.2pt] (9.3,0.5)-- (9.25,0.55);
\draw [line width=0.2pt] (9.3,0.5)-- (9.3,0.55);
\draw [line width=0.2pt] (9.15,0.45)-- (9.15,0.5);
\draw [line width=0.2pt] (9.2,0.5)-- (9.2,0.55);
\draw [line width=0.2pt] (9.0,0.5)-- (9.0,0.55);
\draw [line width=0.2pt] (8.95,0.45)-- (8.95,0.5);
\draw [line width=0.2pt] (8.850000000000001,0.49999999999999956)-- (8.85,0.55);
\draw [line width=0.2pt] (8.85,0.55)-- (8.85,0.6);
\draw [line width=0.2pt] (9.049999999999999,0.6500000000000015)-- (9.0,0.7);
\draw [line width=0.2pt] (9.0,0.7)-- (9.05,0.75);
\draw [line width=0.2pt] (9.0,0.7)-- (9.0,0.75);
\draw [line width=0.2pt] (9.149999999999999,0.7500000000000017)-- (9.1,0.8);
\draw [line width=0.2pt] (9.1,0.8)-- (9.05,0.85);
\draw [line width=0.2pt] (9.1,0.8)-- (9.15,0.85);
\draw [line width=0.2pt] (9.1,0.8)-- (9.1,0.85);
\draw [line width=0.2pt] (9.2,0.7)-- (9.25,0.75);
\draw [line width=0.2pt] (9.2,0.7)-- (9.2,0.75);
\draw [line width=0.2pt] (9.15,0.5500000000000007)-- (9.15,0.6);
\draw [line width=0.2pt] (9.15,0.6)-- (9.1,0.65);
\draw [line width=0.2pt] (9.15,0.6)-- (9.15,0.65);
\draw [line width=0.2pt] (9.25,0.55)-- (9.3,0.6);
\draw [line width=0.2pt] (9.25,0.55)-- (9.25,0.6);
\draw [line width=0.2pt] (9.25,0.6)-- (9.3,0.65);
\draw [line width=0.2pt] (9.25,0.6)-- (9.25,0.65);
\draw [line width=0.2pt] (8.8,0.7)-- (8.75,0.75);
\draw [line width=0.2pt] (8.8,0.7)-- (8.8,0.75);
\draw [line width=0.2pt] (8.8,0.7)-- (8.85,0.75);
\draw [line width=0.2pt] (8.85,0.7)-- (8.85,0.6499999999999987);
\draw [line width=0.2pt] (8.749999999999998,0.6500000000000017)-- (8.75,0.7);
\draw [line width=0.2pt] (8.7,0.7)-- (8.65,0.75);
\draw [line width=0.2pt] (8.7,0.7)-- (8.7,0.75);
\draw [line width=0.2pt] (8.8,0.5)-- (8.75,0.55);
\draw [line width=0.2pt] (8.75,0.55)-- (8.75,0.6);
\draw [line width=0.2pt] (8.85,0.6)-- (8.9,0.65);
\draw [line width=0.2pt] (8.7,0.6000000000000003)-- (8.7,0.65);
\draw [line width=0.2pt] (8.9,0.7)-- (8.9,0.75);
\draw [line width=0.2pt] (8.9,0.75)-- (8.95,0.8);
\draw [line width=0.2pt] (8.9,0.75)-- (8.9,0.8);
\draw [line width=0.2pt] (8.9,0.8)-- (8.85,0.85);
\draw [line width=0.2pt] (8.9,0.8)-- (8.95,0.85);
\draw [line width=0.2pt] (8.9,0.8)-- (8.9,0.85);
\draw [line width=0.2pt] (8.8,0.75)-- (8.75,0.8);
\draw [line width=0.2pt] (8.8,0.75)-- (8.8,0.8);
\draw [line width=0.2pt] (8.8,0.8)-- (8.75,0.85);
\draw [line width=0.2pt] (8.8,0.8)-- (8.8,0.85);
\draw [line width=0.2pt] (9.400000000000002,0.3999999999999996)-- (9.45,0.5);
\draw [line width=0.2pt] (9.45,0.5)-- (9.45,0.55);
\draw [line width=0.2pt] (9.45,0.5)-- (9.5,0.55);
\draw [line width=0.2pt] (9.5,0.55)-- (9.45,0.6);
\draw [line width=0.2pt] (9.5,0.55)-- (9.55,0.6);
\draw [line width=0.2pt] (9.5,0.55)-- (9.5,0.6);
\draw [line width=0.2pt] (9.45,0.65)-- (9.45,0.7);
\draw [line width=0.2pt] (9.45,0.7)-- (9.45,0.75);
\draw [line width=0.2pt] (9.45,0.7)-- (9.5,0.75);
\draw [line width=0.2pt] (9.5,0.75)-- (9.45,0.8);
\draw [line width=0.2pt] (9.5,0.75)-- (9.5,0.8);
\draw [line width=0.2pt] (9.5,0.75)-- (9.55,0.8);
\draw [line width=0.2pt] (9.4,0.6000000000000014)-- (9.35,0.65);
\draw [line width=0.2pt] (9.35,0.65)-- (9.35,0.7);
\draw [line width=0.2pt] (9.35,0.7)-- (9.3,0.75);
\draw [line width=0.2pt] (9.35,0.7)-- (9.35,0.75);
\draw [line width=0.2pt] (9.3,0.75)-- (9.3,0.8);
\draw [line width=0.2pt] (9.3,0.75)-- (9.35,0.8);
\draw [line width=0.2pt] (9.35,0.8)-- (9.3,0.85);
\draw [line width=0.2pt] (9.35,0.8)-- (9.4,0.85);
\draw [line width=0.2pt] (9.35,0.8)-- (9.35,0.85);
\draw [line width=0.2pt] (8.950000000000001,0.549999999999999)-- (8.95,0.6);
\draw [line width=0.2pt] (8.95,0.6)-- (8.95,0.65);
\draw [line width=0.2pt] (8.95,0.6)-- (9.0,0.65);
\draw [line width=0.2pt] (9.049999999999999,0.450000000000001)-- (9.05,0.5);
\draw [line width=0.2pt] (9.25,0.45)-- (9.25,0.5);
\draw [line width=0.2pt] (9.25,0.35)-- (9.25,0.4);
\draw [dotted,color=cqcqcq] (11.95,1.05)-- (6.5,1.05);
\draw [dotted,color=cqcqcq] (11.95,2.0)-- (6.5,2.0);
\draw (12,2.0) node[anchor=west,color=red] {$\min_{x\in\mathscr{Z}_{n_j}} |x|$};
\draw (12,1.05) node[anchor=west,color=red] {$\max_{0\leq i\leq\overline{T}_{n_{j-1}}} |X_i|$};
\draw [color=red,<->] (11.8,1.9)-- (11.8,1.1);
\begin{scriptsize}
\draw [fill=blue] (7.0000000000000036,2.1999999999999855) circle (2pt);
\draw [fill=blue] (7.59999999999999,2.200000000000013) circle (2pt);
\draw [fill=blue] (7.8,2.2) circle (2pt);
\draw [fill=blue] (8.400000000000007,1.9999999999999916) circle (2pt);
\draw [fill=blue] (8.8,2.4000000000000066) circle (2pt);
\draw [fill=blue] (8.399999999999984,2.600000000000019) circle (2pt);
\draw [fill=blue] (10.8,2.200000000000001) circle (2pt);
\draw [fill=blue] (11.59999999999998,2.3999999999999813) circle (2pt);
\draw [fill=blue] (7.200000000000001,1.9999999999999982) circle (2pt);
\draw [fill=blue] (8.200000000000001,1.9999999999999956) circle (2pt);
\draw [fill=blue] (10.75,0.45) circle (2pt);
\draw [fill=red] (7.39999999999999,3.00000000000001) circle (2pt);
\end{scriptsize}
\end{tikzpicture}
\caption{The event $B_j$ (in red).}
\label{f:Bj}
\end{figure}

Recall the definition of $\T^{(n_{j-1})}$ from~\eqref{def-Tk}. Remark that $\max_{0\le i \le \overline T_{n_{j-1}}} | X_i| = \max_{x \in \T^{(n_{j-1})}} |x|$. 
Then $B_j$ is adapted to the filtration ${\cal G}_j$ defined as follows: 
 $${\cal G}_j
 :=
\sigma\Big( \overline L_{\overline T_k}(x), 1\le k \le n_j, \, x \in \T^{(n_j)} \Big), \qquad j\ge1. $$

 For any $x \in \ZZ_{n_j}$, let  \begin{equation}\label{def-etax}
 \zeta_x:= \max_{y\in \T, y \ge x}  {\overline L}_{ {\overline T}_{n_j}}(y) .
 \end{equation}

 \noindent Then for any event $A \in {\cal G}_{j-1}$, we have $$ \p \Big( B_j \cap A\Big) 
\ge
\e\Big( 1_{A \cap \{ \max_{x \in \T^{(n_{j-1})}} |x| < \min_{x \in  \ZZ_{n_j}} |x| , \, \# \ZZ_{n_j} \ge \varepsilon n_j \} } \big( 1- \prod_{x\in   \ZZ_{n_j}} 1_{\{ \zeta_x <K n_j\}} \big)\Big).
$$

Define  $ {\cal G}^*_j
 :=
\sigma\Big( \overline L_{\overline T_k}(x), 1\le k \le n_j, \, x \in \T^{(n_j)} , x \le \ZZ_{n_j} \Big)   $ which is a $\sigma$-field smaller than ${\cal G}_j$. Remark that $A \cap \{ \max_{x \in \T^{(n_{j-1})}} |x| < \min_{x \in  \ZZ_{n_j}} |x|\}$ is ${\cal G}^*_j$-measurable, as well as $ \ZZ_{n_j}$. On the other hand, the process $( \T^{(n_j)}, ( \overline L_{\overline T_k}(x), 1\le k \le n_j )_{x \in \T^{(n_j)}})$ is again a multi-type Galton-Watson tree; then we may apply the branching property to see that conditioned on ${\cal G}^*_j$,  $(\zeta_x)_{x \in \ZZ_{n_j}}$ are i.i.d.\ and are distributed as $\max_{z \in \T} \overline L_{\overline T_1}(z)$, see e.g. Jagers~\cite{jagers89}, Theorem 3.1 for the justification of the use of branching property along an optional line. It follows that \begin{eqnarray*}
\p(B_j \cap A)
&\ge &
\e\Big( 1_{A \cap \{ \max_{x \in \T^{(n_{j-1})}} |x| < \min_{x \in  \ZZ_{n_j}} |x| , \, \# \ZZ_{n_j} \ge \varepsilon n_j \} } \big( 1- \prod_{x\in   \ZZ_{n_j}} \p( \max_{z \in \T} \overline L_{\overline T_1}(z) <K n_j)\big)\Big)
\\
&\ge& 
\big( 1-  \p( \max_{z \in \T} \overline L_{\overline T_1}(z) <K n_j) ^{\varepsilon n_j }\big)\, \p\Big(  A \cap \{ \max_{x \in \T^{(n_{j-1})}} |x| < \min_{x \in  \ZZ_{n_j}} |x| , \, \# \ZZ_{n_j} \ge \varepsilon n_j \}  \Big) .
\end{eqnarray*}

\noindent
By Theorem~\ref{t:tailmaxlocaltime}, $$1-  \p( \max_{z \in \T} \overline L_{\overline T_1}(z) <K n_j) ^{\varepsilon n_j }
\ge
1- (1- \frac{c}{K n_j})^{\varepsilon n_j } 
\ge c_{\varepsilon, K},$$

\noindent with some positive constant $c_{\varepsilon, K}$ only depending on $\varepsilon$ and $K$. Thus we have proved that \begin{equation}\label{pbj} \p\Big(B_j \,|\, {\cal G}_{j-1} \Big)
\ge
c_{\varepsilon, K}\, \p\Big( \max_{0\le i \le \overline T_{n_{j-1}}} | X_i| < \min_{x \in \ZZ_{n_j}} |x| , \, \# \ZZ_{n_j} \ge \varepsilon n_j \big| {\cal G}_{j-1} \Big),\end{equation}

\noindent where we have used again the fact that $ \max_{0\le i \le \overline T_{n_{j-1}}} | X_i|= \max_{x \in \T^{(n_{j-1})}} |x| $.

 By Fact~\ref{f:1} (combining~\eqref{maxlimsup2} with~\eqref{aslocaltime}), we easily get that  $\p^*$-a.s.\ for all large $j$, $ \max_{0\le i \le \overline T_{n_{j-1}}} | X_i| \le  n_{j-1}^{\kappa-1 +o(1)} < n_{j-1}$ for all large $j$ [here $\kappa<2$].   Notice that $n_{j-1} <  a_{n_j}$ which is in turn smaller than  $ \min_{x \in \ZZ_{n_j}} |x|$ by~\eqref{l:zzk2}.  Hence $\p^*$-a.s.\ for all large $j$, $$ \max_{0\le i \le \overline T_{n_{j-1}}} | X_i| < \min_{x \in \ZZ_{n_j}} |x| .$$
 
 To treat $\# \ZZ_{n_j} $, we apply Fact~\ref{f:loic} and see that $$
 \limsup_{j\to \infty} \frac{\# \ZZ_{n_j} }{n_j} \ge W_\infty, \qquad \p^*\mbox{-a.s.}$$

\noindent It follows that $\p$-a.s.\ on $\{ W_\infty > \varepsilon\}$, $$ \sum_{j} 1_{ \{ \max_{0\le i \le \overline T_{n_{j-1}}} | X_i| < \min_{x \in \ZZ_{n_j}} |x| , \, \# \ZZ_{n_j} \ge \varepsilon n_j \}}= \infty.$$
 
 \noindent Applying  Fact~\ref{f:xi} gives that $\p$-a.s.\ on $\{ W_\infty > \varepsilon\}$, 
 $$\sum_j \p\big( \max_{0\le i \le \overline T_{n_{j-1}}} | X_i| < \min_{x \in \ZZ_{n_j}} |x| , \, \# \ZZ_{n_j} \ge \varepsilon n_j \big| {\cal G}_{j-1} \big)= \infty,$$
 
 \noindent which, again in view of Fact~\ref{f:xi} and~\eqref{pbj}, yields that  $ \sum_j 1_{B_j} =\infty. $ It follows that $\p$-a.s.\ on $\{ W_\infty > \varepsilon\}$, $ \limsup_{n \to \infty} \frac1n\, \max_{x \in \T} \overline L_{\overline T_n}(x)  \ge K$. Recall that $\{\T=\infty\}=\{W_\infty>0\}=\bigcup_{\varepsilon>0}\{W_\infty>\varepsilon\} $. Letting $K \to \infty$, we get Proposition~\ref{p:kappa<2}. 
 \hfill$\Box$

 \medskip
 Now we are ready to give the proofs of~\eqref{t:kappa<2} and~\eqref{t:kappa<2b} in Theorem~\ref{t:main} for the case $1< \kappa <2$:
 \medskip
 
 {\noindent\it Proof of~\eqref{t:kappa<2} in Theorem~\ref{t:main}.} Remark that for any $x \in \T$, $E_\omega (L_{\overline T_1}(x))= \ee^{-U(x)}$.  By the law of large numbers we see that for any fixed $K\ge1$,  $P_\omega$-a.s., $$ 
\lim_{n \to \infty} \frac1n \max_{| x | \le K} L_{\overline T_n} (x)= \max_{|x | \le K} \ee^{- U(x)} \le \max_{x \in \T} \ee^{- U(x)}  < \infty,$$

 \noindent by Lemma~\ref{l:U}. In view of Proposition \ref{p:kappa<2} we get  that $\p^*$-a.s., $$ \limsup_{n \to \infty} \inf_{x\in {\mathbb F}(n)} | x| \ge K.$$ 
 
 \noindent The above limsup equals in fact infinity because $K$ can be chosen arbitrarily large. This together with~\eqref{favempty} imply~\eqref{t:kappa<2}. \hfill$\Box$

 \medskip
  {\noindent\it Proof of~\eqref{t:kappa<2b} in Theorem~\ref{t:main}.} 
  Let $A>1$ be a large constant. Recall the definition of $\zeta_x$ from~\eqref{def-etax}. Define for any $j\ge1$, 
 \begin{eqnarray*} 
 C_j
 &:=& \Big\{ \max_{x \in \T: x \ge \ZZ_{n_j}}  L_{\overline T_{n_j}}(x) \le \frac{n_j}2 ,\, \max_{0\le i \le \overline T_{n_{j-1}}} | X_i| < \min_{x \in \ZZ_{n_j}} |x| , \, \# \ZZ_{n_j} \le  A\, n_j \Big\} \\
 &=&
\Big\{ \max_{x \in \ZZ_{n_j}}  \zeta_x \le \frac{n_j}2 ,\, \max_{0\le i \le \overline T_{n_{j-1}}} | X_i| < \min_{x \in \ZZ_{n_j}} |x| , \, \# \ZZ_{n_j} \le  A\, n_j 
\Big\} ,
\end{eqnarray*}

\noindent where $n_j$ and  $\zeta_x$ are defined in~\eqref{defnj} and~\eqref{def-etax}  respectively, and we denote by $ x \ge \ZZ_{n_j}$ if there exists $y \in \ZZ_{n_j}$ such that $y \le x$.  The same argument as that in the proof of Proposition~\ref{p:kappa<2} (with the same ${\cal G}_j$ as there) yields that $$
\p\Big( C_j \, \big| {\cal G}_{j-1}\Big)
\ge
 \p\big( \max_{z \in \T} \overline L_{\overline T_1} (z) < \frac{n_j}2 \big) ^{ A\, n_j } \, \p\Big( \max_{0\le i \le \overline T_{n_{j-1}}} | X_i| < \min_{x \in \ZZ_{n_j}} |x| , \, \# \ZZ_{n_j} \le  A\, n_j \big| {\cal G}_{j-1}\Big).$$

\noindent Applying Theorem~\ref{t:tailmaxlocaltime} to $\max_{x \in \T}  L_{\overline T_1}(x)$ gives that $\p\big( \max_{z \in \T} \overline L_{\overline T_1}(z)  < \frac{n_j}2 \big) ^{A\, n_j } \ge c_A$ with some positive constant $c_A$.  By Fact~\ref{f:loic}, $$
 \liminf_{j\to \infty} \frac{\# \ZZ_{n_j} }{n_j} \le W_\infty, \qquad \p^*\mbox{-a.s.}$$
 
\noindent Therefore $\p$-a.s.\ on $\{ W_\infty < A\}$,  $\sum_j 1_{C_j}
=\infty,$ which implies that there are infinitely many $j\to \infty$ such that $ \max_{x \in \T: x \ge \ZZ_{n_j}}  L_{\overline T_{n_j}}(x) \le \frac{n_j}2 < L_{\overline T_{n_j}}(\varnothing) $,  a fortiori, any favorite site $x \in {\mathbb F}(\overline T_{n_j})$ must satisfy that $x < \ZZ_{n_j}$. By~\eqref{l:zzk1}, $|x| \le c_3 \log n_j $ which by~\eqref{aslocaltime}, is in turn smaller than $(\overline T_{n_j})^b$ for any contant $ b \in (0, \frac{\kappa-1}{\kappa})$ for all $j$ large enough. In view of~\eqref{favempty}, we conclude that for those $j \to \infty$, any favorite site $x \in {\mathbb F}(\overline T_{n_j})$ must satisfy that $|x| \le K_\varepsilon$ for any fixed $0< \varepsilon< \frac12$. Hence $\p$-a.s.\ on $\{ W_\infty < A \}$, $\liminf_{n\to \infty} \max_{x \in {\mathbb F}(n)} |x| \le K_\varepsilon$,  yielding~\eqref{t:kappa<2b} by letting $A\to \infty$. \hfill $\Box$

\subsection{Proof of~\eqref{t:kappa=2}: Case $\kappa = 2$}

First we remark that by~\eqref{limsupltn} and the upper bound~\eqref{l:zzk1} in Lemma~\ref{l:zzk},  for any $\varepsilon>0$, \begin{equation}\label{limsupltn2}
 \limsup_{n\to \infty} \frac1n\, 
 \,  \max_{ |x|\ge K_\varepsilon, \, x \le \ZZ_n}
  L_{\overline T_n} (x) 
  < 
  \varepsilon, \qquad \p^*\mbox{-a.s.},
  \end{equation}

\noindent  Hence for any $\delta >0$ and $\varepsilon >0$, there is some integer $n_0(\varepsilon, \delta)$ such that for all  $n\ge n_0(\varepsilon, \delta)$, we have \begin{equation}\label{limsupltn3}
\p^*\Big( \max_{ |x|\ge K_\varepsilon, \, x \le \ZZ_n}
  L_{\overline T_n} (x) \le 2 \varepsilon n\Big) 
  \ge 
  1- \delta. 
  \end{equation}

\noindent Now we recall from ~\cite{subdiff}  that when $\kappa=2$, as $j \to \infty$, \begin{equation}\label{tight2}
  \frac{L_j({\buildrel\leftarrow \over \varnothing}) }{\sqrt{j \log j}} \,  {\buildrel \mbox{\small (law)} \over \longrightarrow} \, \mbox{ a non-degenerate distribution  on $(0, \infty)$ under $\p^*$}.
   \end{equation}

 %  \noindent In fact, it was shown in~\cite{subdiff} that under $\p^*$, if we assume an additional condition of non-lattice on the distribution of $(S_1, \P)$, $L_j(\varnothing)/\sqrt{j \log j}$ converges in law towards a non-degenerate distribution on $(0, \infty)$. It is fairly easy to adapt the arguments therein to get~\eqref{tight2}.  

   Then for any $\delta>0$, there is a sufficiently small constant $c_\delta>0$ such that for all large $j \ge j_0$,  \begin{equation}\label{tight2a}
  \p^*\Big( n_1\le L_j({\buildrel\leftarrow \over \varnothing}) \le  n_2 \Big) \ge 1- \delta,
     \end{equation}

  \noindent where $n_1 \equiv n_1(j):= \lceil c_\delta \sqrt{j \log j} \rceil$ and $n_2\equiv n_2(j):= \lceil \frac1{c_\delta} \sqrt{j \log j} \rceil$. Observe that by the Markov property, for each $n \ge1$, the following identity in law holds under the annealed probability measure $\p$:
 $$
  \max_{x \ge \ZZ_n}  L_{\overline T_n}(x) 
  \law
 \max_{1\le i \le \#\ZZ_n} \widehat L^{*, i} ,$$
  
 \noindent where $\widehat L^{*, i}, i\ge1$, are i.i.d.\ copies of $ \max_{x\in \T}  L_{\overline T_1}(x)$, independent of $ \#\ZZ_n$. Choose and fix $0<\varepsilon< \frac{c_\delta^2}{4}$.  It follows from the tail estimate in Theorem~\ref{t:tailmaxlocaltime} (with $\kappa=2$) that for all large $n \ge n_0(\varepsilon)$, \begin{eqnarray*} 
 \p\Big( \max_{x \ge \ZZ_n}  L_{\overline T_n}(x) \ge \varepsilon n \Big)
 &\le&
  1- \e \Big( 1- \frac{c}{\varepsilon n \sqrt{ \log (\varepsilon n)}}\Big)^{\#\ZZ_n}
 \\
 &\le&
   \e \frac{c\,  \#\ZZ_n}{\varepsilon n \sqrt{ \log (\varepsilon n)}} 
   \\
   &\le&
   c_\varepsilon\, (\log n)^{-1/2},
   \end{eqnarray*}
   
   \noindent by applying Fact~\ref{f:loic} to get the last inequality. In particular for all large $n$, $$\p^*\Big( \max_{x \ge \ZZ_n}  L_{\overline T_n}(x) \ge \varepsilon n \Big) \le \delta.$$
   
 Let $j $ be large. On the event of the probability term in~\eqref{tight2a}, we have $ \overline T_{n_1} \le j \le \overline T_{n_2}$. If there is some $x \in {\mathbb F}(j)$ such that $|x| > K_\varepsilon$, then $L_j(x) \ge L_j(\varnothing)$, which implies that $L_{\overline T_{n_2}}(x) \ge L_{\overline T_{n_1}}(\varnothing)\ge n_1 > 2 \, \varepsilon\,  n_2$ by the choice of $\varepsilon$. It follows that \begin{eqnarray*} && \Big\{ \max_{x \in {\mathbb F}(j)} |x| > K_\varepsilon \Big\} \cap \Big\{ n_1\le L_j({\buildrel\leftarrow \over \varnothing}) \le  n_2 \Big\} \cap \Big\{ \max_{ |x|\ge K_\varepsilon, \, x \le \ZZ_{n_2}}
  L_{\overline T_{n_2}} (x) \le 2 \, \varepsilon \, n_2\Big\}
  \\
  &&
  \subset
  \Big\{ \max_{x \ge \ZZ_{n_2}}  L_{\overline T_{n_2}}(x) \ge 2 \, \varepsilon \, n_2\Big\},
  \end{eqnarray*}
   
  \noindent whose probability (under $\p^*$) is less than $\delta$. This together with~\eqref{limsupltn3} and~\eqref{tight2a} imply that $$ \p^*\Big( \max_{x \in {\mathbb F}(j)} |x| > K_\varepsilon\Big) \le 3\,  \delta,$$
  \noindent for all large $j$, which yields~\eqref{t:kappa=2} and proves that $(\sup_{x\in\mathbb{F}(n)})_{n\ge 1}$ is tight. In order to localize the tightness, we need the following lemma:

\begin{lemma}\label{lemma:excfinite}
Let $S$ be a finite subset of $\T$ disjoint from $\mathscr{M}$. Then  $\p$-almost surely for all $n$ large enough, $\mathbb{F}(n)\cap S=\emptyset$.  
\end{lemma}
\noindent {\it Proof. }
Let $x\in S$  and $y\in\mathscr{M}$. Recall that under $P_\omega$, for any $n\ge 1$, $L_{\overline{T}_n}(x)$ (resp.  $L_{\overline{T}_n}(y)$) has the law of a sum of i.i.d. random variables of law $L_{\overline{T}_1}(x)$ (resp. $ L_{\overline{T}_1}(y)$). According to the strong law of large numbers, $n^{-1}L_{\overline{T}_n}(x)\substack{ \\ \longrightarrow \\ {n\to\infty}}  E_\omega[L_{\overline{T}_1}(x)]=\ee^{-U(x)}$ and $n^{-1}L_{\overline{T}_n}(y)\substack{ \\ \longrightarrow \\ {n\to\infty}} E_\omega[L_{\overline{T}_1}(y)]=\ee^{-U(y)}$, $P_\omega$- almost surely.

By definition of $\mathscr{M}$ and since $S\cap \mathscr{M}=\emptyset$, we have $\ee^{-U(x)}<\ee^{-U(y)}$. Therefore, $P_\omega$-almost surely there exists $n_0$ large enough such that for any $n\ge n_0$, we have $L_{\overline{T}_{n+1}}(x)<L_{\overline{T}_n}(y)$.  By the monotonicity,  for all $k\in[\overline{T}_n;\overline{T}_{n+1}]$, $L_k(x)\le L_{\overline{T}_{n+1}}(x)<L_{\overline{T}_{n}}(y)\le L_k(y)$, thus $L_k(x)<L_k(y)$ for all $k\ge \overline{T}_{n_0}$. Hence $P_\omega$-a.s.\ for all $n$ large enough $x\notin \mathbb{F}(n)$ and therefore, $S$ being finite, almost surely for all $n$ large enough $\mathbb{F}(n)\cap S=\emptyset$. \hfill$\Box$\\

   Let $0< \varepsilon < \frac12$.  Recall the definition of $K_\varepsilon$ from \eqref{Komega}.   For any $n\ge 1$, $A>0$,
\begin{equation*}
\p\big( \sup_{x\in\mathbb{F}(n)}|x|\ge K_\varepsilon \big)\le \p\big( \sup_{x\in\mathbb{F}(n)}|x|\in[K_\varepsilon, A] , K_\varepsilon < A \big)+\p\big( \sup_{x\in\mathbb{F}(n)}|x|\ge A \big).
\end{equation*}
The sequence $(\sup_{x\in\mathbb{F}(n)}|x|)_{n\ge 1}$ being tight,   there exists  $A$ large enough such that for all large $n\ge n_0$, 
\begin{equation*}
\p\big( \sup_{x\in\mathbb{F}(n)}|x|\ge K_\varepsilon \big)\le \p\big( \sup_{x\in\mathbb{F}(n)}|x|\in[K_\varepsilon, A]  , K_\varepsilon < A  \big)+\varepsilon < 2 \varepsilon,
\end{equation*}

\noindent by applying  Proposition \ref{p:main2}:  $\lim_{n\to\infty} \p\big( \sup_{x\in\mathbb{F}(n)}|x|\in[K_\varepsilon, A]  , K_\varepsilon < A  \big) =0$.  

%%$\p\big( \mathbb{F}(n)\subset\{x\in\T\; : \;  |x|<K_\varepsilon\} \big)\substack{ \\ \longrightarrow \\ {n\to\infty}} 1$.

Therefore  for all large $n$, \begin{equation}\label{keptight} \p\big( \mathbb{F}(n)\subset\{x\in\T\; : \;  |x|<K_\varepsilon\} \big) \ge 1- 2\varepsilon. 
\end{equation}

   As $\{x\in\T  :   |x| <  K_\varepsilon\} \cap \mathscr{M}^c $ is almost surely finite and disjoint from $\mathscr{M}$, according to Lemma~\ref{lemma:excfinite}, we have 
\begin{equation*}
\p\big( \mathbb{F}(n)\cap \{x\in\T  :   |x| <  K_\varepsilon\} \cap \mathscr{M}^c =\varnothing \big)\substack{ \\ \longrightarrow \\ {n\to\infty}} 1, 
\end{equation*}
which together with \eqref{keptight} yield   \eqref{t:kappa=2} and complete  the proof of Theorem~\ref{t:main}. \hfill$\Box$
  
\subsection{Proofs of Corollaries~\ref{cor:cardinal} and~\ref{cor:limsup}}

\noindent {\it Proof of Corollary~\ref{cor:cardinal}.}  At first, we shall   prove a preliminary result which will also be used in the proof of  Corollary \ref{cor:limsup}.  Let $U_{min}:=\min_{x\in\T}U(x)$. Under $P_\omega$,  the process $\big(n^{-1/2}\big(L_{\overline{T}_n}(x)-n\ee^{-U_{min}}\big)_{x\in\mathscr{M}}\big)_{n\ge 0}$ is a lattice random walk with covariance matrix equal to that of $(L_{\overline{T}_1}(x))_{x\in\mathscr{M}}$. The coefficients of this matrix are finite, as for any $x\in\mathscr{M}$, $L_{\overline{T}_1}(x)$ is stochastically smaller than a geometric random variable.  We claim that  the matrix is of rank $M:=\#\mathscr{M}$, so that  the random walk  is  a {\it genuinely} $M$-dimensional random walk.  Indeed,   suppose that there exists $(a_x)_{x\in\mathscr{M}}$ a sequence of real numbers (which may depend on the environment) not all zero such that $P_\omega$-a.s., $\sum_{x\in\mathscr{M}}a_x L_{\overline{T}_1}(x)=0$. Let $x_0$ be such that $a_{x_0}\neq 0$ and for any other $y\in\mathscr{M}$ such that $a_y\neq 0$, $|x_0|\le |y|$. As $L_{\overline{T}_1}(x_0)=\overline{L}_{\overline{T}_1}(x_0)+\sum_{\parent{z}=x_0}\overline{L}_{\overline{T}_1}(z)$ (see \eqref{siteedge}), we have that 
\begin{equation*}
\overline{L}_{\overline{T}_1}(x_0) =- \sum_{\parent{z}=x_0}\overline{L}_{\overline{T}_1}(z)+\frac{1}{a_{x_0}}\sum_{y\in\mathscr{M},y\neq x_0} a_y  \Big(\overline{L}_{\overline{T}_1}(y)+\sum_{\parent{z}=y}\overline{L}_{\overline{T}_1}(z)\Big), 
\end{equation*}
almost surely. This is absurd as conditionally on the environment the left member of this equation is a non-trivial random variable only depending  on the behaviour of the walk on the edge $(\parent{x_0},x_0)$, whereas  the right member only depends   on the behaviour of the walk on strictly distinct edges. Therefore there exists no such family $(a_x)_{x\in\mathscr{M}}$, and the covariance matrix of $(L_{\overline{T}_n}(x))_{x\in\mathscr{M}}$ for any $n$ is of rank $M$.

Now we are ready to prove  Corollary~\ref{cor:cardinal}.   Suppose that $\#\mathscr{M}\ge 4$ (otherwise the result is immediate, as $\mathbb{F}(n)\subset\mathscr{M}$ for all large $n$), and let $w,x,y,z$ be any four distinct vertices of $\mathscr{M}$. Let 
\begin{align*}
\mathcal{A}^{w,x,y,z}_n:&=\{ L_{\overline{T}_n}(w)\le L_{\overline{T}_n}(x)\le L_{\overline{T}_n}(y)\le L_{\overline{T}_n}(z)\quad\textrm{and}\quad L_{\overline{T}_{n+1}}(w)\ge L_{\overline{T}_n}(z) \}\\ 
&=\{0\le L_{\overline{T}_n}(x)-L_{\overline{T}_n}(w)\le L_{\overline{T}_n}(y)-L_{\overline{T}_n}(w)\le L_{\overline{T}_n}(z)-L_{\overline{T}_n}(w)\le L_{\overline{T}_{n+1}}(w)-L_{\overline{T}_n}(w)\}. 
\end{align*}
Now as under $P_\omega$,  the $L_{\overline{T}_{n+1}}(w)-L_{\overline{T}_n}(w)$ are i.i.d. random variables dominated by a geometric  variable, there is a $\lambda>0$ such that almost surely for all $n$ large enough, $L_{\overline{T}_{n+1}}(w)-L_{\overline{T}_n}(w)<\lambda\log(1+n)$. But as proved previously, $\big(L_{\overline{T}_n}(x)-L_{\overline{T}_n}(w),L_{\overline{T}_n}(y)-L_{\overline{T}_n}(w),L_{\overline{T}_n}(z)-L_{\overline{T}_n}(w)\big)_{n\ge 0}$ is a genuinely $3$-dimensional random walk with finite covariance matrix. Hence, equation~(2) p.~313 of~\cite{spitzer} ensures that for all $\varepsilon>0$, almost surely for all $n$ large enough this random walk is out of the ball in $\r^3$, centered at the origin and of radius $n^{1/2-\varepsilon}$. This ensures that $P_\omega\big(  \mathcal{A}_n^{w,x,y,z}\quad \mbox{i.o. as $n \to \infty$} \big)=0$. Hence almost surely, for all but a finite number of $n$, if $L_{\overline{T}_n}(w)\le L_{\overline{T}_n}(x)\le L_{\overline{T}_n}(y)\le L_{\overline{T}_n}(z)$, then for all $k\in[\overline{T}_n, \overline{T}_{n+1}]$ we have $L_k(w)\le L_{\overline{T}_{n+1}}(w) < L_{\overline{T}_n} (z)\le L_{k}(z)$. 

The vertices $w,x,y,z$ playing symmetrical roles, we have the same result for any ordering of $w,x,y,z$, and so for all but a finite number of $k$, we have that $L_k(w),L_k(x),L_k(y)$ and $L_k(z)$ are not all equal. 

The set $\mathscr{M}$ being finite, this result stands simultaneously for all quadruplets of vertices of~$\mathscr{M}$: there is only a finite number of times at which four vertices of $\mathscr{M}$ have the same   local time. In particular, there is only a finite number of times $k$ at which three vertices of $\mathscr{M}$ are in $\mathbb{F}(k)$ (since vertices of $\mathbb{F}(k)$ have the same local time), and since $\mathbb{F}(k)\subset \mathscr{M}$ for all $k$ large enough ($\kappa\in (2, \infty]$), this yields \eqref{card3}.

  Now we show that $\p$-a.s., \begin{equation}\label{limsupcard}
  \limsup_{n\to\infty} \#\mathbb{F}(n) = \min (3, M),
  \end{equation}
  
  \noindent where   $M:=\#\mathscr{M}$.  In view of the \eqref{card3} and \eqref{t:kappa>2} , it suffices  to  prove that $ \limsup_{n\to\infty} \#\mathbb{F}(n)  \ge \min (3, M)$.  When $M=1$, there is nothing to prove.  Let $M\ge2$.  Let $\mathscr{M}= \{x_i, 1\le i \le M\}$.    Observe that  $(L_{\overline{T}_n}(x_i)-L_{\overline{T}_n}(x_1) )_{2\le i \le M}$ is a genuinely $(M-1)$-dimensional   centered random walk with finite covariance matrix, almost surely it  returns to  the origin if $M\in \{2, 3\}$, and to $\{0\}\times\{0\}\times \z_-^{M-3}$ if $M\ge4$,      infinitely many times  [We have used Uchiyama \cite{Uchiyama} for the case $M\ge4$].  Now, as a.s.\ $\mathbb{F}(k)\subset\mathscr{M}$ for all $k$ large enough, we deduce that a.s.\ $\mathbb{F}(\overline{T}_n) \supset \{x_i, 1\le i \le \min(3, M) \}$ for infinitely many $n$.   This proves \eqref{limsupcard} when $M \ge 2$.  \hfill $\Box$
\medskip

%We first give the proof of Corollary~\ref{cor:limsup}, which relies among others on Lemma~\ref{lemma:excfinite}. \\

\noindent {\it Proof of Corollary~\ref{cor:limsup}.}
Let $x\in\T\backslash\mathscr{M}$; according to Lemma~\ref{lemma:excfinite}, $\p$-a.s.\ for all large $n$, $\mathbb{F}(n)\cap\{x\}=\varnothing$. Hence a.s.\ $x\notin\limsup_{n\to\infty} \mathbb{F}(n)$ and therefore $\limsup_{n\to\infty} \mathbb{F}(n)\subset\mathscr{M}$.  

Moreover, notice that~\eqref{t:kappa>2},~\eqref{t:kappa=2} and~\eqref{t:kappa<2b} ensure that almost surely, there exist $A\ge 1$ and an increasing  sequence $(\phi_n)_{n\ge 1}$ of   integers  such that for any $n\ge 1$, there exists a vertex $x\in\mathbb{F}(\phi_n)$ such that $|x|<A$.  The set $\{x\in\T\; : \;|x|<A\}$ being finite,  a fortiori there exists an $|x'|<A$ such that $x'\in\mathbb{F}(\phi_n)$ infinitely many times, i.e.\ $x'\in\limsup_{n\to\infty} \mathbb{F}(n)$; Moreover $x'$ must be a site in $\mathscr{M}$ as   $\limsup_{n\to\infty}\mathbb{F}(n)\subset \mathscr{M}$.

Let us now consider  the case $\kappa \in [ 2, \infty]$ and prove that  $\limsup_{n\to\infty}\mathbb{F}(n)=\mathscr{M}$. If $\kappa>2$, then according to~\eqref{t:kappa>2}, almost surely for all $n$ large enough we have $\mathbb{F}(n)\subset\mathscr{M}$. If $\kappa=2$, then~\eqref{t:kappa=2} together with Borel-Cantelli's  lemma ensure that there exists a {\it deterministic}  increasing sequence $(t_n)_{n\ge 1}$  of   integers such that  $\p$-a.s.\ for all $n$ large enough, $\mathbb{F}(t_n)\in\mathscr{M}$.

Now recalling from the proof of Corollary~\ref{cor:cardinal} that  under $P_\omega$,  the process $\big(\Xi_n:= n^{-1/2}\big(L_{\overline{T}_n}(x)-n\ee^{-U_{min}}\big)_{x\in\mathscr{M}}\big)_{n\ge 0}$ is  a genuinely $M$-dimensional random walk, where  $U_{min}:=\min_{x\in\T}U(x)$. Hence for  any set of the form $\{(k_1,\dots,k_{M})\in\mathbb{Z}^{M}  :    k_1<0,\dots,k_{i-1}<0,k_i>0,k_{i+1}<0,\dots k_M <0\}$  (for any fixed $1\le i \le k$) will be   recurrent for  $\Xi_{t_n}$ [in fact, take for example $i=M$ and let  $A_n:= \{ \Xi_{t_n} \in \{(k_1,\dots,k_{M})\in\mathbb{Z}^{M}  :    k_1<0,\dots,k_{M-1}<0,k_M>0 \}$.  By the central limit theorem, there exists some positive constant $c$ such that $P_\omega(A_n) \to  c$ as  $n \to \infty $, hence $P_\omega(\limsup_{n\to \infty}A_n)\ge c$. This probability is in fact equal to $1$ thanks to the    Hewitt-Savage zero-one law]; consequently,  $P_\omega$-a.s.\ for any $x\in\mathscr{M}$ there exist infinitely many $n\ge 1$ such that $L_{\overline{T}_{\phi_n}}(x)=\max_{y\in\mathscr{M}} L_{\overline{T}_{\phi_n}}(y)$, and so $x\in\limsup_{n\to\infty} \mathbb{F}(\phi_n)\subset \limsup_{n\to\infty}\mathbb{F}(n)$.

This completes  the proof of Corollary~\ref{cor:limsup}. \hfill $\Box$\\
  
 \appendix
 \section{Appendix}

 {\noindent\it Proof of~\eqref{pfi}.}  Let $i$ be a large integer. By the definition of the transition matrix $P$ in~\eqref{pij}, we see that $$
 Pf(i)
 =
\E \sum_{n=0}^\infty C_{i+n}^i \, \frac{\ee^{- n S_1}}{(1+\ee^{-S_1})^{i+n+1}} \, \frac{\Gamma(n+\gamma+1)}{\Gamma(n+1)}.
 $$
 
% \noindent Observe that for all $0< x < 1$, $$
% \sum_{n=0}^\infty C_{i+n}^i \, \frac{\Gamma(n+\gamma+1)}{\Gamma(n+1)} \, x^n 
 %=
% \frac{\Gamma(1+\gamma)}{i!} \, \empty_2F_1(1+i, 1+\gamma; 1, x),$$

% \noindent where $ \empty_2F_1(1+i, 1+\gamma; 1, x)$ is a hypergeometric function with parameter $(1+i, 1+\gamma; 1)$. By Lebedev~\cite{Lebedev}, pp.248 formula (9.5.3), $ \empty_2F_1(1+i, 1+\gamma; 1, x)= (1-x)^{-i-1-\gamma} \, \empty_2F_1(-i, -\gamma; 1, x).$ Notice that $\empty_2F_1(-i, -\gamma; 1, x)$ reduce to a polynomial, we obtain that 

 We claim that for any $\varepsilon >0$, there exist some $i_\varepsilon$, $c_\gamma>0$ such that for all $i \ge i_\varepsilon$ and $0< x< 1$, \begin{equation}\label{hypergeometric}
  \sum_{n=0}^\infty C_{i+n}^i \, \frac{\Gamma(n+\gamma+1)}{\Gamma(n+1)} \, x^n
  \le
  (1-x)^{-i-1-\gamma}\, f(i) \, \Big( (1+\varepsilon) \, x^\gamma + c_\gamma\, \varepsilon (1+x^{-1} (1-x)^{1+\gamma})\Big).
  \end{equation}

\noindent Let us admit for the moment~\eqref{hypergeometric} and finish the proof of~\eqref{pfi}. Indeed, by~\eqref{hypergeometric}, we get that 
\begin{eqnarray*}
Pf(i) 
&\le& 
f(i) \, \E \Big( (1+\varepsilon) \ee^{-\gamma S_1} + c_\gamma \varepsilon\, [(1+\ee^{-S_1})^\gamma + \ee^{S_1} ]\Big)
\\
&\le&
f(i) \, \Big( (1+\varepsilon) \E \ee^{-\gamma S_1} + c' \, \varepsilon\Big),
\end{eqnarray*}

\noindent Since $0< \gamma< \kappa-1$, $\E \ee^{-\gamma S_1}< 1$, and we can choose a sufficiently small $\varepsilon$ such that $(1+\varepsilon) \E \ee^{-\gamma S_1} + c' \, \varepsilon < 1$ and get~\eqref{pfi}.

 It remain to prove~\eqref{hypergeometric}. Firstly we  observe that for all $0< x < 1$, $$
 \sum_{n=0}^\infty C_{i+n}^i \, \frac{\Gamma(n+\gamma+1)}{\Gamma(n+1)} \, x^n 
 =
\frac{\Gamma(1+\gamma)}{i!} \frac{\d^i}{\d x^i} \Big( x^i (1-x)^{-(1+\gamma)}\Big),$$
 
 \noindent By using the change of variable: $y=1-x$, we get that for any $0< x <1$, 
 \begin{eqnarray*}
&&
 \frac{\d^i}{\d x^i} \Big( x^i (1-x)^{-(1+\gamma)}\Big)
\\&=&
 (-1)^i \frac{\d^i}{\d y^i} \Big( (1-y)^i y^{-1-\gamma}\Big) \big|_{y=1-x}
\\ &=&
 (-1)^i \frac{\d^i}{\d y^i} \Big( \sum_{k=0}^i C_i^k \, (-1)^k \, y^{k-1-\gamma}\Big)\big|_{y=1-x}
\\ &=&
 y^{-i-1-\gamma}\, \sum_{k=0}^i C_i^k \, (-1)^i\, (k-1-\gamma)\times\cdots\times(k-i-\gamma) \times (-y)^k \big|_{y=1-x}
\\&=&
i! \, y^{-i-1-\gamma}\, \sum_{k=0}^i \frac{(i+\gamma-k)\cdots(\gamma+1)}{(i-k)!} \frac{\gamma(\gamma-1)\cdots(\gamma- (k-1))}{k!} (-y)^k \big|_{y=1-x}.
\end{eqnarray*}

Recalling that  $f(i-k+1)= \frac{\Gamma(i+\gamma-k+1)}{\Gamma(i-k+1)}$ for any $0 \le k \le i$, we get that for any $0< x< 1$, \begin{eqnarray}
&& \sum_{n=0}^\infty C_{i+n}^i \, \frac{\Gamma(n+\gamma+1)}{\Gamma(n+1)} \, x^n
\nonumber\\
& =&
 (1-x)^{-i-1-\gamma}\, \sum_{k=0}^i \, f(i-k+1) \, \frac{\gamma(\gamma-1)\cdots(\gamma- (k-1))}{k!} (x-1)^k 
 \nonumber
\\
&=:&
 (1-x)^{-i-1-\gamma}\,  \Big( I_{\eqref{fgamma}} + J_{\eqref{fgamma}}\Big), \label{fgamma}
\end{eqnarray}

 \noindent where \begin{eqnarray*}
 I_{\eqref{fgamma}} &:=& \sum_{k=0}^l \, f(i-k+1) \, \frac{\gamma(\gamma-1)\cdots(\gamma- (k-1))}{k!} (x-1)^k,
 \\
 J_{\eqref{fgamma}} &:=& \sum_{k=l+1}^i \, f(i-k+1) \, \frac{\gamma(\gamma-1)\cdots(\gamma- (k-1))}{k!} (x-1)^k,
 \end{eqnarray*}
 
 \noindent with $l=l(\gamma)$ the unique integer such that $l-1 \le \gamma < l$. When $i \to \infty$, $f(i-k+1) \sim f(i)$ uniformly on $k\in \{0, 1, .., l\}$, it follows that for any $\varepsilon>0$, there exists some $i_\varepsilon>0$ such that for all $i\ge i_\varepsilon$ and all $0< x< 1$, $$
 I_{\eqref{fgamma}}= f(i) \, \sum_{k=0}^l \, \frac{\gamma(\gamma-1)\cdots(\gamma- (k-1))}{k!} (x-1)^k + \varepsilon f(i).$$

 For the term $ J_{\eqref{fgamma}}$, we discuss two cases according to the parity of $l$:
 
 (i) If $l$ is even, then every term in the sum in $J_{\eqref{fgamma}}$ is positive (for any $0< x <1$). By using the fact that $f$ is increasing, we get that \begin{eqnarray*} J_{\eqref{fgamma}} 
& \le &
 f(i) \, \sum_{k=l+1}^i  \, \frac{\gamma(\gamma-1)\cdots(\gamma- (k-1))}{k!} (x-1)^k
 \\
 &\le&
  f(i) \, \sum_{k=l+1}^\infty  \, \frac{\gamma(\gamma-1)\cdots(\gamma- (k-1))}{k!} (x-1)^k,
  \end{eqnarray*}
  
  \noindent (because every term in the above sum is positive). It follows that for any $i \ge i_\varepsilon$, and all $0< x <1$, 
  $$
  \sum_{n=0}^\infty C_{i+n}^i \, \frac{\Gamma(n+\gamma+1)}{\Gamma(n+1)} \, x^n
  \le
  (1-x)^{-i-1-\gamma}\, f(i) \, x^\gamma + \varepsilon \, (1-x)^{-i-1-\gamma}\, f(i).
  $$

 (ii) If $l$ is odd, then every term in the sum in $J_{\eqref{fgamma}}$ is negative.  For all large $i \ge i_\varepsilon$, \begin{eqnarray*}
 J_{\eqref{fgamma}} 
 &\le& \sum_{k=l+1}^{i^{1/2}} \, f(i-k+1) \, \frac{\gamma(\gamma-1)\cdots(\gamma- (k-1))}{k!} (x-1)^k
 \\
 &\le&
 (1-\varepsilon) f(i) \sum_{k=l+1}^{i^{1/2}} \, \frac{\gamma(\gamma-1)\cdots(\gamma- (k-1))}{k!} (x-1)^k,
 \end{eqnarray*}
 
 \noindent by using the fact that as $i \to \infty$, $f(i-k+1) \sim f(i)$ uniformly on $l < k\le i^{1/2}$.  Observe that for some positive constant $c_\gamma$, $\big| \frac{\gamma(\gamma-1)\cdots(\gamma- (k-1))}{k!} \big| \le c_\gamma\, i^{-\ell/2} $ for any $k > i^{1/2}$ (which is larger than $\ell$). It follows that for all $0< x <1$, $ \sum_{k> i^{1/2}} \big| \frac{\gamma(\gamma-1)\cdots(\gamma- (k-1))}{k!} (x-1)^k\big| \le  c_\gamma\, i^{-\ell/2} x^{-1} (1-x)^{i^{1/2}} \le \varepsilon \, x^{-1} (1-x)^{1+\gamma}$.  
 
 Therefore when $\ell$ is odd, we get that for all $i \ge i_\varepsilon$ and all $0<x< 1$, $$
  \sum_{n=0}^\infty C_{i+n}^i \, \frac{\Gamma(n+\gamma+1)}{\Gamma(n+1)} \, x^n
  \le
  (1-x)^{-i-1-\gamma}\, f(i) \, \Big( (1+\varepsilon) \, x^\gamma + c\, \varepsilon + \,  c_\gamma\, i^{-\ell/2} x^{-1} (1-x)^{i^{1/2}}\Big).
  $$

\noindent Therefore we get~\eqref{hypergeometric} and complete the proof of~\eqref{pfi}. \hfill$\Box$

     \medskip
  
{\noindent\bf Acknowledgments.} We are  grateful to Elie A\"\i d\'ekon for stimulating remarks and to Vladimir Vatutin for helpful  references. Special thanks go to Zhan Shi for his enthusiastic and intensive  discussions  on the topic during the whole preparation of the manuscript.

\end{document}